\documentclass[10pt,a4paper]{article}

\usepackage{amsmath,amssymb, bbm}
\usepackage{color}
\usepackage{float}
\newcommand{\R}{\mathbb{R}}

\newcommand{\Z}{\mathbb{Z}}

\newcommand{\Co}{\mathbb{C}_0}
\newcommand{\Do}{\mathbb{D}_0}
\newcommand{\Eo}{\mathbb{E}_0}
\newcommand{\Fo}{\mathbb{F}_0}
\newcommand{\Ko}{\mathbb{K}_0}
\newcommand{\bP}{\mathbb{P}}

\def\cC{{\mathcal C}}
\def\cS{{\mathcal S}}
\def\cF{{\mathcal F}}

\newcommand{\ee}{\varepsilon}
\renewcommand{\aa}{\alpha}
\newcommand{\bb}{\beta}

\renewcommand{\div}{{\rm div}\,}

\newcommand{\Sum}{\displaystyle \sum}
\newcommand{\Hs}{\dot{H^s}}

\def\d{\partial}
\def\ddj{\dot \Delta_j}
\def\ddq{\dot \Delta_q}
\def\ddqa{\dot \Delta_{q+\alpha}}
\def\ddl{\dot \Delta_l}
\def\tilde{\widetilde}
\def\hat{\widehat}
\newcommand{\D}{\Delta}

\newcommand{\n}{\nabla}

\newcommand{\Ge}{G_{ext}}
\newcommand{\Fe}{F_{ext}}

\newcommand{\om}{\omega}

\newcommand{\foe}{f_{0,\ee}}

\newcommand{\Ue}{U_\ee}
\newcommand{\ue}{u_\ee}
\newcommand{\be}{b_\ee}
\newcommand{\ve}{v_\ee}
\newcommand{\vea}{v_\ee^1}
\newcommand{\veb}{v_\ee^2}
\renewcommand{\vec}{v_\ee^3}
\newcommand{\ce}{c_\ee}
\newcommand{\de}{d_\ee}
\newcommand{\dde}{\delta_{\ee}}
\newcommand{\Uoe}{U_{0,\ee}}
\newcommand{\uoe}{u_{0,\ee}}
\newcommand{\voe}{v_{0,\ee}}
\newcommand{\boe}{b_{0,\ee}}
\newcommand{\coe}{c_{0,\ee}}

\newcommand{\De}{D_\ee}

\newcommand{\tDe}{\tilde{D}_\ee}

\newcommand{\pe}{p_\ee}
\newcommand{\tp}{\tilde{p}}
\newcommand{\tq}{\tilde{q}}

\newcommand{\tu}{\tilde{u}}
\newcommand{\tb}{\tilde{b}}

\newcommand{\tuo}{\tilde{u}_{0}}
\newcommand{\tbo}{\tilde{b}_{0}}
\newcommand{\tvo}{\tilde{v}_{0}}
\newcommand{\co}{c_0}

\newcommand{\tA}{\tilde{A}}
\newcommand{\tB}{\tilde{B}}
\newcommand{\tC}{\tilde{C}}
\newcommand{\tG}{\tilde{G}}

\newcommand{\tg}{\tilde{g}}

\newcommand{\qe}{q_{\ee}}

\newcommand{\We}{W_\ee}

\renewcommand{\Re}{R_\ee}
\newcommand{\re}{r_\ee}

\newtheorem{thm}{Theorem}[section]
\newtheorem{lem}{Lemma}[section]

\newtheorem{prop}{Proposition}[section]
\newtheorem{defi}{Definition}[section]

\newtheorem{rem}{Remark}[section]

\usepackage{graphicx}
\usepackage{tikz}
\usepackage{scalerel}
\usepackage{pict2e}
\usepackage{tkz-euclide}
\usetikzlibrary{calc}
\usetikzlibrary{patterns,arrows.meta}
\usetikzlibrary{shadows}
\usetikzlibrary{external}
\usepackage{pgfplots}
\pgfplotsset{compat=newest}
\usepgfplotslibrary{statistics}
\usepgfplotslibrary{fillbetween}
\usepackage{xcolor}
\usepackage{nicefrac}
\usepackage{pgfplots}

\newcommand{\abs}[1]{\left\vert#1\right\vert}

\newcommand{\norm}[1]{\left\Vert#1\right\Vert}
\newcommand{\pare}[1]{\left(#1\right)}
\newcommand{\eps}{\varepsilon}

\newcommand{\RR}{\mathbb{R}}

\newcommand{\TT}{\mathbb{T}}

\newcommand{\PP}{\mathbb{P}}

\newcommand{\dd}{\partial}
\newcommand{\divv}{\mbox{div\,}}

\title{Mixed 3D-2D asymptotics for the weak and strong solutions of the rotating magnetohydrodynamic system.}

\author{Fr\'ed\'eric Charve\footnote{Univ Paris Est Creteil, Univ Gustave Eiffel, CNRS, LAMA UMR8050, F-94010 Creteil, France. E-mail: frederic.charve@u-pec.fr}, Van-Sang Ngo\footnote{Univ Rouen Normandie, CNRS, Normandie Univ, LMRS UMR 6085, F-76000 Rouen, France. E-mail: van-sang.ngo@univ-rouen.fr}}

\date{}
\begin{document}

\maketitle

\begin{abstract} In this article, we consider the 3D-rotating magnetohydrodynamic (MHD) system when the initial velocity and magnetic field both feature some 2D-part (i.-e. depending only on the horizontal space variables).

We prove for weak and strong, solutions, that the limit system, when the Rossby number goes to zero (i.-e. for strong rotation), is a 2D-MHD system \emph{with three components}. Moreover we are able to provide explicit global-in-time convergence rates thanks to adapted Strichartz estimates and with the help of an additional 3D magnetic field transported by the 2D limit velocity.
\end{abstract}
\textbf{MSC: } 35B40, 35Q35, 76D03, 76U05, 76W05\\
\textbf{Keywords: }Magnetohydrodynamics, rotating fluids, Strichartz estimates, Besov and Sobolev spaces.

\section{Introduction}

\subsection{The magnetohydrodynamic system}

The magnetohydrodynamic (MHD) systems describe the evolution of the flow of charged particles under the influence of the magnetic field induced by its own motion. We can find examples of these flows in the form of plasmas in nuclear fusion reactors or in the form of charged liquid metals in the core of the Earth, which are responsible for the Earth's magnetic field. These systems show the different ways to combine the Navier-Stokes equations, which describes the behavior of a viscous incompressible fluid, and the Maxwell equations, which form the the foundation of classical electromagnetism together with the Lorentz force law. We refer the reader to the short survey \cite{A2019} and the full monograph \cite{ASR2019} for the diverse couplings of these two systems of partial differential equations and for the rigorous derivation of MHD systems from the Vlasov-Maxwell-Boltzmann system.

Let us briefly recall the derivation of the MHD system which is considered in this paper. The evolution of a viscous fluid is governed by the classical Navier-Stokes equations
\begin{equation}
	\label{eq:NS3D}
	\left\{
	\begin{aligned}
		&\dd_t U + U\cdot \nabla U - \nu \Delta U = -\nabla P + F,
		\\
		&\divv U = 0,
	\end{aligned}
	\right. 
\end{equation}
where $\nu > 0$ is the viscosity of the fluid, the velocity $U(t,x)$, the pressure $p(t,x)$ are unknown, the forcing term $F(t,x)$ is given and where $t \in \RR_+$ and $x \in \RR^3$ stand for the time and the space variables. We also suppose that the fluid evolves in the whole space $\RR^3$. For a conducting fluid, $F$ represents the Lorentz force produced by the movements of its charged particles
\begin{align*}
	F = nE + j \wedge B,
\end{align*}
where $n(t,x)$ is the electric charge density, $j(t,x)$ the electric current and where $E$ and $B$ are respectively the electric and the magnetic fields, which are determined by the quasi-static Maxwell's equations
\begin{align*}
	\left\{
	\begin{aligned}
		\dd_t B \ + \ &\nabla \wedge E = 0, && \divv B = 0,
		\\
		&\nabla \wedge B = j, && \divv E = 0.
	\end{aligned}
	\right.
\end{align*} 
Note that, in this quasi-static approximation of classical Maxwell's equations, the displacement current density $\dd_t E$ is neglected, which is reasonable in physics (see \cite{D2001}). Together with this consideration, Amp\`ere's law $j = \nabla \wedge B$ implies that $\divv j = 0$ and that the electric charge density $n$ is constant for all time $t$. Therefore, for the sake of simplicity, we suppose that $n=0$. Finally, we recall that the current density $j$ is also related to the electric and the magnetic fields through Ohm's law
\begin{align*}
	j = \sigma \pare{E + U \wedge B},
\end{align*}
where $\sigma > 0$ denotes the electrical conductivity.

Putting all the these considerations together, we obtain the following MHD equations
\begin{equation*}
	\left\{
	\begin{aligned}
		&\dd_t U - \nu \Delta U + U\cdot\nabla U = - \nabla P + (\nabla \wedge B) \wedge B,
		\\
		&\dd_t B - \frac{1}{\sigma} \Delta B = \nabla \wedge (U \wedge B),
		\\
		&\divv U = 0, \ \divv B = 0.
	\end{aligned}
	\right.
\end{equation*}
Now, classical identities show that
\begin{align*}
	(\nabla \wedge B) \wedge B = B \cdot \nabla B - \frac12 \nabla \abs{B}^2,
\end{align*}
and 
\begin{align*}
	\nabla \wedge (U \wedge B) = B \cdot \nabla U - U \cdot \nabla B + (\divv B) U - (\divv U) B.
\end{align*}
Thus, taking into account the fact that $\divv U = 0$ and $\divv B = 0$ and setting $\nu' = \frac{1}{\sigma}$, we can rewrite the MHD system as
\begin{equation}
	\label{eq:MHD}
	\tag{MHD}
	\left\{
	\begin{aligned}
		&\dd_t U - \nu \Delta U + U\cdot\nabla U - B \cdot \nabla B = - \nabla \tilde{P},
		\\
		&\dd_t B - \nu' \Delta B + U\cdot \nabla B - B\cdot\nabla U = 0,
		\\
		&\divv U = 0, \ \divv B = 0.
	\end{aligned}
	\right.
\end{equation}

\begin{rem} \sl{The condition $\div B=0$ seems to make System \eqref{eq:MHD} over-determined. However, thanks to the incompressibility condition $\divv U = 0$, we remark that $\divv B$ satisfies a linear equation of the type
		\begin{equation}
			\d_t f - \nu' \Delta f+ U\cdot \n f=0,
		\end{equation}
		and so, if $\divv B(0)=0$ then the condition $\div B(t)=0$ automatically holds true for any time $t$.}
	\label{Overdet}
\end{rem}

There is a vast literature devoted to the study of various models taking into account the MHD coupling, especially to the study of the classical MHD system \eqref{eq:MHD}. Remark that like the Navier-Stokes system \eqref{eq:NS3D}, \eqref{eq:MHD} is also parabolic and due to a good symmetry of the nonlinear terms, one can easily prove the same type of energy inequality as for \eqref{eq:NS3D}, say
\begin{equation*}
	\frac12 \pare{\norm{U(t)}^2_{L^2} + \norm{B(t)}^2_{L^2}} + \int_0^t \pare{\nu \norm{\nabla U(s)}^2_{L^2} + \frac{1}{\sigma} \norm{\nabla B(s)}^2_{L^2}} ds \leq \frac12 \pare{\norm{U(0)}^2_{L^2} + \norm{B(0)}^2_{L^2}},
\end{equation*}
which enables the classical construction of Leray weak solutions. In \cite{DL1972}, Duvaut and Lions proved the global existence of weak solutions to \eqref{eq:MHD} if the initial data $U_0,B_0 \in L^2(\RR^3)$ and the local existence of a unique strong solution to \eqref{eq:MHD} if $U_0,B_0 \in H^m(\RR^3)$, with $m \geq 3$. In \cite{ST1983}, the regularity of weak and strong solutions of \eqref{eq:MHD} were studied by Sermange and Temam in $H^1(\RR^3)$. For more results concerning the regularity of solutions of \eqref{eq:MHD}, we refer to \cite{AP2008, CW2010, CW2002, HX2005, HW2010, LZ2014}.

In the inviscid case where there are no diffusion ($\nu = \nu' = 0$), the theory for hyperbolic systems applies to \eqref{eq:MHD} and we can obtain similar results as in the case of Euler equations (System \eqref{eq:NS3D} with $\nu = 0$). We refer to \cite{CKS1997, W1997} for more details concerning the inviscid case and the inviscid limit of the MHD system. Contrary to the case where there are full diffusion in both equations of $U$ and of $B$ or the case where there are no diffusion at all, the case where there are only partial diffusion in one or both equations is a challenging problem in both physical and mathematical points of view (see for instance \cite{CW2011, CWY2014, FMRR2014, LZ2014}).

\subsection{The rotating MHD system}

In this paper, we consider the case of MHD flows in fast rotation. Such models are interesting to describe the evolution of the plasmas in the core of the Earth or in the core of stars (see for instance \cite{DDG1999} and the references therein for physical motivations). Recall that the evolution of a vector field $V$ in an absolute reference frame and a rotating frame attached to the Earth (of angular velocity $\Omega$) is given by the formula
\begin{align*}
	\pare{\partial_t V}_A = \pare{\partial_t V}_R + \Omega \wedge V.
\end{align*}
Direct calculations show that
\begin{align*}
	U_A = U_R + \Omega \wedge r,
\end{align*}
where $r$ is the position of a fluid particle and $U_A$, $U_R$ its velocity written respectively in the two references. The acceleration is then given by
\begin{align*}
	\pare{\dd_t U_A}_A = \pare{\dd_t U_R}_R + 2\Omega \wedge U_R + \Omega \wedge (\Omega \wedge r),
\end{align*}
where $2\Omega \wedge U_R$ describes the Coriolis acceleration and $\Omega \wedge (\Omega \wedge r) = -\frac12 \nabla \abs{\Omega \wedge r}^2$ describes the centrifugal acceleration, which can be combined with the pressure. Let $\eps = \frac{1}{2\abs{\Omega}}$ the Rossby number and suppose that $\Omega = \abs{\Omega} e_3$, where $e_3 = (0,0,1)$. By replacing $U_R$ by $u_\eps$, $B$ by $b_\eps$, we come to the following rotating MHD system
\begin{equation}
	\left\{
	\begin{aligned}
		&\d_t \ue -\nu \D \ue +\frac{1}{\ee} \ue\wedge e_3 +\ue\cdot \n \ue -\be\cdot \n \be=-\n \pe,
		\\
		&\d_t \be -\nu' \D \be +\ue\cdot \n \be -\be\cdot \n \ue=0,
		\\
		&\div \ue=0,
		\\
		&{(\ue,\be)}_{|t=0}=(\uoe, \boe).
	\end{aligned}
	\label{MHD}
	\tag{$MHD_\ee$}
	\right.
\end{equation}
The unknowns are $\Ue =(\ue, \be)=(\ue^1, \ue^2, \ue^3, \be^1, \be^2, \be^3)$, where we recall that $\ue$ denotes the velocity of the fluid, $\be$ the magnetic field, and $\pe$ the pressure, which gathers the classical pressure term and the centrifugal force. The diffusion coefficients $\nu,\nu'$ are positive ($\nu$ is the kinematic viscosity).

Let us recall previous results focusing on the asymptotics when the Rossby number goes to zero. In the case where there is no magnetic field ($b_\eps \equiv 0$), \eqref{MHD} reduces to the classical rotating fluid system (also called Navier-Stokes-Coriolis system)
\begin{equation*}
	\left\{
	\begin{aligned}
		&\d_t \ue - \nu \Delta \ue + \frac{1}{\ee} \ue\wedge e_3 + \ue\cdot \n \ue = -\n \pe,
		\\
		&\div \ue=0,
		\\
		&{\ue}_{|t=0}=\uoe.
	\end{aligned}
	\label{NSCe}
	\tag{$NSC_\ee$}
	\right.
\end{equation*}
This system has been intensively studied over the last thirty years so we will only briefly recall some of the most significant results in our opinion. For a complete survey on the Navier-Stokes and the Navier-Stokes-Coriolis systems, we point the reader to the book \cite{LR2002, CDGGbook}. We remark that the skew-symmetric operator $\cdot \wedge e_3$ does not contribute to classical energy estimates and so one can prove the global existence of Leray weak solutions and the local existence (global for small data) and uniqueness of Fujita-Kato strong solutions in a similar way as for the 3D Navier-Stokes system.

Let us remark that mathematically, the only way to deal with the singular perturbation $\frac{1}{\ee} \ue\wedge e_3$ is to balance it with the pressure. This formally implies that at the limit $\eps \to 0$, the pressure is independent of the vertical variable $x_3$, which in turn implies that the fluid velocity is independent of $x_3$, taking into account that the fluid is incompressible. The fluid has tendency to move in ``columns'', as predicted by what geophysicists call the Taylor-Proudman theorem. In the case of periodic domains, this singular perturbation was first studied by Babin, Mahalov and Nicolaenko \cite{BMN1} and Grenier \cite{G1997} and then also by Gallagher \cite{IG1} using Schochet's method. In these works, it is proved that the fast rotation stabilizes the system and when $\eps \to 0$, \eqref{NSCe} tends towards a 2D Navier-Stokes system with 3 components (sometimes also called 2D + $\frac12$). Moreover, when $\eps$ is small enough, the strong solution of \eqref{NSCe} globally exists in time, even for large data.

In the case of the whole space $\RR^3$, things are quite different. If we consider the wave equations associated to the Coriolis term $\frac{1}{\ee} \ue\wedge e_3$
\begin{equation*}
	\dd_t W + \frac{1}{\eps} \PP(W \wedge e_3) = 0,
\end{equation*}
where $\PP$ stands for the Leray projection from $L^2$ onto the subspace of divergence-free vector fields then, physically, this perturbation creates waves of non-vanishing amplitude and high speed $\mathcal{O}(\eps^{-1})$ which propagate through the medium (the so-called Rossby waves in physics). This waves go to infinity and carry out the energy of the system. In \cite{CDGG, CDGG2}, Chemin \textit{et al.} proved similar result as in the periodic case by using Strichartz estimates on the associated wave system, which indeed characterizes the energy dispersion previously described. Note that in accordance to the Taylor-Proudman theorem, the fluid has a 2D behavior at the limit $\eps \to 0$, and the only vector field, which is of $L^2$ finite energy in $\RR^3$ and independent of $x_3$, is zero.

Thus, in \cite{CDGG, CDGG2}, to put in evidence the 2D limit system with 3 components, the authors considered initial data of the form $2D + 3D$, which we call \textbf{non-conventional} data. Similar results of convergence are obtained through the use of Strichartz estimates for various systems: among them the primitive system and the strongly stratified Boussinesq system. We refer to the following subsection for details.

The rotating MHD system \eqref{MHD} is more difficult to deal with because the magnetic field $b_\eps$ does not \textit{a priori} possess any dispersion. 
Inspired by the work of Desjardins, Dormy and Grenier \cite{DDG1999}, by considering the case where $b_\eps$ is a small perturbation of $e_3$, in \cite{BIM2005}, Benameur \textit{et al.} studied the convergence of strong solutions of \eqref{MHD} when $\eps \to 0$ in both cases where the domain is $\TT^3$ or $\RR^3$ and for evanscent viscosity ($\nu=\ee$). In \cite{N2017}, Ngo considered \eqref{MHD} with no vertical diffusion and only a small horizontal diffusion, in the case where $b_\eps$ is a small perturbation of $e_3$ and proved that for $\eps$ close to $0$, \eqref{MHD} is globally well-posed for large initial data, using meticulous Strichartz estimates on the associated wave system. Also using Strichartz estimates, in \cite{AKL2021}, Ahn, Kim and Lee proved similar results as in \cite{N2017} for system \eqref{MHD} without the hypothesis that $b_\eps$ is a small perturbation of $e_3$ (see also \cite{KimJ, TY}). We also refer to the work of Cobb and Fanelli \cite{CF2021} concerning the asymptotics of the rotating MHD system in the density-dependent case. Finally, we want to mention the works of Desjardins, Dormy and Grenier \cite{DDG1999} and of Rousset \cite{R2005} on the stability of the Ekman-Hartmann boundary layers for \eqref{MHD}.

\subsection{Statement of the results}

\subsubsection{Hidden asymptotics}

As was first observed in \cite{CDGG}, systems (in the whole space $\R^3$) featuring penalized terms and subject to dispersion may converge to a incomplete limit model when appropriate initial data are not considered.

The first example is the Rotating fluids system \eqref{NSCe}:  when the initial data is a \textbf{conventional} data for the 3D-Navier-Stokes model (that is $\voe(x)=\voe(x_1,x_2,x_3)$), then the solutions converge to zero. As proved in \cite{CDGG}, richer asymptotics are reached for \textbf{non-conventional} initial data: if the initial data is of the form $\voe(x_h,x_3)+\tu_0(x_h)$ (where $x_h=(x_1,x_2)$ stands for the horizontal space variables) the solutions converge, as the Rossby number goes to zero, to the unique solution of the 2D-Navier-Stokes system with three components (and whose initial data is $\tu_0$). This is the reason why we can speak of hidden asymptotics. We refer to \cite{CDGG2, CDGGbook, FCRF, MuWei} for refined studies in this direction.

Another example is the strongly stratified Boussinesq system which only takes into account the influence of gravity (the small parameter is called the Froude number and is linked to the buoyancy effect). For a conventional initial data $\Uoe(x)=(\voe(x),\theta_{0,\ee}(x))$ the solutions converge to the unique solution of the $3D$-Navier-Stokes system with two components (see \cite{Scro3}). Surprizingly, this limit model does not depend anymore on the thermal diffusivity $\nu'$ involved in the equation on the potential temperature $\theta_{\ee}$, which  suggests that there are some hidden asymptotics. It is the object of \cite{FCStratif1, FCStratif2} to study this phenomenon: getting a complete limit system requires us to consider \textbf{non-conventional} initial data of the form $\Uoe(x)+(0,0,0,\theta_{0,\ee}(x_3))$ and the previous limit system is then completed by a 1D-heat equation on $\theta_\ee$ whose diffusion term involves $\nu'$.

Concerning the primitive system, which takes into account both the rotation of the Earth and the influence of gravity, there does not seem to be some hidden asymptotics as the limit system, called the 3D-quasigeostrophic system, genuinely depends on each parameter. In this case, any result will strongly rely on the structure of the limit system. We refer to \cite{Chemin2, FC1, FC2, FCVSN, FCRF, FCPAA, FCcompl, FCRF} for detailed studies.

We will not give much details about the analoguous works in a torus or between planes because they do not involve dispersion but different mechanisms (and precise study of resonences).\\

In the present work, we focus on System \eqref{MHD} and we wish to provide a limit system as general as possible. This will be possible when considering the following initial data:
\begin{equation}
{(\ue,\be)}_{|t=0}=(\tuo(x_h) +\voe(x), \tbo(x_h) +\co(x)).
 \label{Initdata}
\end{equation}
Assuming that $\voe$ is bounded we obtain as a limit the following 2D MHD system (depending on the horizontal variables) with three components:
\begin{equation}
\begin{cases}
\d_t \tu -\nu \D_h \tu +\tu^h\cdot \n_h \tu -\tb^h\cdot \n_h \tb=-\left(\begin{array}{c}\n_h \tq^0 \\ 0 \end{array}\right),\\
\d_t \tb -\nu' \D_h \tb +\tu^h\cdot \n_h \tb -\tb_h\cdot \n_h \tu=0,\\
\div_h \tu=0,\\
{(\tu,\tb)}_{|t=0}=(\tuo, \tbo),
\end{cases}
\label{MHD2D}
\tag{$2D-MHD^3$}
\end{equation}
completed by a second magnetic field transported by the previous velocity $\tu$:
\begin{equation}
\begin{cases}
\d_t c -\nu' \D c +\tu\cdot \n c -c\cdot \n \tu=0,\\
c_{|t=0}=\co.
\end{cases}
\label{Mag}
\tag{$3D-M$}
\end{equation}
\begin{rem}
 \sl{As in Remark \ref{Overdet}, as $\div \co=0$, then $\div c(t)=0$ for all time.}
 \label{Remdivc}
\end{rem}
Let us first introduce the spaces we will use in this article: for $s\in \R$, $d\in \{2,3\}$ and $T>0$ we define:
$$
 \dot{E}_T^s(\R^3)=\Big[\mathcal{C}_T(\Hs (\R^3)) \cap L_T^2(\dot{H}^{s+1}(\R^3))\Big]^3, \quad \dot{F}_T^s(\R^d)=\Big[\mathcal{C}_T(\Hs (\R^d)) \cap L_T^2(\dot{H}^{s+1}(\R^d))\Big]^6
$$
whose norms are defined as follows:
\begin{equation}
\begin{cases}
\|u\|_{\dot{E}_T^s}^2 \overset{def}{=}\|u\|_{L_T^\infty \Hs }^2 +\nu \int_0^T \|u(\tau)\|_{\dot{H}^{s+1}}^2 d\tau, \\
\|(u,b)\|_{\dot{F}_T^s}^2 \overset{def}{=}\|u\|_{L_T^\infty \Hs }^2 +\|b\|_{L_T^\infty \Hs }^2 +\nu \int_0^T \|u(\tau)\|_{\dot{H}^{s+1}}^2 d\tau +\nu' \int_0^T \|b(\tau)\|_{\dot{H}^{s+1}}^2 d\tau,
\end{cases}
\end{equation}
where $H^s(\R^d)$ and $\dot{H}^s(\R^d)$ respectively denote the inhomogeneous and homogeneous Sobolev spaces of index $s\in \R$. When $T=\infty$ we simply write $\dot{E}^s$ or $\dot{F}^s$ and the corresponding norms are understood as taken over $\R_+$ in time.

We can now state the main results of the present article in the following sections.

\subsubsection{Existence and convergence for Weak solutions (in the sense of Leray)}

\begin{thm} (Existence of Leray weak solutions)
 \sl{Let $\ee>0$ be fixed. For any $\tuo,\tbo \in L^2(\R^2)$, any $\voe\in L^2(\R^3)$ and any $\co \in H^\frac12(\R^3)$ (with $\div_h \tuo =\div_h \tbo =\div \voe =\div \co =0$), the respective solutions $(\tu,\tb)$ and $c$ of Systems \eqref{MHD2D} and \eqref{Mag} globally exist in $\dot{F}^0(\R^2)$ and $\dot{E}^0 \cap \dot{E}^{\frac12}(\R^3)$ and there exists a global weak solution $\Ue=(\ue,\be) \in \dot{F}^0(\R^3)$ of \eqref{MHD} with initial data as in \eqref{Initdata}. The pressure satisfies $\pe \in L^\frac83 (\R_+,L^2 (\R^3))+L^2 (\R_+,L^2 (\R^3))+L^2 (\R_+,L^2 (\R^2)) +\dot{E}^1(\R^3)$.

 Moreover there exists a constant $\Do>0$ depending on $\nu,\nu', \|\tuo\|_{L^2(\R^2)}, \|\tbo\|_{L^2(\R^2)}$ and $\|\co\|_{H^\frac12(\R^3)}$ such that if we define $\De \overset{def}{=} \Ue-(\tu,\tb+c)$, then for all $t\geq 0$,
 \begin{equation}
  \|\De(t)\|_{\dot{E}^0} \leq \Do (\|\voe\|_{L^2}^2 +1).
\label{estimaprioriLeray}
 \end{equation}
}
 \label{ThLeray}
\end{thm}
The convergence is described as follows:
\begin{thm}(Convergence)
 \sl{For any $\Co\geq 1$, any divergence-free vectorfields $\tuo,\tbo \in L^2(\R^2)$, $\co \in H^\frac12(\R^3)$ and (for all $\ee>0$) $\voe\in L^2(\R^3)$ such that:
 \begin{equation}
  \max\left(\|(\tuo, \tbo)\|_{L^2(\R^2)},\; \|\co\|_{H^\frac12(\R^3)},\; sup_{\ee>0} \|\voe\|_{L^2(\R^3)}\right) \leq \Co,
  \label{Hyp:Maj}
\end{equation}
 then any weak solution $\Ue$ of \eqref{MHD} (with initial data as in \eqref{Initdata}) such that $\De=\Ue-(\tu,\tb+c)$ satisfies \eqref{estimaprioriLeray} converges (as $\ee$ goes to zero) to $(\tu,\tb+c)$ in the following sense:
 \begin{enumerate}
  \item For any $r \in]2,6[$ there exists a constant $\Fo$ (depending on $\nu, \nu',\Co, r$) such that:
\begin{equation}
\|\ue-\tu\|_{\tilde{L}^{\max\left(2,\frac{4r}{5(r-2)}\right)} \dot{B}_{r,2}^0 +\tilde{L}^{\max\left(1,\frac{2r}{3r-5}\right)} \dot{B}_{r,2}^0} \leq \Fo \ee^{\frac1{2r} \min(r-2,3-\frac{r}2)}.
\label{ThCV-estim1}
\end{equation}
Moreover we have the simpler (but local-in-time) estimates: for any $r\in]2,6]$ there exists a similar constant $\Fo$ such that for any $t\geq 0$,
\begin{equation}
 \|\ue-\tu\|_{L_t^2 L^r} \leq \Fo \max(1,t)^{\delta_r} \ee^{m_r},
 \label{ThCV-estim2}
\end{equation}
where $\delta_r\overset{def}{=}\frac1{4r} \max(10-3r,\frac57 (6-r))$ and $m_r\overset{def}{=}\frac1{2r}\min(r-2,\frac17(6-r))$.
\item Seconds, $\be$ converges towards $\tb+c$ in the following sense: for any $r\in]2,6[$,
$$
 \|\be-\tb-c\|_{L_{loc}^2(\R_+, L_{loc}^r(\R^3))}\underset{\ee \rightarrow 0}{\longrightarrow}0.
 $$
 \end{enumerate}
}
\label{ThCV}
\end{thm}
\begin{rem}
 \sl{\begin{enumerate}
      \item We can say the previous results reach \textbf{hidden asymptotics} in the sense that when there are no 2D-parts in the initial data ($\tuo=\tvo=0$), the solution of system \eqref{MHD2Db} reduces to zero, and System \eqref{Mag} reduces to the heat equation (transport by the 2D-limit velocity $\tu$ vanish). This is the context of \cite{BIM2005, AKL2021, KimJ, TY}, see also \cite{BIM2005}.
      \item The previous results can be extended to initial data as follows:
      \begin{equation}
{(\ue,\be)}_{|t=0}=(\tuo(x_h) +\voe(x), \tbo(x_h) +\coe(x)).
 \label{Initdata2}
\end{equation}
      If $(\coe)_{\ee>0}$ is uniformly bounded, existence remains unchanged and convergence concerns $\Ue-(\tu, \tb+\ce)$ where $\ce$ solves \eqref{Mag} with initial data $\coe$.
      \item We can also consider an initial 2D-part $(\tu_{0,\ee}, \tb_{0,\ee})$ depending on $\ee$ and converging to $(\tuo, \tbo)$.
       \item With slight adaptation, we can add in the first three lines of System \eqref{MHD} an external force term of the form $\tilde{f}(x_h)+f_\ee(x)$ where $\tilde{f}\in L^1 \dot{H}^{-1}(\R^2)$ and $(f_\ee)_{\ee>0}$ uniformly bounded in $\ee$ in $L^1 \dot{H}^{-1}(\R^3)$.
      \item With another slight adaptation of our proof, we can allow $\|\voe\|_{L^2} \leq \Co \ee^{-\gamma}$ (with $\gamma$ bounded by a constant depending on the chosen $r$) and still have as convergence rate a power of $\ee$.
     \end{enumerate}
}
\label{Remweak}
\end{rem}

\subsubsection{Existence and convergence for Strong solutions (in the sense of Fujita and Kato)}

In this part, we directly consider a general initial data as in \eqref{Initdata2}, and $c$ is now replaced by $\ce$, solving:
\begin{equation}
\begin{cases}
\d_t \ce -\nu' \D \ce +\tu\cdot \n \ce -\ce\cdot \n \tu=0,\\
{\ce}_{|t=0}=\coe.
\end{cases}
\label{Mag2}
\tag{$3D-M_\ee$}
\end{equation}

\begin{thm}
 \sl{For any $\Co\geq 1$ (size), any $\delta\in]0,\frac16]$ (extra regularity), any $\gamma \in [0, \frac5{12}\delta]$,  there exists $\ee_0>0$ and $\Ko,\Fo\geq 1$ (depending on $\nu,\nu',\Co,\delta, \gamma,s$) such that for any $\ee\in]0,\ee_0]$ and any initial data as in \eqref{Initdata2} such that $\tuo, \tbo \in H^\delta(\R^2)$, $\voe, \coe \in (\dot{H}^{\frac12-\delta}(\R^3) \cap \dot{H}^{\frac12+\delta}(\R^3)) \times H^{\frac12+\delta}(\R^3)$ with (for some $\mu\in[0,1[$ as close to 1 as we wish):
 $$
  \|\voe\|_{\dot{H}^{\frac12+\mu\delta} \cap \dot{H}^{\frac12+\delta}}\leq \Co \ee^{-\gamma}, \mbox{ and } \|\coe\|_{H^{\frac12+\delta}}\leq \left(\Ko |\ln \ee|\right)^{\frac14},
 $$
then $T_\ee^*=+\infty$ and we have the following convergence rate:
$$
 \|\Ue-(\tu,\tb+\ce)\|_{L^2 L^\infty} \leq \Fo \ee^{\frac{1}{18} (\frac{\delta}2-\gamma)}.
$$
 }
\label{ThCVStrongsimplif}
\end{thm}

\begin{rem}
 \sl{\begin{enumerate}
 \item  As previously, with slight adaptations of the proof,  our result can be extended when:
 \begin{itemize}
  \item the initial 2D-part depends on $\ee$ and $(\tu_{0,\ee}, \tb_{0,\ee}) \underset{\ee \rightarrow 0}{\longrightarrow} (\tuo, \tbo)$ in $H^{\delta}(\R^2)$.
       \item the velocity part of System \eqref{MHD} features an external force term of the form $\tilde{f}(x_h)+f_\ee(x)$ where $\tilde{f}\in L^1 \dot{H}^{-1}\cap \dot{H}^{-1+\delta}(\R^2)$ and $(f_\ee)_{\ee>0}$ uniformly bounded in $\ee$ in $L^1 \dot{H}^{-\frac12}\cap \dot{H}^{-\frac12+\delta}(\R^3)$.
 \end{itemize}
 \item It is natural to compare our results with \cite{N2017, AKL2021, KimJ, TY} which focus on global existence when $\ee$ is small enough (but do not study limits) of strong solutions for small \textbf{conventional} inital data. The first article focusses on a slightly different MHD model, with anisotropic evanescent viscosities.

 Our results consider both weak and strong solutions and generalize the previous cited works as we get global existence and reach hidden limits for \textbf{non-conventional} large initial data.

 More precisely, if we reword the result from \cite{AKL2021} with our notations, without bidimensional parts in the initial data ($\tuo=\tbo=0$), for any $\ee>0$, $\delta\in]0,\frac2{17}[$ there exists $C>0$ small enough such that if $\uoe=u_0$ and $\boe=b_0$ are in $H^{\frac12+\delta}(\R^3)$ (inhomogeneous space) and satisfy the following bounds (we choose $q=\frac3{1-\frac{\delta}2}$ and the conditions on $q,s$ in Theorem 1 from \cite{AKL2021} give the bound $\frac2{17}$)
 $$
 \begin{cases}
   \|u_0\|_{\dot{H}^{\frac12+\delta}} + \|b_0\|_{\dot{H}^{\frac12+\delta}} \leq C \ee^{-\frac{\delta}2},\\
(\|u_0\|_{L^2}^2+\|u_0\|_{L^2}^2)  (\|u_0\|_{\dot{H}^{\frac12+\delta}}^{1+\frac2{\delta}} +\|b_0\|_{\dot{H}^{\frac12+\delta}}^{2(1+\frac2{\delta})}) \leq C \ee^{-\frac{1+\delta}2},
 \end{cases}
 $$
 then there exists a global mild solution. In \cite{KimJ} the previous conditions are improved as follows: for any $\delta\in]0, \frac12]$, there exists some $\aa\in]0,2\delta[$ such that for any $u_0,b_0 \in H^{\frac12+\delta}$ with
 $$
 \begin{cases}
  \|(u_0,b_0)\|_{\dot{H}^{\frac12+\delta}} \leq C \ee^{-\frac{\delta}2},\\
  \|(u_0,b_0)\|_{L^2}^\aa \|(u_0,b_0)\|_{\dot{H}^{\frac12+\delta}} \leq C \ee^{\frac{\aa}4-\frac{\delta}2},
 \end{cases}
 $$
 there exists a unique global mild solution. They also provide time-estimates.

 Finally, in \cite{TY} the authors prove (also by fixed point arguments) that for any $u_0,b_0\in \dot{H}^\frac12(\R^3)$ there exists $\ee_0$ such that for any $\ee\in]0,\ee_0]$, we have a unique global strong solution. Note that as in \cite{N2017}, in order to reach global existence, the authors introduce an auxiliary system (following the ideas from \cite{FC2, N2017, FCPAA, FCRF, FCStratif2}) that also features the heat equation on the magnetic field component (see also \cite{BIM2005}). We refer to Remark \ref{Remweak} where we explain that the hidden part of this heat equation is System \eqref{Mag2}.
\end{enumerate}
}
\end{rem}

\subsection{Outline of the article}

The present article unfolds as follows: in the next section, we explain through a formal argument how we may have an intuition about the limit system and the more general initial data (see \eqref{Initdata2}) that has to be considered to reach this limit. Then we perform a change of unknown functions, get the system on which we will work, and study the limit systems.

Sections 3 and 4 are devoted to prove existence and convergence results related to weak solutions (Theorems \ref{ThLeray} and \ref{ThCV}). In Section 5, we first state a more precise version for the existence and convergence results for strong solutions (Theorems \ref{Th0FK} and \ref{ThCVStrong}) and give their proofs.

In the appendix, we state general notations and results, and give proofs of several results used in the previous sections: estimates for external force terms, isotropic and anisotropic Strichartz estimates.

\section{The limit system}

\subsection{A candidate for the limit}

In this section we explain (in the same spirit as in \cite{FCStratif1}) how we are formally lead to consider Systems \eqref{MHD2D} and \eqref{Mag} as a limit. Let us begin with the pressure term: considering the divergence of the first three lines in \eqref{MHD}, as $\div \ue=0$ we can decompose $\pe$ and identify its truly penalized part:
\begin{equation}
\pe=\pe^0+\frac1{\ee} \pe^1 \quad \mbox{with }
 \begin{cases}
 \vspace{0.1cm}
  \D \pe^0= -\sum_{i,j=1}^3 \d_i \d_j (\ue^i \ue^j-\be^i \be^j),\\
  \D \pe^1=\d_2 \ue^1 -\d_1 \ue^2.
 \end{cases}
\end{equation}
We can now rewrite System \eqref{MHD} with well-identified penalized terms:
$$
 \begin{cases}
 \vspace{0.1cm}
  \d_t \ue^1-\nu \D \ue^1+\ue\cdot \n \ue^1 -\be\cdot \n \be^1=-\d_1 \pe^0 -\frac1{\ee}(\d_1 \pe^1+\ue^2),\\
  \vspace{0.1cm}
  \d_t \ue^2-\nu \D \ue^2+\ue\cdot \n \ue^2 -\be\cdot \n \be^2=-\d_2 \pe^0-\frac1{\ee}(\d_2 \pe^1-\ue^1),\\
  \d_t \ue^3-\nu \D \ue^3+\ue\cdot \n \ue^3 -\be\cdot \n \be^3=-\d_3 {\pe}^0-\frac1{\ee} \d_3 \pe^1,\\
  \d_t \be-\nu' \D \be+\ue\cdot \n \be -\be\cdot \n \ue=0.
 \end{cases}
$$
The announced formal argument is the following: let us assume that $(\ue,\be,\pe^0,\pe^1)\underset{\ee \rightarrow 0}{\longrightarrow} (\tu, \tb, \tp^0, \tp^1)$ strongly enough so that we have convergence of the derivatives and non-linear terms on one hand and convergence of the penalized terms on the other hand:
\begin{equation}
 \begin{cases}
  -\frac1{\ee}(\d_1 \pe^1+\ue^2) \underset{\ee \rightarrow 0}{\longrightarrow} \tA,\\
  -\frac1{\ee}(\d_2 \pe^1-\ue^1) \underset{\ee \rightarrow 0}{\longrightarrow} \tB,\\
  -\frac1{\ee} \d_3 \pe^1 \underset{\ee \rightarrow 0}{\longrightarrow} \tC.
 \end{cases}
 \label{limites}
\end{equation}
This immediately implies that the formal limits not only are solutions of the following limit system:
\begin{equation}
 \begin{cases}
  \d_t \tu-\nu \D \tu+\tu\cdot \n \tu -\tb\cdot \n \tb= -\n \tp^0 +(\tA,\tB,\tC),\\
  \d_t \tb-\nu' \D \tb+\tu\cdot \n \tb -\tb\cdot \n \tu= 0,
 \end{cases}
 \label{Systlim1}
\end{equation}
but also satisfy:
$$
 \begin{cases}
  \d_1 \tp^1+\tu^2=0,\\
  \d_2 \tp^1-\tu^1=0,\\
  \d_3 \tp^1=0,
 \end{cases}
 \quad \mbox{and }
 \begin{cases}
 \div \tu=0=\div \tb,\\
 \vspace{0.1cm}
  \D \tp^0= -\sum_{i,j=1}^3 \d_i \d_j (\tu^i \tu^j-\tb^i \tb^j),\\
  \D \tp^1=\d_2 \tu^1 -\d_1 \tu^2.
 \end{cases}
$$
It follows that $\tp^1$ and $(\tu^1, \tu^2)=(\d_2 \tp^1, -\d_1 \tp^1)$ are independant of $x_3$. As $\div \tu=0$ we get that $\tu$ is also independant of $x_3$. Moreover introducing the vorticity $\om (\tu)=\d_1 \tu^2-\d_2 \tu^1$ we get $\om (\tu)=-\D_h \tp^1= -\D \tp^1$, so that $ \tp^1=-\D_h^{-1} \om(\tu)$ and:
\begin{equation}
 \tu^h=\left(\begin{array}{c}\tu^1\\ \tu^2\end{array}\right)=\left(\begin{array}{c}-\d_2\\ \d_1\end{array}\right) \D_h^{-1} \om(\tu)= \n_h^\perp \D_h^{-1} \om(\tu).
 \label{Gradorth}
\end{equation}
Computing $\d_1 \eqref{limites}_2-\d_2 \eqref{limites}_1$ we obtain in addition that
\begin{equation}
 -\frac1{\ee}\d_3 \ue^3 \underset{\ee \rightarrow 0}{\longrightarrow} \d_1 \tB-\d_2 \tA \overset{def}{=} \tG.
\end{equation}
Next, computing $\d_1 \eqref{limites}_1+\d_2 \eqref{limites}_2+\d_3 \eqref{limites}_3$, we obtain that (when $\ee \rightarrow 0$):
$$
 \d_1 \tA+\d_2 \tB+\d_3 \tC=0.
$$
Formally solving the resulting system
$$
\begin{cases}
 \d_1 \tA +\d_2 \tB=-\d_3 \tC,\\
 -\d_2 \tA+\d_1 \tB= \tG,
\end{cases}
$$
we end-up with:
$$
\left(\begin{array}{c}\tA\\ \tB\end{array} \right)=-\n_h \D_h^{-1} \d_3 \tC +\n_h^{\perp} \D_h^{-1} \tG,
$$
which turns the velocity equation from \eqref{Systlim1} into:
\begin{equation}
 \begin{cases}
  \d_t \tu^1 -\nu \D_h \tu^1 +\tu^h\cdot \n_h \tu^1 =-\d_1 \tp^0 +\tb\cdot \n \tb^1 -\d_1 \D_h^{-1} \d_3 \tC -\d_2 \D_h^{-1} \tG,\\
  \d_t \tu^2 -\nu \D_h \tu^2 +\tu^h\cdot \n_h \tu^2 =-\d_2 \tp^0 +\tb\cdot \n \tb^2 -\d_2 \D_h^{-1} \d_3 \tC +\d_1 \D_h^{-1} \tG\\
  \d_t \tu^3 -\nu \D_h \tu^3 +\tu^h\cdot \n_h \tu^3 =-\d_3 \tp^0 +\tb\cdot \n \tb^3 +\tC.
 \end{cases}
\end{equation}
Using that $\d_3 \tu=0$, we can make two useful observations: first, in the third equation we can write
$$
\d_3(-\d_3 \tp^0 +\tb\cdot \n \tb^3 +\tC)=0.
$$
Seconds, considering $\d_3 \eqref{Systlim1}_2$, $\d_3 \tb$ obviously satisfies the following linear system:
\begin{equation}
\begin{cases}
\d_t \d_3 \tb -\nu' \D \d_3 \tb +\tu\cdot \n \d_3 \tb -\d_3 \tb\cdot \n \tu=0,\\
\d_3 \tb_{|t=0}=\d_3 \tbo.
\end{cases}
\end{equation}
As a consequence, if $\d_3 \tbo=0$, then for all time $\d_3 \tb(t)=0$ and the previous condition turns into
\begin{equation}
 \d_3(-\d_3 \tp^0 +\tC)=0.
\end{equation}
This suggests us to impose $\d_3 \tbo=0$ and $(\tC,\tG)=(\d_3 \tp^0,0)$, which turns the previous velocity equation into:
$$
\d_t \tu -\nu \D_h \tu +\tu^h\cdot \n_h \tu- \tb^h\cdot \n \tb =-\left(\begin{array}{c} \n^h(p^0+\Delta_h^{-1}\d_3^2 \tp^0) \\0\end{array}\right) =-\left(\begin{array}{c} \n^h\D_h^{-1}\D\tp^0 \\0\end{array}\right).
$$
As $\d_3 \tu=\d_3 \tb=0$, we can introduce $\tq^0$ (only depending on $x_h$) such that:
$$
\D \tp^0= -\sum_{i,j=1}^2 \d_i \d_j (\tu^i \tu^j-\tb^i \tb^j)=-\D_h \tq^0,
$$
which is exactly what we needed to obtain System \eqref{MHD2D}, bidimensional MHD system (as functions only depend on $x_h$) but with three components.

\subsection{Change of unknown function}

As in \cite{CDGGbook, FCStratif1, FCStratif2}, we first need to rewrite the previous limit model to correctly define the system we will study. Thanks to \eqref{Gradorth} the following rotation term can be written as an horizontal gradient: $\tu \wedge e_3=-\n_h \tp^1$ and as $\tu,\tb$ do not depend on $x_3$ we can rewrite \eqref{MHD2D} as follows:
\begin{equation}
\begin{cases}
\d_t \tu -\nu \D \tu+\frac1{\ee} \tu \wedge e_3 +\tu\cdot \n \tu -\tb\cdot \n \tb=-\n \tq^0-\frac1{\ee} \n\tp^1,\\
\d_t \tb -\nu' \D \tb +\tu\cdot \n \tb -\tb\cdot \n \tu=0,\\
\div \tu=0,\\
{(\tu,\tb)}_{|t=0}=(\tuo, \tbo).
\end{cases}
\label{MHD2Db}
\tag{$2D-MHD^3$}
\end{equation}
If we consider $\Ue-(\tu, \tb)$, the initial data becomes $(\voe,\co)$ (we are now able to adapt the methods to prove the Leray theorem) but some linear terms (advected by $\tu$) appear. As dispersion only affects the velocity part, we thus need to filter these terms and this is the reason why we introduced System \eqref{Mag} and define $\De$ as follows:
\begin{equation}
 \De = \left(\begin{array}{c}\ve \\ \de\end{array}\right) \overset{def}{=} \Ue-\left(\begin{array}{c}\tu \\ \tb +c\end{array}\right) .
\end{equation}
From now on we will study $\De$ which satisfies the following modified MHD system:
\begin{equation}
\begin{cases}
\d_t \ve -\nu \D \ve +\frac{1}{\ee} \ve\wedge e_3 +\ve\cdot \n \ve +\ve\cdot \n \tu +\tu\cdot \n \ve -\de\cdot \n \de -\de\cdot \n (\tb +c) -(\tb +c)\cdot \n \de\\
\hspace{9.5cm}
 =-\n \qe + \tb\cdot \n c +c\cdot \n \tb +c \cdot \n c,\\
\d_t \de -\nu' \D \de +\ve\cdot \n \de -\de\cdot \n \ve + \tu\cdot \n \de -\de\cdot \n \tu +\ve\cdot \n (\tb +c) -(\tb +c)\cdot \n \ve=0,\\
\div \ve=0,\\
{(\ve,\de)}_{|t=0}=(\voe,0).
\end{cases}
\label{MHDmodif}
\end{equation}
\begin{rem}
 \sl{\begin{enumerate}
      \item Notice that we "pay" the use of $c$ in the magnetic field equation with the presence in the velocity equation of a new external force term: $\tb\cdot \n c +c \cdot \n \tb +c\cdot \n c$.
      \item Once more, for all time $\div \de=0$.
     \end{enumerate}
}
\end{rem}

\subsection{Study of the limit systems}

Let us begin with System \eqref{MHD2D}, which is the three components version of the classical 2D-MHD system.

\begin{thm}
 \sl{For any $\tuo, \tbo\in L^2(\R^2)^3$ there exists a unique global solution $(\tu,\tb)\in \dot{F}^0(\R^2)$. Moreover for any $t\geq 0$,
 \begin{multline}
 \|\tu(t)\|_{L^2(\R^2)}^2 +\|\tb(t)\|_{L^2(\R^2)}^2 +2\nu \int_0^t \|\n_h\tu(\tau)\|_{L^2(\R^2)}^2 d\tau +2\nu' \int_0^t \|\n_h\tb(\tau)\|_{L^2(\R^2)}^2 d\tau\\
 \leq \|\tuo\|_{L^2(\R^2)}^2 +\|\tbo\|_{L^2(\R^2)}^2,
 \label{estimtutb}
\end{multline}
Propagation of regularity: moreover, if in addition $\tu_0,\tb_0 \in \dot{H}^s$ for some $s\in]-1,1[$, then for any $t\geq 0$ we have (for $C_{\nu,\nu'}=C(\frac1{\nu}+\frac1{\nu'})^2$ where $C$ is a universal constant),
\begin{multline}
 \|\tu(t)\|_{\dot{H}^s(\R^2)}^2 +\|\tb(t)\|_{\dot{H}^s(\R^2)}^2 +\int_0^t \left(\nu \|\n_h\tu(\tau)\|_{\dot{H}^s(\R^2)}^2 +\nu' \|\n_h\tb(\tau)\|_{\dot{H}^s(\R^2)}^2\right) d\tau\\
 \leq \left(\|\tuo\|_{\dot{H}^s(\R^2)}^2 +\|\tbo\|_{\dot{H}^s(\R^2)}^2\right) e^{C_{\nu, \nu'}\left(\|\tu_0\|_{L^2(\R^2)}^2 +\|\tb_0\|_{L^2(\R^2)}^2 \right)},
 \label{estimtutbHs}
\end{multline}
\label{ThExistlim}
 }
\end{thm}
\textbf{Proof:} first we recall that the classical 2D-MHD system (without rotation and with values in $\R^2$) can be solved using the same ideas as in the case of the Navier-Stokes system: a priori estimates and Friedrichs-type scheme coupled with the fact that for any functions $f,g,h$, if $\div f=0$ then
 \begin{equation}
 (f\cdot \n g,h)_{L^2} +(f\cdot \n h,g)_{L^2}=0.
  \label{simplif}
 \end{equation}
This is also what we will do in details for System \eqref{MHDmodif} below. Similarly as for the limit system of the rotating fluids system, the first two lines of the velocity and magnetic field equations in System \eqref{MHD2D} (2D-MHD with values in $\R^3$) form the usual 2D-MHD system with $(\tuo^h,\tbo^h)$ as initial data. Classical results immediately give unique global existence of $(\tu^h,\tb^h)$ (in $2D$ the Leray and Fujita-Kato solutions coincide). Then, the third lines of each equation are linked through the following linear system (transported by $\tu^h,\tb^h$) that we easily globally solve:
$$
\begin{cases}
\d_t \tu^3 -\nu \D_h \tu^3 +\tu^h\cdot \n_h \tu^3 -\tb^h\cdot \n_h \tb^3=0,\\
\d_t \tb^3 -\nu' \D_h \tb^3 +\tu^h\cdot \n_h \tb^3 -\tb_h\cdot \n_h \tu^3=0,\\
{(\tu^3,\tb^3)}_{|t=0}=(\tuo^3, \tbo^3),
\end{cases}
$$
The estimates in $\dot{H}^s$ are obtained very similarly as for Navier-Stokes. $\blacksquare$

Let us now turn to System \eqref{Mag}:
\begin{thm}
 \sl{Let $\nu'>0$ and $\tu \in \dot{E}^0(\R^2)$. For any $s\in ]-1,1[$ and any $\co \in \dot{H}^s(\R^3)$, there exists a unique global solution $c\in \dot{E}^s(\R^3)$ of \eqref{Mag}. Moreover there exists a constant $C'>0$ (depending on $s$) such that $c$ satisfies for any $t\geq 0$:
 \begin{multline}
  \|c(t)\|_{\dot{H}^s(\R^3)}^2 +\nu' \int_0^t \|\n c(\tau)\|_{\dot{H}^s(\R^3)}^2 d\tau\\
 \leq \|\co\|_{\dot{H}^s(\R^3)}^2 \exp\left\{\frac{C'}{\nu'} \int_0^t (1+\frac1{(\nu')^2} \|\tu (\tau)\|_{L^2(\R^2)}^2) \|\n_h \tu (\tau)\|_{L^2(\R^2)}^2 d\tau\right\}.
 \label{estimcHs}
 \end{multline}
\label{Thc}
 }
\end{thm}
\textbf{Proof:} this linear equation is easily solved (see also Remark \ref{Remdivc}). To obtain the estimate, we compute the innerproduct in $\dot{H}^s$ with $c$:
\begin{equation}
 \frac12 \frac{d}{dt} \|c\|_{\dot{H}^s}^2 +\nu' \|\n c\|_{\dot{H}^s}^2 \leq A_1 +A_2 \overset{def}{=} |(\tu\cdot \n c|c)_{\dot{H}^s}| +|(c\cdot \n \tu|c)_{\dot{H}^s}|.
\end{equation}
To deal with products of $2D$ functions (i.-e. depending on $x_h$) with $3D$ functions, we use Proposition \ref{prod2D3D} (that will impose restrictions for $s$, see also Remark \eqref{ProdRq}). Similarly to what we did in \cite{FCRF}, we obtain that when $s\in]-1,\frac32[$:
\begin{multline}
 A_1 \leq \|\tu\cdot \n c\|_{\dot{H}^{s-1}} \|c\|_{\dot{H}^{s+1}} \leq C\|\tu\|_{\dot{H}^\frac12(\R^2)} \|\n c\|_{\dot{H}^{s-\frac12}(\R^3)} \|c\|_{\dot{H}^{s+1}(\R^3)}\\
 \leq C\|\tu\|_{L^2}^\frac12 \|\n_h\tu\|_{L^2}^\frac12 \|c\|_{\dot{H}^s}^\frac12 \|\n c\|_{\dot{H}^s}^\frac32 \leq \frac{\nu'}4 \|\n c\|_{\dot{H}^s}^2 +\frac{C'}{(\nu')^3} \|\tu\|_{L^2}^2 \|\n_h \tu\|_{L^2}^2 \|c\|_{\dot{H}^s}^2,
\end{multline}
and when $s\in]-1,1[$:
\begin{multline}
 A_2 \leq \|c\cdot \n \tu\|_{\dot{H}^{s-1}} \|c\|_{\dot{H}^{s+1}} \leq C\|c\|_{\dot{H}^s(\R^3)}\|\n_h\tu\|_{L^2(\R^2)} \|\n c\|_{\dot{H}^s(\R^3)}\\
 \leq \frac{\nu'}4 \|\n c\|_{\dot{H}^s}^2 +\frac{C'}{\nu'} \|\n_h \tu\|_{L^2}^2 \|c\|_{\dot{H}^s}^2.
\end{multline}
Using the Gronwall lemma implies \eqref{estimcHs}. $\blacksquare$

\section{Existence of global weak solutions}

\subsection{A priori estimates}

In this section we prove:

\begin{prop}
 \sl{Let $\ee>0$ fixed. There exists a universal constant $C>0$ and a constant $\Eo=\Eo(\nu, \nu',\|\tuo\|_{L^2}, \|\tbo\|_{L^2})>0$ such that for any solution $\De=(\ve,\de) \in \dot{F}^0(\R^3)$ of \eqref{MHDmodif}, we have for all $t\geq 0$,
 \begin{multline}
\|\De\|_{\dot{F}_t^0}^2 \overset{def}{=}\|\ve\|_{L_t^\infty L^2}^2 +\|\de\|_{L_t^\infty L^2}^2 +\nu \|\n \ve\|_{L_t^2 L^2}^2 +\nu' \|\n \de\|_{L_t^2 L^2}^2\\
\leq \left[\|\voe\|_{L^2}^2+C(1+\frac1{\nu})\Big((1+\|\tb\|_{L_t^\infty L^2}^2) \|\n c\|_{L_t^2 L^2}^2 +\|\n_h \tb\|_{L_t^2 L^2}^2 \|c\|_{L_t^\infty L^2}^2\Big)\right]\\
\times \exp \left\{C(1+\frac1{\nu} +\frac1{\nu'}) \Big(\|\n_h \tu\|_{L_t^2 L^2}^2 +\|\n_h \tb\|_{L_t^2 L^2}^2 +\|\n c\|_{L_t^2 \dot{H}^\frac12}^2\Big) \right\}\\
\leq \Eo \left[\|\voe\|_{L^2}^2 +\|\co\|_{L^2}^2\right]e^{\Eo \|\co\|_{\dot{H}^\frac12}^2}
\label{estimaprioriL2}
 \end{multline}
 }
\end{prop}
\textbf{Proof:} computing the innerproduct in $L^2$ of the velocity part of \eqref{MHDmodif} with $\ve$, and of the magnetic part with $\de$, thanks to \eqref{simplif}, we obtain that for all $t\geq 0$,
\begin{equation}
 \frac12 \frac{d}{dt} (\|\ve\|_{L^2}^2 +\|\de\|_{L^2}^2) +\nu \|\n \ve\|_{L^2}^2 +\nu' \|\n \de\|_{L^2}^2 \leq \Sum_{i=1}^9 |B_i|,
\label{estimL2}
 \end{equation}
where we put
$$
\begin{cases}
 B_1 \overset{def}{=} (\ve \cdot \n \tu| \ve)_{L^2}, \quad B_2 \overset{def}{=} (\de \cdot \n \tb| \ve)_{L^2}, \quad B_3 \overset{def}{=} (\de \cdot \n c| \ve)_{L^2},\\
 B_4 \overset{def}{=} (\de \cdot \n \tu| \de)_{L^2}, \quad B_5 \overset{def}{=} (\ve \cdot \n \tb| \de)_{L^2}, \quad B_6 \overset{def}{=} (\ve \cdot \n c| \de)_{L^2},\\
 B_7 \overset{def}{=} (\tb \cdot \n c| \ve)_{L^2}, \quad B_8 \overset{def}{=} (c \cdot \n \tb| \ve)_{L^2}, \quad B_9 \overset{def}{=} (c \cdot \n c| \ve)_{L^2}.
\end{cases}
$$
With the same $2D \times 3D$ product laws used to prove \eqref{estimcHs} but for $s=0$, we obtain that
\begin{multline}
 |B_1| \leq \|\ve\cdot \n \tu\|_{\dot{H}^{-1}} \|\ve\|_{\dot{H}^1} \leq \|\ve\|_{L^2(\R^3)}\|\n_h\tu\|_{L^2(\R^2)} \|\n \ve\|_{L^2(\R^3)}\\
 \leq \frac{\nu}6 \|\n \ve\|_{L^2}^2 +\frac{C}{\nu} \|\n_h \tu\|_{L^2}^2 \|\ve\|_{L^2}^2.
 \label{estimB1}
\end{multline}
Similarly, we obtain:
$$
\begin{cases}
\vspace{0.2cm}
   |B_2| \leq \frac{\nu}6 \|\n \ve\|_{L^2}^2 +\frac{C}{\nu} \|\n_h \tb\|_{L^2}^2 \|\de\|_{L^2}^2,\\
   \vspace{0.2cm}
   |B_4|+|B_5| \leq \frac{\nu'}3 \|\n \de\|_{L^2}^2 +\frac{C}{\nu'} \|\n_h \tu\|_{L^2}^2 \|\de\|_{L^2}^2 +\frac{C}{\nu'} \|\n_h \tb\|_{L^2}^2 \|\ve\|_{L^2}^2,\\
   |B_7|+|B_8| \leq \frac{\nu}6 \|\n \ve\|_{L^2}^2 +\frac{C}{\nu} \|\tb\|_{L^2}^2 \|\n c\|_{L^2}^2 +\frac{C}{\nu} \|\n_h \tb\|_{L^2}^2 \|c\|_{L^2}^2.
\end{cases}
$$
The last three terms do not feature products with 2D functions and will require additional regularity from $c$, we obtain:
$$
 \begin{cases}
 \vspace{0.2cm}
 |B_3|+|B_6| \leq \|\de\|_{L^6} \|\n c\|_{L^3} \|\ve\|_{L^2} \leq \frac{\nu'}6 \|\n \de\|_{L^2}^2 +\frac{C}{\nu'} \|\n c\|_{\dot{H}^\frac12}^2 \|\ve\|_{L^2}^2,\\
|B_9| \leq \|c\|_{L^6} \|\n c\|_{L^3} \|\ve\|_{L^2} \leq \frac12 \|\n c\|_{L^2}^2 +\frac12 \|\n c\|_{\dot{H}^\frac12}^2 \|\ve\|_{L^2}^2.
\end{cases}
$$
Gathering the previous estimates into \eqref{estimL2}, we obtain that:
\begin{multline}
\frac{d}{dt} (\|\ve\|_{L^2}^2 +\|\de\|_{L^2}^2) +\nu \|\n \ve\|_{L^2}^2 +\nu' \|\n \de\|_{L^2}^2 \leq\\
C(1+\frac1{\nu} +\frac1{\nu'}) (\|\ve\|_{L^2}^2 +\|\de\|_{L^2}^2)\Big(\|\n_h \tu\|_{L^2}^2 +\|\n_h \tb\|_{L^2}^2 +\|\n c\|_{\dot{H}^\frac12}^2\Big)\\
+C(1+\frac1{\nu}) \Big(\|\tb\|_{L^2}^2 \|\n c\|_{L^2}^2 +\|\n_h \tb\|_{L^2}^2 \|c\|_{L^2}^2 + \|\n c\|_{L^2}^2\Big),
 \end{multline}
which gives \eqref{estimaprioriL2} thanks to the Gronwall lemma, \eqref{estimtutb} and \eqref{estimcHs}. $\blacksquare$

\subsection{Proof of Theorem \ref{ThLeray}}

We will not give details here and refer to \cite{FCStratif1} where we adapt classical arguments to prove the existence of Leray-type weak solutions obtained as limits of a Friedrichs scheme: taking advantage of the divergence-free condition, we express the pressure term which allows to define and solve a scheme involving frequency truncation operators $P_n$. Using the a priori estimates, we are able to prove the solutions are bounded and global. Then we prove that the solution of the scheme converge (as $n$ goes to infinity) towards the desired solution.

Compared to \cite{FCStratif1} where we had to deal with products of the type $f(x_1,x_2,x_3)g(x_3)$ (requiring regularity indices $s<\frac12$ in product laws), the arguments are simpler as we multiply by functions of the type $g(x_1,x_2)$ (which only requires $s<1$). In particular we outline that in \cite{FCStratif2} the restrictions on coefficients forced us to ask for additional low frequency assumptions on the initial data to obtain existence of strong Fujita-Kato solutions.

Let us outline that the change in the pressure comes from the fact that its expression features additional terms compared to the case of Navier-Stokes: the $2D$-pressure (in $L^2L^2(\R^2)$), the rotational term (in $\dot{E}^1(\R^3)$), and the bilinear terms involving products of $3D$ by $2D$-functions (in $L^2L^2(\R^3)$).

\section{Convergence}

This section is devoted to the proof of Theorem \ref{ThCV}. We first prove that the velocity strongly converges to zero. Then we use it to prove convergence of the magnetic field part.

\subsection{Velocity}

The equation satisfied by the velocity is the following:
\begin{multline}
 \d_t \ve -\nu \D \ve +\frac{1}{\ee} \bP(\ve\wedge e_3) = \bP \Big[-\ve\cdot \n \ve -\ve\cdot \n \tu -\tu\cdot \n \ve\\
 +\de\cdot \n \de +\de\cdot \n (\tb +c) +(\tb +c)\cdot \n \de + \tb\cdot \n c +c\cdot \n \tb +c \cdot \n c \Big]
\end{multline}
We wish here to adapt the argument from \cite{FCStratif1} and obtain a global-in-time estimate for $\ve$. As in \cite{FCRF, FCStratif1} as the initial data and various external force terms do not share the same regularity, we will split $\ve=\vea+\veb+\vec$ where the new functions satisfy:
\begin{equation}
 \begin{cases}
   \d_t \vea -\nu \D \vea +\frac{1}{\ee} \bP(\vea\wedge e_3) = \bP \Big[-\ve\cdot \n \ve +\de\cdot \n \de +\de\cdot \n c +c\cdot \n \de +c\cdot \n c\Big] \overset{def}{=}G_{\ee,1}\\
   {\vea}_{|t=0}=0,
 \end{cases}
\label{Systv1}
\end{equation}
\begin{equation}
 \begin{cases}
   \d_t \veb -\nu \D \veb +\frac{1}{\ee} \bP(\veb\wedge e_3) = \bP \Big[-\ve\cdot \n \tu -\tu\cdot \n \ve +\de\cdot \n \tb +\tb \cdot \n \de +\tb\cdot \n c +c \cdot \n \tb \Big] \overset{def}{=}G_{\ee,2}\\
   {\veb}_{|t=0}=0,
 \end{cases}
\label{Systv2}
\end{equation}
\begin{equation}
 \begin{cases}
   \d_t \vec -\nu \D \vec +\frac{1}{\ee} \bP(\vec\wedge e_3) = 0\\
   {\vec}_{|t=0}=\voe,
 \end{cases}
\label{Systv3}
\end{equation}
Thanks to the first point of Proposition \ref{EstimStri} (with $\Fe=\Ge=0$, $\theta\in[0,1]$ and $\sigma_1=0$), we obtain that for any $r\geq 2$, any $p\in [1,\frac2{\theta(1-\frac2{r})}]$ such that
\begin{equation}
 \frac2{p}=\frac32-\frac3{r}+\theta\big(1-\frac2{r}\big),
 \label{Cond:p3}
\end{equation}
there exists a constant $C_{\theta,r,p,\nu}>0$ such that we have:
$$
 \|\ve^3\|_{\tilde{L}^p\dot{B}_{r,2}^0} \leq C_{\theta,r,p,\nu} \ee^{\frac{\theta}2 (1-\frac2{r})}\|\voe\|_{L^2}.
$$
Moreover for any $r\geq2$, $\theta\in]0,1]$, and any $p$ satisfying \eqref{Cond:p3}, we have:
$$
 p\in [1,\frac2{\theta(1-\frac2{r})}]\; \Longleftrightarrow\; \frac{\theta}2(1-\frac2{r}) \leq \frac1{p}\leq 1\; \Longleftrightarrow\; r\geq 2 \mbox{ and }\frac{3+2\theta}{r}\geq \theta-\frac12.
$$
And finally, for any $\theta \in]0,1]$ we have:
\begin{equation}
 \|\ve^3\|_{\tilde{L}^{\frac4{(3+2\theta)(1-\frac2{r})}}\dot{B}_{r,2}^0} \leq C_{\theta,r,\nu} \ee^{\frac{\theta}2 (1-\frac2{r})}\|\voe\|_{L^2},
 \label{Estimv3}
 \mbox{ for any }
 \begin{cases}
 \vspace{1mm}
 r\in[2,\infty] & \mbox{if }\theta\in[0,\frac12],\\
 r\in[2,2+\frac4{\theta-\frac12}] & \mbox{if }\theta\in]\frac12,1].
\end{cases}
\end{equation}
We now turn to $v_\ee^2$: using the divergence-free condition and the modified product laws from Proposition \ref{prod2D3D}, we have for some small $\eta>0$:
\begin{multline}
 \|\ve\cdot \n \tu\|_{\dot{H}^{-\frac12}} =\|\div(\ve\otimes \tu)\|_{\dot{H}^{-\frac12}}\leq \|\ve\otimes \tu\|_{\dot{H}^\frac12} \leq C\|\ve\|_{\dot{H}^{1-\eta}} \|\tu\|_{\dot{H}^{\frac12+\eta}}\\
 \leq C\|\ve\|_{L^2(\R^3)}^\eta \|\n\ve\|_{L^2(\R^3)}^{1-\eta} \|\tu\|_{L^2(\R^2)}^{\frac12-\eta} \|\n_h \tu\|_{L^2(\R^2)}^{\frac12+\eta}.
\end{multline}
The second term in the product is $L^\frac2{1-\eta}$ in time, the last one is $L^\frac2{\frac12+\eta}$, the rest being $L^\infty$, we obtain (thanks to \eqref{estimaprioriL2}) that for some constant depending on $\nu$:
\begin{equation}
 \|\ve\cdot \n \tu\|_{L^\frac43 \dot{H}^{-\frac12}} \leq C\|\ve\|_{L^\infty L^2}^\eta \|\n\ve\|_{L^2 L^2}^{1-\eta} \|\tu\|_{L^\infty L^2}^{\frac12-\eta} \|\n_h \tu\|_{L^2 L^2}^{\frac12+\eta} \leq C\|\De\|_{\dot{F}^0(\R^3)} \|\tu\|_{\dot{E}^0(\R^2)}.
\end{equation}
The other terms in $G_{\ee,2}$ are estimated similarly, so that
for some constant $\Eo>0$ (depending on $\nu, \nu',\|\tuo\|_{L^2}, \|\tbo\|_{L^2}$) we have:
\begin{equation}
 \|G_{\ee,2}\|_{L^\frac43 \dot{H}^{-\frac12}} \leq\Eo \left[\|\voe\|_{L^2}^2 +\|\co\|_{L^2}^2\right]^\frac12 e^{\Eo \|\co\|_{\dot{H}^\frac12}^2}
\label{estim:Ge2}
\end{equation}
Thanks to the second point of Proposition \ref{EstimStri} (with $(\Fe,\Ge)=(0,G_{\ee,2})$, $\theta\in]0,1]$, $k=\frac43$ and $\sigma_2=-\frac12$), we obtain that for any $r\geq 2$, any $p\in [2,\frac2{\theta(1-\frac2{r})}[$ also satisfying \eqref{Cond:p3} there exists a constant $C_{\theta,r,p,\nu}>0$ such that:
$$
 \|\ve^2\|_{\tilde{L}^{p}\dot{B}_{r,2}^0} \leq C_{\theta,r,p,\nu} \ee^{\frac{\theta}2 (1-\frac2{r})} \|G_{\ee,2}\|_{L^\frac43 \dot{H}^{-\frac12}}.
$$
Similarly, for any $r\geq2$, $\theta\in]0,1]$, and any $p$ defined as in \eqref{Cond:p3}, we have:
$$
p\in [2,\frac2{\theta(1-\frac2{r})}[\; \Longleftrightarrow\; r\in]2,2+\frac2{\theta+\frac12}],
$$
therefore, for any $\theta\in]0,1]$, any $r\in]2,2+\frac2{\theta+\frac12}]$ we have:
\begin{equation}
\|\ve^2\|_{\tilde{L}^{\frac4{(3+2\theta)(1-\frac2{r})}}\dot{B}_{r,2}^0} \leq C_{\theta,r,\nu} \ee^{\frac{\theta}2 (1-\frac2{r})}\|G_{\ee,2}\|_{L^\frac43 \dot{H}^{-\frac12}}.
\label{Estimv2}
\end{equation}
Gathering \eqref{Estimv2}, \eqref{Estimv3} and \eqref{estim:Ge2} (and using the assumption \eqref{Hyp:Maj}) we obtain that for any $\theta\in[0,1]$, any $r \in]2,2+\frac2{\theta+\frac12}]$ there exists a constant $\Fo$ (depending on $\nu, \nu',\Co,\theta, r $) such that:
$$
\|\ve^2+\ve^3\|_{\tilde{L}^{\frac4{(3+2\theta)(1-\frac2{r})}}\dot{B}_{r,2}^0} \leq \Fo \ee^{\frac{\theta}2 (1-\frac2{r})}.
$$
As $0<\theta\leq 1$, the upper bound $2+\frac2{\theta+\frac12}$ lies in $[\frac{10}3,6]$, so for any $r\in]2,6]$, there exists some $\theta\in]0,1]$ such that $r\in ]2,2+\frac2{\theta+\frac12}]$ if, and only if, $\theta\leq \frac{6-r}{2(r-2)}$. Of course, this latter quantity goes to infinity when $r>2$ goes to $2$, so introducing
$$
 \theta_r\overset{def}{=} \min(1,\frac{6-r}{2(r-2)})
$$
and the following quantity (lying in $[2,\infty[$):
\begin{equation}
  \aa_r\overset{def}{=} \frac4{(3+2\theta_r)(1-\frac2{r})}= \max(2,\frac4{5(1-\frac2{r})}),
  \label{def:alphar}
\end{equation}
then for any $r\in]2,6]$ there exists a constant $\Fo$ (now depending on $r, \nu, \nu',\Co$) such that:
\begin{equation}
\|\ve^2+\ve^3\|_{\tilde{L}^{\aa_r} \dot{B}_{r,2}^0} \leq \Fo \ee^{\frac{\theta_r}2 (1-\frac2{r})}= \Fo \ee^{\frac1{2r} \min(r-2,3-\frac{r}2)}.
\label{Estimv2v3}
\end{equation}
Turning to $\ve^1$, thanks to the classical Sobolev product laws, we have:
$$
\|\ve \cdot \n \ve\|_{\dot{H}^{-\frac12}} \leq C\|\ve\|_{\dot{H}^1} \|\n \ve\|_{L^2} \leq C\|\n \ve\|_{L^2}^2.
$$
The other terms in $G_{\ee,1}$ are estimated similarly, and we obtain that:
$$
\|G_{\ee,1}\|_{L^1 \dot{H}^{-\frac12}} \leq C(\|\n \ve\|_{L^2 L^2}^2 +\|\n \de\|_{L^2 L^2}^2 +\|\n c\|_{L^2 L^2}^2) \leq\Eo \left[\|\voe\|_{L^2}^2 +\|\co\|_{L^2}^2\right] e^{\Eo \|\co\|_{\dot{H}^\frac12}^2}.
$$
Using once more the first point of Proposition \ref{EstimStri} (with $(\Fe,\Ge)=(G_{\ee,1},0)$, $\theta\in]0,1]$ and $\sigma_1=-\frac12$), we obtain that for any $r\geq 2$, any $p\in [1,\frac2{\theta(1-\frac2{r})}]$ such that
\begin{equation}
 \frac2{p}=2-\frac3{r}+\theta(1-\frac2{r}),
 \label{condp1r1}
\end{equation}
there exists a constant $C_{\theta,r,p,\nu}>0$ such that we have:
$$
 \|\ve^1\|_{\tilde{L}^{p}\dot{B}_{r,2}^0} \leq C_{\theta,r,p,\nu} \ee^{\frac{\theta}2 (1-\frac2{r})}\|G_{\ee,1}\|_{L^1 \dot{H}^{-\frac12}}.
$$
Similarly, for any $r\geq 2$, $\theta\in[0,1]$ and $p$ satisfying Condition \eqref{condp1r1}, then
$$
 p\in [1,\frac2{\theta(1-\frac2{r})}]\; \Longleftrightarrow\; r\in[2,2+\frac3{\theta}].
$$
We can conclude similarly as before that for any $\theta\in[0,1]$ and $r\in[2,2+\frac3{\theta}]$, we have :
\begin{equation}
\|\ve^1\|_{\tilde{L}^\frac4{1+(3+2\theta)(1-\frac2{r})}\dot{B}_{r,2}^0}\\
\leq \Fo \ee^{\frac{\theta}2 (1-\frac2{r})},
 \label{Estimv1}
\end{equation}
Thus for any $r>2$, there exists some $\theta\in[0,1]$ such that $r\leq 2+\frac3{\theta}$ if, and only if $\theta\leq \frac3{r-2}$ so that defining
$$
\theta_r'\overset{def}{=} \min(1,\frac3{r-2})
$$
and the following quantity (lying in $[1,4[$):
\begin{equation}
  \bb_r\overset{def}{=} \frac4{1+(3+2\theta_r')(1-\frac2{r})}= \max(1,\frac2{3-\frac5{r}}),
  \label{def:betar}
\end{equation}
then for any $r>2$, there exists a constant $\Fo$ (depending on $r, \nu, \nu',\Co$) such that:
\begin{equation}
\|\ve^1\|_{\tilde{L}^{\bb_r} \dot{B}_{r,2}^0} \leq \Fo \ee^{\frac{\theta_r'}2 (1-\frac2{r})}= \Fo \ee^{\frac1{2r} \min(r-2,3)}.
\label{Estimv1b}
\end{equation}
Gathering \eqref{Estimv2v3} and \eqref{Estimv1b} for $r\in ]2,6[$ gives \eqref{ThCV-estim1}.
Let us now turn to the second simpler (and local-in-time) estimates from Theorem \ref{ThCV}.
\\
Coming back to \eqref{def:alphar} and \eqref{Estimv2v3}, for any $r\in]2,6]$ we have $\aa_r\geq2$ and thanks to Propositions \ref{injectionLr} and \ref{Propermut} we obtain that for any $r\in]2,6]$:
\begin{equation}
\|\ve^2+\ve^3\|_{L^{\aa_r} L^r} \leq C \|\ve^2+\ve^3\|_{L^{\aa_r} \dot{B}_{r,2}^0} \leq C \|\ve^2+\ve^3\|_{\tilde{L}^{\aa_r} \dot{B}_{r,2}^0} \leq \Fo \ee^{\frac1{2r} \min(r-2,3-\frac{r}2)},
\end{equation}
which leads, for any $r\in]2,6]$ and $t\geq 0$ to:
\begin{multline}
\|\ve^2+\ve^3\|_{L_t^2 L^r}\leq t^{\frac12-\frac1{\aa_r}} \|\ve^2+\ve^3\|_{L^{\aa_r} L^r}\\
\leq \Fo t^{\frac12-\frac1{\aa_r}} \ee^{\frac1{2r} \min(r-2,3-\frac{r}2)} =\Fo t^{\max(0, \frac5{2r}-\frac34)} \ee^{\frac1{2r} \min(r-2,3-\frac{r}2)}.
\label{Estimv2v3b}
\end{multline}
Let us recall that, thanks to Proposition \ref{Propermut}, the injection $\tilde{L}_t^p \dot{B}_{r,2}^0 \hookrightarrow  L_t^p \dot{B}_{r,2}^0$ is true if, and only if $p\geq 2$. Considering \eqref{def:betar}, we have to restrict the domain in $r$: $\bb_r\geq 2 \; \Longleftrightarrow\; r\in]2, \frac52]$, and in this case we have $\bb_r=\frac2{3-\frac5{r}}$ and thanks to \eqref{Estimv1b},
$$
 \|\ve^1\|_{L_t^2 L^r} \leq t^{\frac12-\frac1{\bb_r}}\|\ve^1\|_{L_t^{\bb_r} L^r} \leq t^{\frac12-\frac1{\bb_r}} \|\ve^1\|_{\tilde{L}^{\bb_r} \dot{B}_{r,2}^0} \leq \Fo t^{\frac5{2r}-1} \ee^{{\frac12-\frac1{r}}}.
$$
Gathering this with \eqref{Estimv2v3b}, we obtain that for any $r\in]2, \frac52]$, there exists a constant $\Fo$ (depending on $\nu, \nu',\Co,r$) such that for any $t\geq 0$,
\begin{multline}
 \|\ve\|_{L_t^2 L^r} \leq \|\ve^1\|_{L_t^2 L^r} +\|\ve^2+\ve^3\|_{L_t^2 L^r} \leq \Fo t^{\frac5{2r}-1} \ee^{{\frac12-\frac1{r}}} +\Fo t^{\frac5{2r}-\frac34} \ee^{\frac1{2r} \min(r-2,3-\frac{r}2)}\\
 \leq \Fo (t^{\frac5{2r}-1} +t^{\frac5{2r}-\frac34}) \ee^{{\frac12-\frac1{r}}} \leq \Fo \max(1,t)^{\frac5{2r}-\frac34} \ee^{\frac12-\frac1{r}}
 \label{Estim:intermediaire}
\end{multline}
For $r\in[\frac52,6[$, we interpolate the previous estimates for $r=\frac52$,
$$
\|\ve\|_{L_t^2 L^\frac52} \leq \Fo \max(1,t)^\frac14 \ee^\frac1{10},
$$
with the following estimates (given by \eqref{estimaprioriL2}):
$$
\|\ve\|_{L_t^2 L^6} \leq C\|\ve\|_{L^2 \dot{H}^1} \leq \Fo,
$$
and we obtain that for any $r\in[\frac52,6]$,
\begin{equation}
\|\ve\|_{L_t^2 L^r} \leq \Fo \max(1,t)^{\frac5{28}(\frac6{r}-1)} \ee^{\frac1{14}(\frac6{r}-1)}.
\label{Estim:intermediaire2}
\end{equation}
Combining \eqref{Estim:intermediaire} with \eqref{Estim:intermediaire2}, we finally obtain that for any $r\in]2,6]$,
$$
\|\ve\|_{L_t^2 L^r} \leq \Fo \max(1,t)^{\max\left(\frac5{2r}-\frac34,\frac5{28}(\frac6{r}-1)\right)} \ee^{\min\left(\frac12-\frac1{r},\frac1{14}(\frac6{r}-1)\right)},
$$
which gives \eqref{ThCV-estim2}. $\blacksquare$

\subsection{Magnetic field}
We now turn to the magnetic-field part: thanks to \eqref{estimaprioriL2}, $\de$ is uniformly bounded (with respect to $\ee$) and we can extract a subsequence $d_{\ee'}$ that weakly converges to some $d\in L^\infty L^2 \cap L^2 \dot{H}^1$ and satisfies the following system:
$$
\begin{cases}
 \partial_t d-\nu'\D d+\tu\cdot \n d-d\cdot \n \tu=0,\\
 d_{|t=0}=0
\end{cases}
$$
\begin{rem}
 \sl{We do not give details for the weak convergence of the nonlinear terms, for which the arguments are the same as for the classical Navier-Stokes system.}
\end{rem}
Thanks to Theorem \ref{Thc}, we immediately obtain that $d\equiv 0$, and more generally, we obtained that the only possible weak limit of a subsequence is zero which ensures that the whole sequence itself weakly converges to zero. With similar arguments as in the classical proof for the Navier-Stokes system (see also \cite{FC1} and the end of Section 4.1 in \cite{FCStratif1}), using once more \eqref{estimaprioriL2} we have that for any $r\in]2,6[$,
$$
\de \underset{\ee\rightarrow 0}{\longrightarrow}0 \mbox{ in } L_{loc}^2 L_{loc}^r,
$$
which ends the proof of Theorem \ref{ThCV}. $ \blacksquare$

\section{Strong solutions}

\subsection{Precise statement of the results}

We first recall that for an initial data as in \eqref{Initdata2}, every term $c$ in System \eqref{MHDmodif} is now replaced by $\ce$.

Let us begin with stating more precisely our main results concerning strong solutions (see the simplified Theorem \ref{ThCVStrongsimplif} from the introduction): first local existence theorem, then global existence and convergence theorem for very large ill-prepared initial data.

\begin{thm}(Existence of local Fujita-Kato strong solutions)
 \sl{Let $\ee>0$ fixed, $\tu,\tb\in L^2(\R^3)$, $\voe\in\dot{H}^\frac12$ and $\coe\in H^\frac12$. For any initial data as in \eqref{Initdata2}, there exists a unique local solution $\De$ of \eqref{MHDmodif} with lifespan $T_\ee^*>0$ such that for any $T<T_\ee^*$, $\De\in \dot{F}_T^\frac12$. Moreover, the following properties are true:
 \begin{itemize}
  \item Regularity propagation: if in addition $\Uoe,\coe\in \dot{H}^s$ for some $s\in]-1,1[$ then for any $t<T_\ee^*$, $\De\in \dot{E}_t^\frac12 \cap \dot{E}_t^s$ and for some $C=C_{\nu,\nu',\|\tuo\|_{L^2}, \|\tbo\|_{L^2}}>0$,
\begin{multline}
  \|\De\|_{\dot{F}_t^s} \leq C \Big(\|\voe\|_{\dot{H}^s}^2 +\|\coe\|_{\dot{H}^s}^2 \Big) \times e^{C\Big(1+\|\coe\|_{\dot{H}^\frac12}^2 +\|\coe\|_{\dot{H}^\frac12}^4\big)}\\
  \times \exp \left\{C \int_0^t (\|\n\ve(t')\|_{\dot{H}^\frac12}^2 +(\|\n\de(t')\|_{\dot{H}^\frac12}^2)dt'\big)\right\}
\end{multline}
  \item Continuation criterion: $\int_0^{T_\ee^*} \|\n \De (\tau)\|_{\dot{H}^\frac12}^2 d\tau <\infty \Longrightarrow T_\ee^*=\infty$.
  \item Global existence for small data: there exist a constant $c_0>0$ such that if $\|\voe\|_{\dot{H}^\frac12}\leq c_0 \nu$ (we recall that $d_{\ee|t=0}=0$) then $T_\ee^*=+\infty$ and for any $t\geq 0$ (for the same $C$):
  $$
  \|\De\|_{\dot{F}_t^\frac12} \leq C \Big(\|\voe\|_{\dot{H}^\frac12}^2 +\|\coe\|_{\dot{H}^\frac12}^2 \Big) \times e^{C\Big(1+\|\coe\|_{\dot{H}^\frac12}^2 +\|\coe\|_{\dot{H}^\frac12}^4\big)}
  $$
 \end{itemize}
 }
 \label{Th0FK}
\end{thm}
\textbf{Proof:} the proof involves estimates similar to what we do below (but simpler) so we will only point out that to obtain the estimate for the regularity propagation, we have to bound a little differently the quadratic terms from system \eqref{MHDmodif}, for instance:
$$
 |(\ve\cdot \n \ve|\ve)_{\dot{H}^s}|\leq \frac{\nu}{10} \|\n \ve\|_{\dot{H}^s} +\frac{C}{\nu} \|\ve\|_{\dot{H}^s}^2 \|\n \ve\|_{\dot{H}^\frac12}^2.\blacksquare
$$
If we wish to obtain convergence for strong solutions, as in the cases of the rotating fluids, primitive or stratified Boussinesq systems, any direct approach involving System \eqref{MHDmodif} is bound to fail due to the following term appearing in the velocity equation: $\tb\cdot \n \ce +\ce\cdot \n \tb +\ce \cdot \n \ce$. Indeed, this term is bounded but does not go to zero as $\ee$ goes to zero.\\
The usual way to reach convergence, first consists in introducing  the oscillations $\We$ that can absorb this constant term, and go to zero thanks to dispersive properties as they solve the following linear system:
\begin{equation}
 \begin{cases}
   \d_t \We -\nu \D \We +\frac1{\ee} \bP (\We \wedge e_3)=\bP (\tb\cdot \n \ce +\ce\cdot \n \tb +\ce \cdot \n \ce),\\
 W_{\ee|t=0}=\voe.
 \end{cases}
\label{SystWe}
 \end{equation}
If we set:
$$
\tDe \overset{def}{=}\Ue -\left(\begin{array}{c}\tu +\We \\ \tb +\ce \end{array}\right) =\De-\left(\begin{array}{c}\We \\ 0 \end{array}\right)= \left(\begin{array}{c}\ve-\We \\ \de \end{array}\right) =\left(\begin{array}{c}\dde \\ \de \end{array}\right),
$$
then it satisfies the following system
\begin{equation}
 \begin{cases}
  \d_t \dde -\nu \D\dde +\frac1{\ee} \bP (\We \wedge e_3) \dde = \Sum_{i=1}^{13} F_i,\\
  \d_t \de -\nu' \D \de = \Sum_{i=1}^{14} G_i,\\
  {(\dde,\de)}_{|t=0}= (0,0),
 \end{cases}
\label{Systdde}
\end{equation}
where we define:
\begin{equation}
\begin{cases}
 F_1 \overset{def}{=}-\mathbb{P}(\dde \cdot \n \dde), \quad F_2 \overset{def}{=}-\mathbb{P}(\dde \cdot \n \We), \quad F_3 \overset{def}{=}-\mathbb{P}(\We \cdot \n \dde),\\
 F_4 \overset{def}{=}-\mathbb{P}(\We \cdot \n \We),\quad F_5 \overset{def}{=}-\mathbb{P}(\dde \cdot \n_h \tu),\quad F_6 \overset{def}{=}-\mathbb{P}(\We \cdot \n_h \tu),\\
 F_7 \overset{def}{=}-\mathbb{P}(\tu \cdot \n \dde),\quad F_8 \overset{def}{=}-\mathbb{P}(\tu \cdot \n \We),\quad F_9 \overset{def}{=}\mathbb{P}(\de \cdot \n \de),\\
F_{10} \overset{def}{=}\mathbb{P}(\de \cdot \n_h \tb),\quad
F_{11} \overset{def}{=}\mathbb{P}(\de \cdot \n \ce), \quad F_{12} \overset{def}{=}\mathbb{P}(\tb \cdot \n \de),
\\ F_{13} \overset{def}{=}\mathbb{P}(\ce \cdot \n \de), \quad G_1\overset{def}{=}-\dde\cdot \n \de, \quad G_2\overset{def}{=}-\We \cdot \n \de, \quad G_3\overset{def}{=} \de \cdot \n \dde,\\
G_4\overset{def}{=} \de \cdot \n \We, \quad G_5\overset{def}{=} -\tu\cdot \n \de, \quad G_6\overset{def}{=} \de\cdot \n_h \tu,\\
G_7\overset{def}{=} -\dde \cdot \n_h \tb, \quad G_8\overset{def}{=} -\dde \cdot \n \ce, \quad G_9\overset{def}{=} -\We\cdot \n_h \tb, \quad G_{10}\overset{def}{=} -\We\cdot \n \ce,\\
G_{11}\overset{def}{=} \tb \cdot \n \dde, \quad G_{12}\overset{def}{=} \ce \cdot \n \dde, \quad G_{13}\overset{def}{=} \tb\cdot \n \We, \quad G_{14}\overset{def}{=} \ce\cdot \n \We.
\end{cases}
 \label{systdde}
\end{equation}

\begin{thm}
 \sl{For any $\Co\geq 1$ (size), any $\delta\in]0,\frac16]$ (extra regularity), any $\mu\in[0,1[$ (as close to 1 as we wish), any $\gamma \in [0, \frac5{12}\delta]$, introducing $\eta_0$ such that $\gamma =\frac{\delta}2 (1-2\eta_0)$ (it implies $\eta_0 \in [\frac1{12}, \frac12]$) there exists $\ee_0>0$ and $\Ko,\Fo\geq 1$ (depending on $\nu, \nu', \Co, \delta, \gamma$) such that for any $\ee\in]0,\ee_0]$ and any initial data as in \eqref{Initdata2} with $\tuo, \tbo \in H^\delta(\R^2)$, $(\voe, \coe) \in (\dot{H}^\frac12 (\R^3) \cap \dot{H}^{\frac12+\delta}(\R^3)) \times H^{\frac12+\delta}(\R^3)$ satisfying:
 \begin{equation}
  \|\voe\|_{\dot{H}^{\frac12+\mu\delta} \cap \dot{H}^{\frac12+\delta}}\leq \Co \ee^{-\gamma}, \mbox{ and } \|\coe\|_{H^{\frac12+\delta}}\leq \left(\Ko |\ln \ee|\right)^{\frac14},
  \label{Hypinit}
 \end{equation}
then $T_\ee^*=+\infty$ and for any $s\in[\frac12,\frac12 +\frac{\eta_0}{10}\delta]=[\frac12,\frac12 +\frac1{10}(\frac{\delta}2-\gamma)]$ we have:
\begin{equation}
 \|\tDe\|_{\dot{E}^s} =\|(\dde,\de)\|_{\dot{E}^s} \leq \Fo \ee^{\frac{\eta_0}{18} \delta}= \Fo \ee^{\frac{1}{18} (\frac{\delta}2-\gamma)}.
 \label{EstimTHStrongA}
\end{equation}
Moreover, if $\voe\in\dot{H}^{\frac12-\delta}\cap \dot{H}^{\frac12+\delta}$, the previous estimates remain true for any $s\in[\frac12-\frac{\eta_0}{10}\delta,\frac12 +\frac{\eta_0}{10}\delta]$, we can get rid of the oscillations $\We$ and get:
\begin{equation}
 \|\De\|_{L^2 L^\infty} =\|\Ue-(\tu,\tb+\ce)\|_{L^2 L^\infty} \leq \Fo \ee^{\frac{\eta_0}{18} \delta}= \Fo \ee^{\frac{1}{18} (\frac{\delta}2-\gamma)}.
\label{EstimTHStrongB}
\end{equation}
 }
\label{ThCVStrong}
\end{thm}
\begin{rem}
 \sl{\begin{enumerate}
      \item Due to the new term $G_{14}$, we cannot choose $\gamma$ as close to $\frac{\delta}2$ as we did in \cite{FCcompl, FCRF, FCStratif2}. At best we are able to reach $\gamma <\frac37 \delta$ (that is $\eta_0\in]\frac1{14},\frac12]$) but the price to pay is a much smaller power of $\ee$. On the opposite, if we want the final convergence rate to be $\ee^{\frac{\eta_0}2 \delta}$ (as in \cite{FCcompl,FCRF} for instance) we can consider the case $\gamma \in[0,\frac3{10} \delta]$ ($\eta_0 \in [\frac15,\frac12]$). For the sake of simplicity we restrict to the case $\eta_0\geq \frac1{12}$ ($\gamma \in[0,\frac5{12} \delta]$).
      \item The constant $\Ko$ is taylored for \eqref{EstimaprioriS} to hold: with $\Do$ introduced in \eqref{Estimapriori3}, we have:
      $$
      \Ko =\frac{\eta_0 \delta}{144 \Do}.
      $$
      \item The bound $\ee_0$ is specified in particular in \eqref{Condeps1} and \eqref{Condeps2}.
     \end{enumerate}
     }
\end{rem}

\subsection{Proof of Theorem \ref{ThCVStrong}}

\subsubsection{Estimates in $\dot{H}^s$ of external force terms}

We know from Theorem \ref{Th0FK} that there exists a local solution defined on $[0,T_\ee^*[$. Assume by contradiction that $T_\ee^*<+\infty$, then by the continuation criterion, we have:
\begin{equation}
 \int_0^{T_\ee^*} \|\n \De (t)\|_{\dot{H}^\frac12}^2 dt =\int_0^{T_\ee^*} \left(\|\n \ve (t)\|_{\dot{H}^\frac12}^2 +\|\n \de (t)\|_{\dot{H}^\frac12}^2\right) dt =+\infty.
\end{equation}
We can now define the following time (where the universal constant $C$ is given below in the bounds involving $F_1,F_9$ and $G_1,G_3$):
\begin{equation}
 T_\ee=\sup \{t\in[0,T_\ee^*[, \; \forall t'\leq t, \|\dde(t')\|_{\dot{H}^\frac12} +\|\de(t')\|_{\dot{H}^\frac12}\leq \frac1{8C}\min(\nu,\nu')\}.
 \label{DefTe}
\end{equation}
As the initial data from System \eqref{Systdde} is zero, we know that $T_\ee>0$. Let us now assume by contradiction that
\begin{equation}
 T_\ee<T_\ee^*.
\end{equation}
Now, performing innerproducts in $\dot{H}^s$ of System \eqref{Systdde} with $(\dde,\de)$, we obtain that for all $t\leq T_\ee$:
\begin{equation}
\begin{cases}
  \frac12 \frac{d}{dt}\|\dde(t)\|_{\dot{H}^s}^2 +\nu \|\n \dde(t)\|_{\dot{H}^s}^2 \leq \Sum_{k=1}^{13} |(F_k(t)|\dde(t))_{\dot{H}^s}|,\\
  \frac12 \frac{d}{dt}\|\de(t)\|_{\dot{H}^s}^2 +\nu' \|\n \de(t)\|_{\dot{H}^s}^2 \leq \Sum_{k=1}^{14} |(G_k(t)|\de(t))_{\dot{H}^s}|.
\end{cases}
\label{Systddde}
\end{equation}
Most of the right-hand-side terms are treated like in the case of the rotating fluids, or very similarly (see \cite{FCRF} for details). We will give minimal details for these terms and focus on the new terms, among them $F_8, G_{13}$ and especially $G_{14}$ which will require most of the work.
\begin{rem}
 \sl{Another change is that contrary to the case of rotating fluids, the oscillations $\We$ now satisfy System \eqref{SystWe} featuring several external force terms of different regularities, which forces us to systematically use the variant of our Strichartz estimates for time integration index $k\in]1,2]$ (see the second point of Propositions \ref{EstimStri} and \ref{EstimStrianiso})}
 \label{rem:chgtFextk}
\end{rem}
The following estimates directly come from \cite{FCRF} (see Sections 2.2 and 3.2) using classical Sobolev product laws:
\begin{equation}
\begin{cases}
\vspace{1mm}
 |(F_1|\dde)_{\dot{H}^s}| \leq C\|\dde\|_{\dot{H}^\frac12} \|\n \dde\|_{\dot{H}^s}^2,\\
 \vspace{1mm}
|(F_2|\dde)_{\dot{H}^s}| \leq \frac{\nu}{52} \|\n \dde\|_{\dot{H}^s}^2 +\frac{C}{\nu} \|\n \We\|_{L^3}^2 \|\dde\|_{\dot{H}^s}^2,\\
\vspace{1mm}
|(F_3|\dde)_{\dot{H}^s}| \leq \frac{\nu}{52} \|\n \dde\|_{\dot{H}^s}^2 +\frac{C}{\nu^3} \|\We\|_{L^6}^4 \|\dde\|_{\dot{H}^s}^2,\\
|(F_4|\dde)_{\dot{H}^s}| \leq \frac{\nu}{52} \|\n \dde\|_{\dot{H}^s}^2 +\frac{C}{\nu^\frac{s}{1-s}} \|\We\|_{L^6}^\frac2{1-s} \|\dde\|_{\dot{H}^s}^2 +C\|\n \We\|_{L^3}^2.
\end{cases}
\label{Estim:ext:1}
\end{equation}
Next $G_4$ and $G_2$ are respectively treated like $F_2$ and $F_3$, and $F_9,G_1$ and $G_3$ are dealt with as $F_1$:
\begin{equation}
\begin{cases}
 \vspace{1mm}
|(G_4|\de)_{\dot{H}^s}| \leq \frac{\nu'}{56} \|\n \de\|_{\dot{H}^s}^2 +\frac{C}{\nu'} \|\n \We\|_{L^3}^2 \|\de\|_{\dot{H}^s}^2,\\
\vspace{1mm}
|(G_2|\dde)_{\dot{H}^s}| \leq \frac{\nu'}{56} \|\n \de\|_{\dot{H}^s}^2 +\frac{C}{(\nu')^3} \|\We\|_{L^6}^4 \|\de\|_{\dot{H}^s}^2,\\
\vspace{1mm}
 |(F_9|\dde)_{\dot{H}^s}| +|(G_3|\de)_{\dot{H}^s}| \leq C\|\de\|_{\dot{H}^\frac12} \left(\|\n \dde\|_{\dot{H}^s}^2 +\|\n \de\|_{\dot{H}^s}^2\right),\\
 |(G_1|\de)_{\dot{H}^s}| \leq C\|\dde\|_{\dot{H}^\frac12} \|\n \de\|_{\dot{H}^s}^2.
\end{cases}
\label{Estim:ext:2}
\end{equation}
Now, let us turn to terms involving products of functions depending on the whole space variable $x=(x_1,x_2,x_3)$ by functions only depending on the horizontal variables $x_h=(x_1,x_2)$: $F_{10},G_6,G_7$ on one hand, $G_5$ on the other hand are respectively dealt with as $F_5$ and $F_7$ which are the same as in \cite{FCRF} and we get that (using the $2D\times 3D$ product laws from Proposition \ref{prod2D3D}, which requires $s\in]-1,1[$):
\begin{equation}
 \begin{cases}
  \vspace{1mm}
 |(F_5|\dde)_{\dot{H}^s}| \leq \frac{\nu}{52} \|\n \dde\|_{\dot{H}^s}^2 +\frac{C}{\nu} \|\n_h \tu\|_{L^2(\R^2)}^2 \|\dde\|_{\dot{H}^s}^2,\\
 \vspace{1mm}
 |(G_6|\de)_{\dot{H}^s}| \leq \frac{\nu'}{56} \|\n \de\|_{\dot{H}^s}^2 +\frac{C}{\nu'} \|\n_h \tu\|_{L^2}^2 \|\de\|_{\dot{H}^s}^2,\\
 \vspace{1mm}
 |(F_{10}|\dde)_{\dot{H}^s}| \leq \frac{\nu}{52} \|\n \dde\|_{\dot{H}^s}^2 +\frac{C}{\nu} \|\n_h \tb\|_{L^2}^2 \|\de\|_{\dot{H}^s}^2,\\
 \vspace{1mm}
 |(G_7|\de)_{\dot{H}^s}| \leq \frac{\nu'}{56} \|\n \de\|_{\dot{H}^s}^2 +\frac{C}{\nu'} \|\n_h \tb\|_{L^2}^2 \|\dde\|_{\dot{H}^s}^2,\\
 \end{cases}
 \label{Estim:ext:3}
\end{equation}
and
\begin{equation}
 \begin{cases}
  \vspace{1mm}
 |(F_7|\dde)_{\dot{H}^s}| \leq \frac{\nu}{52} \|\n \dde\|_{\dot{H}^s}^2 +\frac{C}{\nu^3} \|\tu\|_{L^2(\R^2)}^2 \|\n_h \tu\|_{L^2(\R^2)}^2 \|\dde\|_{\dot{H}^s}^2,\\
 \vspace{1mm}
 |(G_5|\de)_{\dot{H}^s}| \leq \frac{\nu'}{56} \|\n \de\|_{\dot{H}^s}^2 +\frac{C}{(\nu')^3} \|\tu\|_{L^2}^2 \|\n_h \tu\|_{L^2}^2 \|\de\|_{\dot{H}^s}^2.
 \end{cases}
 \label{Estim:ext:4}
\end{equation}
Next, we estimate $F_{12}$ with a similar method as we used for $F_7$ (the coefficients in the final Young estimates being $(4,4,2)$):
\begin{multline}
 |(F_{12}|\dde)_{\dot{H}^s}| \leq \|\tb\cdot \n \de\|_{\dot{H}^{s-1}} \|\dde\|_{\dot{H}^{s+1}} \leq C \|\tb\|_{\dot{H}^\frac12(\R^2)} \|\n \de\|_{\dot{H}^{s-\frac12}(\R^3)} \|\n \dde\|_{\dot{H}^s} \\
 \leq C \|\tb\|_{L^2}^\frac12 \|\n_h \tb\|_{L^2}^\frac12 \|\de\|_{\dot{H}^s}^\frac12 \|\de\|_{\dot{H}^{s+1}}^\frac12 \|\n \dde\|_{\dot{H}^s}\\
 \leq \frac{\nu}{52} \|\n \dde\|_{\dot{H}^s}^2 +\frac{\nu'}{56} \|\n \de\|_{\dot{H}^s}^2 +\frac{C}{\nu^2 \nu'} \|\tb\|_{L^2}^2 \|\n_h \tb\|_{L^2}^2 \|\de\|_{\dot{H}^s}^2.
 \label{Estim:ext:5}
\end{multline}
Similarly, we get
\begin{equation}
 |(G_{11}|\dde)_{\dot{H}^s}| \leq \frac{\nu}{52} \|\n \dde\|_{\dot{H}^s}^2 +\frac{\nu'}{56} \|\n \de\|_{\dot{H}^s}^2 +\frac{C}{\nu (\nu')^2} \|\tb\|_{L^2}^2 \|\n_h \tb\|_{L^2}^2 \|\dde\|_{\dot{H}^s}^2.
 \label{Estim:ext:6}
\end{equation}
We skip details for $F_6,G_9$ which are dealt with as in \cite{FCRF} (we refer to the estimates of the term $G_7$ from Section 3.2 therein, and to \eqref{DefLaniso} for definition of anisotropic Lebesgue spaces) and we get that:
\begin{equation}
 \begin{cases}
  \vspace{1mm}
|(F_6|\dde)_{\dot{H}^s}| \leq \frac{\nu}{52} \|\n \dde\|_{\dot{H}^s}^2 +\frac{C}{\nu^\frac{s}{1-s}} \|\n_h \tu\|_{L^2}^2 \|\dde\|_{\dot{H}^s}^2 +C\|\We\|_{L_{h,v}^{\infty,2}}^2 \|\n_h \tu\|_{L^2}^{2s},\\
|(G_9|\de)_{\dot{H}^s}| \leq \frac{\nu'}{56} \|\n \de\|_{\dot{H}^s}^2 +\frac{C}{(\nu')^\frac{s}{1-s}} \|\n_h \tb\|_{L^2}^2 \|\de\|_{\dot{H}^s}^2 +C\|\We\|_{L_{h,v}^{\infty,2}}^2 \|\n_h \tb\|_{L^2}^{2s}.
 \end{cases}
 \label{Estim:ext:7}
\end{equation}
We are now left with the new terms. Some are easily estimated using the previous arguments because they feature $\ce$, which is very regular:
\begin{multline}
 |(F_{11}|\dde)_{\dot{H}^s}| \leq \|\de\cdot \n \ce\|_{\dot{H}^{s-1}} \|\dde\|_{\dot{H}^{s+1}} \leq C \|\de\|_{\dot{H}^s} \|\n \ce\|_{\dot{H}^\frac12} \|\n \dde\|_{\dot{H}^s}\\
 \leq \frac{\nu}{52} \|\n \dde\|_{\dot{H}^s}^2 +\frac{C}{\nu} \|\n \ce\|_{\dot{H}^\frac12}^2 \|\de\|_{\dot{H}^s}^2.
 \label{Estim:ext:8}
\end{multline}
Similarly, we obtain that:
\begin{equation}
 |(G_8|\de)_{\dot{H}^s}| \leq \frac{\nu'}{56} \|\n \de\|_{\dot{H}^s}^2 +\frac{C}{\nu'} \|\n \ce\|_{\dot{H}^\frac12}^2 \|\dde\|_{\dot{H}^s}^2.
 \label{Estim:ext:9}
\end{equation}
Using the classical Sobolev product laws (with $(s_1,s_2)=(1,s-\frac12)$) instead of $2D\times 3D$ laws in the method used for $F_{12}$, we obtain that:
\begin{equation}
 \begin{cases}
 \vspace{1mm}
  |(F_{13}|\dde)_{\dot{H}^s}| \leq \frac{\nu}{52} \|\n \dde\|_{\dot{H}^s}^2 +\frac{\nu'}{56} \|\n \de\|_{\dot{H}^s}^2 +\frac{C}{\nu^2 \nu'} \|\ce\|_{\dot{H}^\frac12}^2 \|\n \ce\|_{\dot{H}^\frac12}^2 \|\de\|_{\dot{H}^s}^2,\\
  |(G_{12}|\de)_{\dot{H}^s}| \leq \frac{\nu}{52} \|\n \dde\|_{\dot{H}^s}^2 +\frac{\nu'}{56} \|\n \de\|_{\dot{H}^s}^2 +\frac{C}{\nu (\nu')^2} \|\ce\|_{\dot{H}^\frac12}^2 \|\n \ce\|_{\dot{H}^\frac12}^2 \|\dde\|_{\dot{H}^s}^2.
 \end{cases}
 \label{Estim:ext:10}
\end{equation}
Slightly adapting the arguments used for $F_6,G_9$ (using Young estimates with $(6,3)$, followed by Sobolev and interpolation estimates, and finally using an algebraic Young estimates with coefficients $(2,\frac2{1-s}, \frac2{s})$) we obtain:
\begin{multline}
 |(G_{10}|\de)_{\dot{H}^s}| \leq \|\We\cdot \n \ce\|_{L^2} \|\de\|_{\dot{H}^{2s}} \leq C \|\We\|_{L^6} \|\n \ce\|_{\dot{H}^\frac12} \|\de\|_{\dot{H}^s}^{1-s} \|\de\|_{\dot{H}^{s+1}}^s\\
 \leq C \left(\|\We\|_{L^6} \|\n \ce\|_{\dot{H}^\frac12}^s\right)\left(\|\n \ce\|_{\dot{H}^\frac12} \|\de\|_{\dot{H}^s}\right)^{1-s} \|\de\|_{\dot{H}^{s+1}}^s\\
 \leq \frac{\nu'}{56} \|\n \de\|_{\dot{H}^s}^2 +\frac{C}{(\nu')^\frac{s}{1-s}} \|\n \ce\|_{\dot{H}^\frac12}^2 \|\de\|_{\dot{H}^s}^2 +C\|\We\|_{L^6}^2 \|\n \ce\|_{\dot{H}^\frac12}^{2s}.
 \label{Estim:ext:11}
\end{multline}
An analogous of $F_8$ was already present in \cite{FCRF} but the estimates that we obtained in this article featured $2s-1$ as an exponent, which required $s\geq \frac12$. In order to obtain estimates in the case $s<\frac12$ we had to use another method, which produced an estimate of $\We$ with time integration index $\frac2{2-s}$. As outlined in Remark \ref{rem:chgtFextk}, contrary to its counterpart $(LRF_\ee)$ from \cite{FCRF}, System \eqref{SystWe} now features several external force terms with different regularities, and we will have to use the modified Strichartz estimates from Section \ref{Sect:AnisoStri}, which require time integration index greater than $2$. As $\frac2{2-s}\geq 2$ is false when $s$ is close to $\frac12$ we need to change our estimates. Adapting the arguments used for $F_6,F_9$ and $G_{10}$, and aiming to use the following Sobolev injection:
$$
\dot{H}^{1-s}(\R^2) \hookrightarrow L^{\frac2{s}}(\R^2),
$$
we have (also using interpolation and the Young estimates with coefficients $(2,\frac2{1-s}, \frac2{s})$):
\begin{multline}
 |(F_8|\dde)_{\dot{H}^s}| \leq \|\tu \cdot \n \We\|_{L^2} \|\dde\|_{\dot{H}^{2s}} \leq C \|\tu\|_{L^\frac2{s}(\R^2)} \|\n \We\|_{L_{h,v}^{\frac2{1-s},2}(\R^3)} \|\dde\|_{\dot{H}^s}^{1-s} \|\dde\|_{\dot{H}^{s+1}}^s\\
 \leq C \left(\|\tu\|_{L^2(\R^2)}^s  \|\n \We\|_{L_{h,v}^{\frac2{1-s},2}(\R^3)}\right) \left(\|\tu\|_{\dot{H}^1(\R^2)} \|\dde\|_{\dot{H}^s}\right)^{1-s} \|\dde\|_{\dot{H}^{s+1}}^s\\
 \leq \frac{\nu}{52} \|\n \dde\|_{\dot{H}^s}^2 +\frac{C}{\nu^\frac{s}{1-s}} \|\n_h \tu\|_{L^2}^2 \|\dde\|_{\dot{H}^s}^2 +C\|\tu\|_{L^2}^{2s} \|\n \We\|_{L_{h,v}^{\frac2{1-s},2}}^2.
 \label{Estim:ext:12}
\end{multline}
\begin{rem}
 \sl{This new estimates is valid independantly of the positions of $s$ and $\frac12$ and would simplify a little the proofs from \cite{FCRF} in the case of rotating fluids.}
\end{rem}
Similarly, we have:
\begin{equation}
 |(G_{13}|\de)_{\dot{H}^s}| \leq \frac{\nu'}{56} \|\n \de\|_{\dot{H}^s}^2 +\frac{C}{(\nu')^\frac{s}{1-s}} \|\n_h \tb\|_{L^2}^2 \|\de\|_{\dot{H}^s}^2 +C\|\tb\|_{L^2}^{2s} \|\n \We\|_{L_{h,v}^{\frac2{1-s},2}}^2.
 \label{Estim:ext:13}
\end{equation}

\subsubsection{Estimating $G_{14}$}

We will find how to correctly estimate the last term $G_{14}$ under the guidance of several constraints:
\begin{itemize}
 \item As outlined in Remark \ref{rem:chgtFextk}, we need to produce Strichartz estimates featuring a time index of integration $p\geq 2$.
 \item We need enough integrability on $\|\ce\|_{\dot{H}^{\sigma}}$ to be able to evacuate $\|\de\|_{\dot{H}^s}^2$ through the Gronwall lemma.
\end{itemize}
First let us state the following result which, in the spirit of Proposition 2 from \cite{FCRF}, details how time-integrable is $\|\ce(t)\|_{\dot{H}^{\sigma}}$ for a given $\sigma$.
\begin{prop}
 \sl{Let $\delta>0$, $\tu,\tb$ satisfying \eqref{estimtutb}, and $\ce$ be the unique global solution of \eqref{Mag2} corresponding to the regular initial $\coe\in H^{\frac12+\delta}$. Then for any $\sigma\in[0,\frac32+\delta]$, we have that $t\mapsto \|\ce(t)\|_{\dot{H}^{\sigma}}$ is in $L^p(\R_+)$ for any $p$ in:
\begin{equation}
  \begin{cases}
  \vspace{1mm}
  [\frac2{\sigma}, +\infty] & \mbox{if }\sigma\in[0,\frac12+\delta],\\
    \vspace{1mm}
   [\frac2{\sigma}, \frac2{\sigma-(\frac12+\delta)}] & \mbox{if }\sigma\in]\frac12+\delta,1],\\
   [2, \frac2{\sigma-(\frac12+\delta)}] & \mbox{if }\sigma\in]1,\frac32+\delta].
  \end{cases}
 \end{equation}
 Moreover in any case, we have for some constant $C$ (depending on $\nu,\nu',\|\tu_0\|_{L^2}, \|\tb_0\|_{L^2}$)
 $$
 \|\ce\|_{L^p\dot{H}^\sigma} \leq C \|\coe\|_{H^{\frac12+\delta}}.
 $$
 }
 \label{Propestimce}
\end{prop}
\textbf{Proof: }thanks to \eqref{estimcHs} and by interpolation, there exists a constant $C=C(\nu,\nu', \|\tuo\|_{L^2}, \|\tbo\|_{L^2})>0$ such that for any $p\in[2,\infty]$, any $s\in[0,\frac12+\delta]$,
$$
\|\ce\|_{L^p \dot{H}^{s+\frac2{p}}} \leq C \|\co\|_{\dot{H}^s} \leq C \|\co\|_{H^{\frac12+\delta}}.
$$
We note that with these assumptions, $\sigma \overset{def}{=} s+\frac2{p} \in [0,\frac32+\delta]$. For any $\sigma \in [0,\frac32+\delta]$ and any $p\in[2,\infty]$ there exists some $s\in[0,\frac12+\delta]$ such that $\sigma=s+\frac2{p}$ if and only if
$$
\frac2{p}\in [\sigma-(\frac12+\delta), \sigma].
$$
From this we have several cases:
\begin{itemize}
 \item If $\sigma\in [0, \frac12+\delta]$, the previous condition is equivalent to  $\frac2{p}\in [0, \sigma]$ that is $p\in[\frac2{\sigma},\infty]$.
 \item If $\sigma\in ]\frac12+\delta, \frac32 +\delta]$, the previous condition is equivalent to $p\in[\frac2{\sigma}, \frac2{\sigma-(\frac12+\delta)}]$. We have to consider the lower bound of this interval keeping in mind that $p\in[2,\infty]$: $\frac2{\sigma}\leq 2 \Leftrightarrow \sigma \geq 1$, which gives the last two cases from the proposition. $\blacksquare$
\end{itemize}
Now we can turn to the estimate for $G_{14}$: let us begin as follows: for any $q_1,q_2\in[2,\infty]$ such that $\frac1{q_1}+\frac1{q_2}=\frac12$, we have
$$
 |(G_{14}|\de)_{\dot{H}^s}| \leq \|\ce\cdot \n \We\|_{L^2} \|\de\|_{\dot{H}^{2s}} \leq C \|\ce\|_{L^{q_1}} \|\n \We\|_{L^{q_2}} \|\de\|_{\dot{H}^s}^{1-s} \|\de\|_{\dot{H}^{s+1}}^s.
$$
Thanks to the Sobolev injections in $\R^3$, if we ask $q_1<\infty$ and put $\sigma'=3(\frac12-\frac1{q_1})$, then $\sigma'\in[0,\frac32[$ and $q_2=\frac3{\sigma'}$ so that (using for the third line the Young estimates with coefficients $(2,\frac2{1-s}, \frac2{s})$) for some constant $C=C_{s,\sigma}$,
\begin{multline}
|(G_{14}|\de)_{\dot{H}^s}| \leq C \|\ce\|_{\dot{H}^{\sigma'}} \|\n \We\|_{L^{\frac3{\sigma'}}} \|\de\|_{\dot{H}^s}^{1-s} \|\de\|_{\dot{H}^{s+1}}^s\\
\leq C \left(\|\ce\|_{\dot{H}^{\sigma'}}^s \|\n \We\|_{L^{\frac3{\sigma'}}}\right) \left(\|\ce\|_{\dot{H}^{\sigma'}} \|\de\|_{\dot{H}^s}\right)^{1-s} \|\de\|_{\dot{H}^{s+1}}^s\\
\leq \frac{\nu'}{56} \|\n \de\|_{\dot{H}^s}^2 +\frac{C}{(\nu')^\frac{s}{1-s}} \|\ce\|_{\dot{H}^{\sigma'}}^2 \|\de\|_{\dot{H}^s}^2 +C\|\ce\|_{\dot{H}^{\sigma'}}^{2s} \|\n \We\|_{L^\frac3{\sigma'}}^2.
\end{multline}
Notice that this estimate is useful when $\sigma'$ is such that $\|\ce\|_{\dot{H}^{\sigma'}}\in L^2$ that is (thanks to Proposition \ref{Propestimce}) when $\sigma'\in[1,\frac32[$. In this case we will have to bound the following integral according to:
$$
\int_0^t \|\ce\|_{\dot{H}^{\sigma'}}^{2s} \|\n \We\|_{\frac3{\sigma'}}^2 d\tau \leq \left(\int_0^t \|\ce\|_{\dot{H}^{\sigma'}}^{2p_1s}d\tau\right)^\frac1{p_1} \left(\int_0^t \|\n \We\|_{\frac3{\sigma'}}^{2p_2} d\tau\right)^\frac1{p_2},
$$
where $\frac1{p_1}+\frac1{p_2}=1$. If we wish to minimize $p_2$ we need to maximize $p_1$ so that thanks to Proposition \ref{Propestimce} we will choose $2p_1 s=\frac2{\sigma'-(\frac12+\delta)}$ and we end-up, for some $\sigma'\in[1,\frac32[$ to specify, with:
\begin{equation}
 \int_0^t \|\ce\|_{\dot{H}^{\sigma'}}^{2s} \|\n \We\|_{\frac3{\sigma'}}^2 d\tau \leq \|\ce\|_{L^\frac2{\sigma'-(\frac12+\delta)} \dot{H}^{\sigma'}}^{2s} \|\n \We\|_{L^\frac2{1-s(\sigma'-(\frac12+\delta))} L^\frac3{\sigma'}}^2.
 \label{Choixp1}
\end{equation}
Thanks to Proposition \ref{EstimStri}, the regularity index that we have to consider applying the first point for $(d,p,r)=(1,\frac2{\sigma'-(\frac12+\delta)},\frac3{\sigma'})$ (and for some $\theta\in]0,1]$) is
\begin{equation}
 \sigma_1=\frac32-\sigma'+s\big(\sigma'-(\frac12+\delta)\big)+\theta(1-\frac23\sigma').
 \label{sigma1}
\end{equation}
Due to the regularities of the external force terms in \eqref{Mag2}, we expect to use the bound $\ce\cdot\n \ce$ in $L^1 \dot{H}^{\frac12+\delta}$ (we refer to the Appendix for details). As $s$ is close to $\frac12$ and $\theta$ will be small, taking $(s,\theta)=(\frac12,0)$ we get that $\sigma_1=\frac54-\frac{\sigma'}2-\frac{\delta}2$ so that
$$
\sigma_1\sim \frac12+\delta \quad \Longleftrightarrow \quad \sigma'\sim \frac32-3\delta,
$$
which suggests us to choose $\sigma'=\frac32-\aa\delta$. Injecting this in \eqref{sigma1} we obtain:
$$
  \sigma_1=\aa \delta +s\Big(1-(\aa+1)\delta\Big)+\frac{2\aa \delta}3 \theta,
$$
and we can find some $\theta\in]0,1]$ such that $\sigma_1=\frac12+\delta$ if and only if:
$$
 \aa \delta +s\Big(1-(\aa+1)\delta\Big)< \frac12+\delta \leq \frac{5\aa \delta}3 +s\Big(1-(\aa+1)\delta\Big),
$$
which is possible for any $s\in[\frac12-\eta \delta, \frac12+\eta \delta]$ if, and only if
$$
\aa \delta +(\frac12+\eta \delta)\Big(1-(\aa+1)\delta\Big) <\frac12+\delta \leq \frac{5\aa \delta}3 +(\frac12-\eta \delta)\Big(1-(\aa+1)\delta\Big),
$$
that is if, and only if
$$
\frac{9+6\eta(1-\delta)}{7+6\eta \delta}\leq \aa < \frac{3-2\eta(1-\delta)}{1-2\eta \delta},
$$
which in turn is possible when $\eta(20-56\delta)<12$, and true as soon as $\eta<\frac35$. We could optimize the choice of $\eta$ but, to simplify the presentation, from now on we will fix $\eta=\frac{\eta_0}{10}$. Moreover, we need $\aa$ to be as small as possible (we refer to the appendix for an explaination), so we also fix
\begin{equation}
 \aa= \aa_0=\frac{9+\frac35 \eta_0(1-\delta)}{7+\frac35 \eta_0\delta},
\label{Choixalpha}
\end{equation}
finally leading to the last estimates
\begin{equation}
 |(G_{14}|\de)_{\dot{H}^s}| \leq \frac{\nu'}{56} \|\n \de\|_{\dot{H}^s}^2 +\frac{C}{\nu^\frac{s}{1-s}} \|\ce\|_{\dot{H}^{\frac32-\aa_0\delta}}^2 \|\de\|_{\dot{H}^s}^2 +C\|\ce\|_{\dot{H}^{\frac32-\aa_0\delta}}^{2s} \|\n \We\|_{L^\frac3{\frac32-\aa_0\delta}}^2.
 \label{Estim:ext:14}
\end{equation}

\subsubsection{End of the proof}

Injecting \eqref{Estim:ext:1} to \eqref{Estim:ext:13} and \eqref{Estim:ext:14} into \eqref{Systddde}, we obtain that for any $t\leq T_\ee<T_\ee^*$, we have
\begin{multline}
 \frac12 \frac{d}{dt}\left(\|\dde(t)\|_{\dot{H}^s}^2 +\|\de(t)\|_{\dot{H}^s}^2\right) +\frac34\nu \|\n \dde(t)\|_{\dot{H}^s}^2 +\frac34\nu' \|\n \de(t)\|_{\dot{H}^s}^2\\
 \leq C\left(\|\dde\|_{\dot{H}^\frac12} +\|\de\|_{\dot{H}^\frac12}\right) \left(\|\n \dde\|_{\dot{H}^s}^2 +\|\n \de\|_{\dot{H}^s}^2\right) +J_1(t) \left(\|\dde(t)\|_{\dot{H}^s}^2 +\|\de(t)\|_{\dot{H}^s}^2\right) +J_2(t),
 \label{estim:Hsgen}
\end{multline}
where
\begin{multline}
 J_1(t) =C_{\nu,\nu',s} \Big\{(1+\|\tu\|_{L^2(\R^2)}^2) \|\n_h \tu\|_{L^2(\R^2)}^2 +(1+\|\tb\|_{L^2(\R^2)}^2) \|\n_h \tb\|_{L^2(\R^2)}^2\\
 +(1+\|\ce\|_{\dot{H}^\frac12}^2) \|\n \ce\|_{\dot{H}^\frac12}^2 +\|\n \ce\|_{\dot{H}^{\frac12-\aa_0\delta}}^2 +\|\n \We\|_{L^3}^2 +\|\We\|_{L^6}^4 +\|\We\|_{L^6}^\frac2{1-s}\Big\},
\end{multline}
and
\begin{multline}
 J_2(t)=C_{\nu,\nu',s} \Big[\|\n \We\|_{L^3}^2 +\|\We\|_{L_{h,v}^{\infty,2}}^2 (\|\n_h \tu\|_{L^2}^{2s} +\|\n_h \tb\|_{L^2}^{2s}) +\|\We\|_{L^6}^2 \|\n \ce\|_{\dot{H}^\frac12}^{2s}\\
 +(\|\tu\|_{L^2}^{2s} +\|\tb\|_{L^2}^{2s}) \|\n \We\|_{L_{h,v}^{\frac2{1-s},2}}^2 +\|\n \ce\|_{\dot{H}^{\frac12-\aa_0\delta}}^{2s} \|\n \We\|_{L^\frac3{\frac32-\aa_0\delta}}^2\Big].
\end{multline}
Thanks to the fact that $t\leq T_\ee$ (we refer to \eqref{DefTe}), we can absorb the first term from the right-hand-side of \eqref{estim:Hsgen}, so that integrating and using the Gronwall estimates, we obtain that (thanks to \eqref{Choixp1} and the choice in \eqref{Choixalpha}) for any $t\leq T_\ee$ (we recall the initial data is zero),
\begin{multline}
 \|\dde(t)\|_{\dot{H}^s}^2 +\|\de(t)\|_{\dot{H}^s}^2 +\nu \int_0^t \|\n \dde(\tau)\|_{\dot{H}^s}^2d\tau +\nu' \int_0^t\|\n \de(\tau)\|_{\dot{H}^s}^2d\tau\\
 \leq \left(\int_0^t J_2(t')dt'\right) e^{\int_0^t J_1(t') dt'},
 \end{multline}
 where (thanks to \eqref{estimtutb} and \eqref{estimcHs}) we have
  \begin{multline}
 \int_0^t J_1(t') dt' \leq C_{\nu,\nu',s} \Big\{(1+\|\tu\|_{L_t^\infty L^2(\R^2)}^2) \|\n_h \tu\|_{L_t^2 L^2(\R^2)}^2 +(1+\|\tb\|_{L_t^\infty L^2(\R^2)}^2) \|\n_h \tb\|_{L_t^2 L^2(\R^2)}^2\\
 +(1+\|\ce\|_{L_t^\infty \dot{H}^\frac12}^2) \|\n \ce\|_{L_t^2 \dot{H}^\frac12}^2 +\|\n \ce\|_{L_t^2 \dot{H}^{\frac12-\aa_0\delta}}^2 +\|\n \We\|_{L_t^2 L^3}^2 +\|\We\|_{L_t^4 L^6}^4 +\|\We\|_{L_t^\frac2{1-s} L^6}^\frac2{1-s}\Big\}\\
\leq \Do \Big\{1+\|\coe\|_{H^{\frac12+\delta}}^4 +\|\n \We\|_{L_t^2 L^3}^2 +\|\We\|_{L_t^4 L^6}^4 +\|\We\|_{L_t^\frac2{1-s} L^6}^\frac2{1-s}\Big\}
\end{multline}
and (thanks to Proposition \ref{Propestimce})
\begin{multline}
 \int_0^t J_2(t')dt' \leq C_{\nu,\nu',s} \Big[\|\n \We\|_{L_t^2 L^3}^2 +\|\We\|_{L_t^\frac2{1-s} L_{h,v}^{\infty,2}}^2 (\|\n_h \tu\|_{L_t^2 L^2}^{2s} +\|\n_h \tb\|_{L_t^2 L^2}^{2s})\\
 +\|\We\|_{L_t^\frac2{1-s} L^6}^2 \|\n \ce\|_{L_t^2 \dot{H}^\frac12}^{2s} +(\|\tu\|_{L_t^\infty L^2}^{2s} +\|\tb\|_{L_t^\infty L^2}^{2s}) \|\n \We\|_{L_t^2 L_{h,v}^{\frac2{1-s},2}}^2\\
 +\|\ce\|_{L_t^\frac2{1-(\aa_0+1)\delta}\dot{H}^{\frac32-\aa_0\delta}}^{2s} \|\n \We\|_{L_t ^\frac2{1-s(1-(\aa_0+1)\delta)}L^\frac3{\frac32-\aa_0\delta}}^2\Big]\\
 \leq \Do \Big[\|\n \We\|_{L_t^2 L^3}^2 +\|\We\|_{L_t^\frac2{1-s} L_{h,v}^{\infty,2}}^2 +\|\n \We\|_{L_t^2 L_{h,v}^{\frac2{1-s},2}}^2 +\|\coe\|_{H^{\frac12+\delta}}^{2s} \|\We\|_{L_t^\frac2{1-s} L^6}^2\\
 +\|\coe\|_{H^{\frac12+\delta}}^{2s} \|\n \We\|_{L_t ^\frac2{1-s(1-(\aa_0+1)\delta)}L^\frac3{\frac32-\aa_0\delta}}^2\big)\Big],
 \end{multline}
for some constant $\Do=\Do(\nu, \nu',\Co,s)\geq 1$.
We finally obtain that for any $t\leq T_\ee$,
\begin{multline}
 \|\dde(t)\|_{\dot{H}^s}^2 +\|\de(t)\|_{\dot{H}^s}^2 +\nu \int_0^t \|\n \dde(\tau)\|_{\dot{H}^s}^2d\tau +\nu' \int_0^t\|\n \de(\tau)\|_{\dot{H}^s}^2d\tau\\
 \leq \Do \Big[\|\n \We\|_{L_t^2 L^3}^2 +\|\We\|_{L_t^\frac2{1-s} L_{h,v}^{\infty,2}}^2 +\|\n \We\|_{L_t^2 L_{h,v}^{\frac2{1-s},2}}^2
 +\|\coe\|_{H^{\frac12+\delta}}^{2s} \|\We\|_{L_t^\frac2{1-s} L^6}^2 \\
 +\|\coe\|_{H^{\frac12+\delta}}^{2s} \|\n \We\|_{L_t ^\frac2{1-s(1-(\aa_0+1)\delta)}L^\frac3{\frac32-\aa_0\delta}}^2\Big]\\
 \times \exp \Do \Big\{1+\|\coe\|_{H^{\frac12+\delta}}^4 +\|\n \We\|_{L_t^2 L^3}^2 +\|\We\|_{L_t^4 L^6}^4 +\|\We\|_{L_t^\frac2{1-s} L^6}^\frac2{1-s}\Big\},
 \label{Estimapriori3}
 \end{multline}

\subsection{Use of Strichartz estimates}

The isotropic and anisotropic Strichartz estimates allow us to state the following result (whose proof is postponed to Sections \ref{Sect:IsoStri} and \ref{Sect:AnisoStri} in the appendix):

\begin{prop}
 \sl{Under the notations from Theorem \ref{ThCVStrong}, there exist positive constants $\ee_0,\Eo$ (both of them depending on $\nu, \nu', \Co, \delta,\eta_0,\mu$) such that if $\delta\leq \frac16$, for any $\ee\in]0,\ee_0]$ and any $s\in[\frac12-\frac{\eta_0}{10}\delta, \frac12+\frac{\eta_0}{10}\delta]$, we have:
 \begin{equation}
  \begin{cases}
  \vspace{2mm}
   \|\n \We\|_{L^2 L^3} +\|\We\|_{L^4 L^6} \leq \Eo \ee^{\eta_0 \delta},\\
 \|\We\|_{L^2 L^\infty} \leq \Eo \ee^{\frac{\eta_0}2 \delta},\\
  \|\n \We\|_{L_t ^\frac2{1-s(1-(\aa_0+1)\delta)}L^\frac3{\frac32-\aa_0\delta}} \leq \Eo \ee^{\frac{\eta_0}{16} \delta},\\
    \|\We\|_{L^\frac2{1-s}L^6} +\|\We\|_{L^\frac2{1-s} L_{h,v}^{\infty,2}} +\|\n \We\|_{L^2 L_{h,v}^{\frac2{1-s},2}} \leq \Eo \ee^{\frac9{10}\eta_0 \delta}.
  \end{cases}
\label{EstimStriStrong}
 \end{equation}
 }
\label{PropestimStriStrong}
\end{prop}
Plugging this into \eqref{Estimapriori3}, if in addition $\ee_0$ is so small that for any $\ee\in]0,\ee_0]$:
\begin{multline}
\|\n \We\|_{L_t^2 L^3}^2 +\|\We\|_{L_t^4 L^6}^4 +\|\We\|_{L_t^\frac2{1-s} L^6}^\frac2{1-s} \leq 2\Eo^2 \ee^{2\eta_0 \delta} +(\Eo \ee^{\frac{9\eta_0}{10}\delta})^\frac2{1-s}\\ \leq 2\Eo^2 \ee^{2\eta_0 \delta} +(\Eo \ee^{\frac{9\eta_0}{10}\delta})^\frac4{1+\frac{\eta_0}5 \delta} \leq 1,
 \label{Condeps1}
\end{multline}
we get that for all $t\leq T_\ee$ and $s\in[\frac12-\frac{\eta_0}{10}\delta, \frac12+\frac{\eta_0}{10}\delta]$,
\begin{multline}
 \|\dde(t)\|_{\dot{H}^s}^2 +\|\de(t)\|_{\dot{H}^s}^2 +\nu \int_0^t \|\n \dde(\tau)\|_{\dot{H}^s}^2d\tau +\nu' \int_0^t\|\n \de(\tau)\|_{\dot{H}^s}^2d\tau\\
 \leq 2\Do \Eo^2\Big[\ee^{\frac{9\eta_0}5 \delta} +\Big(\Ko |\ln \ee|\Big)^\frac{s}2 \ee^{\frac{\eta_0}8 \delta}\Big]
 \exp \Do \Big\{1+\Ko |\ln \ee| +2\Eo^2 \ee^{2\eta_0 \delta} +(\Eo \ee^{\frac{9\eta_0}{10}\delta})^\frac2{1-s} \Big\}\\
 \leq 2\Do \Eo^2\Big[\ee^{\frac{9\eta_0}5 \delta} +\Big(\Ko |\ln \ee|\Big)^\frac{s}2 \ee^{\frac{\eta_0}8 \delta}\Big]e^{\Do \{2-\Ko\ln \ee\}}.
 \end{multline}
If we assume that $\ee_0$ is so small that for any $\ee\in]0,\ee_0]$,
\begin{equation}
 (\Ko|\ln \ee|)^\frac12 \leq \ee^{-\frac{\eta_0 \delta}{96}},
 \label{condeps3}
\end{equation}
then if we choose $\Ko$ according to
\begin{equation}
\Do \Ko\overset{def}{=} \frac{\eta_0\delta}{144},
\label{ChoixKo}
\end{equation}
recalling that $\delta\leq \frac16$, we obtain
 \begin{multline}
  \|\dde(t)\|_{\dot{H}^s}^2 +\|\de(t)\|_{\dot{H}^s}^2 +\nu \int_0^t \|\n \dde(\tau)\|_{\dot{H}^s}^2d\tau +\nu' \int_0^t\|\n \de(\tau)\|_{\dot{H}^s}^2d\tau\\
  \leq 4\Do \Eo^2 \ee^{\eta_0 \delta\Big[\frac18-(\frac12+\frac{\eta_0}{10}\delta)\frac1{96})\Big]} e^{2\Do} \ee^{-\Do \Ko} \leq 4\Do e^{2\Do}\Eo^2 \ee^{\eta_0 \delta\Big[\frac18-(\frac12+\delta)\frac1{96})-\frac1{144}\Big]}\\
  \leq \Fo \ee^{\eta_0 \delta (\frac18 -\frac1{72})} =\Fo \ee^{\frac{\eta_0}9 \delta}.
 \label{EstimaprioriS}
 \end{multline}

\subsection{Conclusion of the proof}

First assuming $\voe \in \dot{H}^\frac12 \cap \dot{H}^{\frac12+\delta}$, using \eqref{EstimaprioriS} in the case $s=\frac12$ we obtain that for any $t\leq T_\ee$ (we refer to \eqref{DefTe} for the definition of $T_\ee$):
$$
 \|\dde(t)\|_{\dot{H}^\frac12}^2 +\|\de(t)\|_{\dot{H}^\frac12}^2 +\nu \int_0^t \|\n \dde(\tau)\|_{\dot{H}^\frac12}^2d\tau +\nu' \int_0^t\|\n \de(\tau)\|_{\dot{H}^\frac12}^2d\tau \leq \Fo \ee^{\frac{\eta_0}9 \delta} \leq \frac1{16C} \min(\nu, \nu'),
$$
for any $\ee\in]0,\ee_0]$ if $\ee_0$ is small enough. This contradicts the definition of $T_\ee$, and allows us to get (by the continuation criterion) that $T_\ee=T_\ee^*=\infty$. Then thanks to the propagation of regularity we are able to obtain that for any $s\in[\frac12, \frac12+\frac{\eta_0}{10}\delta]$ and $t\geq 0$,
$$
 \|\dde(t)\|_{\dot{H}^s}^2 +\|\de(t)\|_{\dot{H}^s}^2 +\nu \int_0^t \|\n \dde(\tau)\|_{\dot{H}^s}^2d\tau +\nu' \int_0^t\|\n \de(\tau)\|_{\dot{H}^s}^2d\tau \leq \Fo \ee^{\frac{\eta_0}9 \delta},
$$
which gives \eqref{EstimTHStrongA}. If we now assume that $\voe \in \dot{H}^{\frac12-\delta} \cap \dot{H}^{\frac12+\delta}$, then the previous estimates is true for any $s\in[\frac12-\frac{\eta_0}{10}\delta, \frac12+\frac{\eta_0}{10}\delta]$ (using that $\delta\leq \frac16$, so that (using once more Proposition \ref{estimBsHs}) we get:
$$
\|\De -\We\|_{L^2 L^\infty} \leq \Fo \ee^{\frac{\eta_0}{18} \delta}.
$$
Combining this with the second line from \eqref{EstimStriStrong} finally gives \eqref{EstimTHStrongB}. $\blacksquare$

\section{Appendix}

\subsection{Notations for general Sobolev and Besov spaces}

We refer to \cite{Dbook} for a complete presentation of Sobolev and Besov spaces through the Littlewood-Paley decompositions. We will only recall some basic fundamental results and give more details about refined results involving products of functions depending on $x=(x_1,x_2,x_3)$ by functions only depending on $x_h=(x_1,x_2)$.\\
Assume $\chi:\R_+\rightarrow \R$ is a smooth function supported in the ball $[0,\frac43]$, which is equal to 1 in a neighborhood of $[0,\frac43]$ and nonincreasing over $\R_+$. Setting $\varphi(r)=\chi(r/2)-\chi(r)$, we observe that $\varphi$ is compactly supported in the interval $\cC=[\frac34, \frac83]$ and we define the homogeneous dyadic blocks in $\R^d$ as follows: for all $j\in \Z$ and any $u$ depending on $x=(x_1,...,x_d)$, and recalling that $\hat{k(D)u}(\xi)=k(\xi) \hat{u} (\xi)$, we define
\begin{equation}
 \begin{cases}
\vspace{2mm}
\ddj u:= \varphi(2^{-j}|D|)u =2^{jd} m(2^j.)* u, \quad \mbox{with } m(x)=\cF^{-1} \big(\varphi(|\xi|)\big),\\
\dot{S}_j u=\Sum_{l\leq j-1} \ddl u=\chi(2^{-j}|D|)u =2^{jd} w(2^j.)* u, \quad \mbox{with } w(x)=\cF^{-1} \big(\chi(|\xi|)\big),\\
\end{cases}
\label{DefBesov}
\end{equation}
We can now define the homogeneous Besov norms and spaces:
\begin{defi}
\sl{For $s\in\R$ and $1\leq p,r\leq\infty,$ we set
$$
\dot{B}_{p,r}^s(\R^d)=\left\{u \in \cS'(\R^d), \mbox{ with } \underset{j\rightarrow -\infty}{\mbox{lim}} \|\dot{S}_j u\|_{L^\infty}=0 \mbox{ and } \|u\|_{\dot{B}_{p,r}^s} \overset{\mbox{def}}{=} \|\big(2^{qs} \|\ddq u\|_{L^p} \big)_{q \in \Z}\|_{\ell^r} < \infty\right\}.
$$
}
\end{defi}
\begin{rem}
 \sl{The inhomogeneous Sobolev space injects into its homogeneous couterpart as we have:
 $$
 H^s(\R^d)=L^2(\R^d)\cap \dot{H}^s(\R^d)=L^2(\R^d)\cap \dot{B}_{2,2}^s(\R^d).
 $$}
\end{rem}

As in \cite{FCPAA, FCcompl, FCRF, FCStratif1, FCStratif2} we will use spaces that are slight modifications of $L_t^p \dot{B}_{q,r}^s$: namely the Chemin-Lerner time-space Besov spaces for which the integration in time is performed before the summation with respect to the frequency decomposition index:
\begin{defi} (\cite{Dbook} section 2.6.3)
 \sl{For $s,t\in \R$ and $a,b,c\in[1,\infty]$, we define the following norm
 $$
 \|u\|_{\tilde{L}_t^a \dot{B}_{b,c}^s}= \Big\| \left(2^{js}\|\ddj u\|_{L_t^a L^b}\right)_{j\in \Z}\Big\|_{l^c(\Z)}.
 $$
 The space $\tilde{L}_t^a \dot{B}_{b,c}^s$ is defined as the set of tempered distributions $u$ such that $\lim_{j \rightarrow -\infty} S_j u=0$ in $L^a([0,t],L^\infty(\R^d))$ and $\|u\|_{\tilde{L}_t^a \dot{B}_{b,c}^s} <\infty$.
 }
 \label{deftilde}
\end{defi}
Let us also recall the following proposition:
\begin{prop}
\sl{
For all $a,b,c\in [1,\infty]$ and $s\in \R$:
     $$
     \begin{cases}
    \mbox{if } a\leq c,& \forall u\in L_t^a \dot{B}_{b,c}^s, \quad \|u\|_{\tilde{L}_t^a \dot{B}_{b,c}^s} \leq \|u\|_{L_t^a \dot{B}_{b,c}^s}\\
    \mbox{if } a\geq c,& \forall u\in\tilde{L}_t^a \dot{B}_{b,c}^s, \quad \|u\|_{\tilde{L}_t^a \dot{B}_{b,c}^s} \geq \|u\|_{L_t^a \dot{B}_{b,c}^s}.
     \end{cases}
     $$
     \label{Propermut}
     }
\end{prop}
Let us continue with classical injections:
\begin{prop}
 \sl{(\cite{Dbook} Chapter 2) We have:
$$
 \begin{cases}
\mbox{For any } p\geq 1, & \dot{B}_{p,1}^0 \hookrightarrow L^p,\\
\mbox{For any } p\in[2,\infty[, & \dot{B}_{p,2}^0 \hookrightarrow L^p,\\
\mbox{For any } p\in[1,2], & \dot{B}_{p,p}^0 \hookrightarrow L^p.
\end{cases}
$$
}
 \label{injectionLr}
\end{prop}
\begin{lem} (\cite{Dbook} Section 2.11, see also Lemma $5$ from \cite{FCestimLp})
 \sl{For any $\aa, \beta>0$ there exists a constant $C_{\aa, \beta}>0$ such that for any $u\in \dot{H}^{s-\aa}(\R^d) \cap \dot{H}^{s+\beta}(\R^d)$, then $u\in\dot{B}_{2,1}^s(\R^d)$ and:
\begin{equation}
 \|u\|_{\dot{B}_{2,1}^s} \leq C_{\aa, \beta} \|u\|_{\dot{H}^{s-\aa}}^{\frac{\beta}{\aa + \beta}} \|u\|_{\dot{H}^{s+\beta}}^{\frac{\aa}{\aa + \beta}}.
\end{equation}
 }
\label{estimBsHs}
 \end{lem}
Let us end this section with a consequence of the previous result involving anisotropic Lebesgue spaces, that we define for $a,b\in [1,\infty]$ as follows:
\begin{equation}
  \|f\|_{L_{h,v}^{a,b}} \overset{def}{=} \big\|\|f(x_h,.)_{L^b(\R_v)}\|\big\|_{L^a(\R_h^2)}.
  \label{DefLaniso}
\end{equation}
\begin{lem}
 \sl{For any $\aa,\bb>0$, there exists a constant $C_{\aa,\bb}$ such that for any for any $c\in \dot{H}^{1-\aa}(\R^3) \cap \dot{H}^{1+\beta}(\R^3)$, then $u\in L_{h,v}^{\infty,2}(\R^3)$ and:
 $$
 \|c\|_{L_{h,v}^{\infty,2}(\R^3)} \leq C_{\aa,\bb} \|c\|_{\dot{H}^{1-\aa}(\R^3)}^\frac{\bb}{\aa+\bb} \|c\|_{\dot{H}^{1+\bb}(\R^3)}^\frac{\aa}{\aa+\bb}.
 $$
 }
 \label{estimBsHsAniso}
\end{lem}
\textbf{Proof:} for any function $c$
$$
 \|c\|_{L_{h,v}^{\infty,2}(\R^3)} \overset{def}{=} \underset{x_h\in\R^2}{\sup} \left(\int_\R |c(x_h,x_3|^2 dx_3\right)^\frac12 \leq \left(\int_\R \Big(\underset{x_h\in\R^2}{\sup}|c(x_h,x_3|\Big)^2 dx_3\right)^\frac12\\
$$
Thanks to the Bernstein lemma and to Lemma \ref{estimBsHs}, for any $\aa,\bb>0$ we have:
\begin{multline}
  \underset{x_h\in\R^2}{\sup}|c(x_h,x_3)|= \|c(.,x_3)\|_{L^\infty(\R^2)} \leq \|c(.,x_3)\|_{\dot{B}_{\infty,1}^0(\R^2)} \leq C \|c(.,x_3)\|_{\dot{B}_{2,1}^1(\R^2)}\\
 \leq C_{\aa,\bb} \|c(.,x_3)\|_{\dot{H}^{1-\aa}(\R^2)}^\frac{\bb}{\aa+\bb} \|c(.,x_3)\|_{\dot{H}^{1+\bb}(\R^2)}^\frac{\aa}{\aa+\bb}.
\end{multline}
This implies that
\begin{multline}
 \|c\|_{L_{h,v}^{\infty,2}(\R^3)} \leq C_{\aa,\bb} \left(\int_\R \||D_h|^{1-\aa} c(.,x_3)\|_{L^2(\R^2)}^\frac{2\bb}{\aa+\bb}\cdot \||D_h|^{1-\bb} c(.,x_3)\|_{L^2(\R^2)}^\frac{2\aa}{\aa+\bb} dx_3\right)^\frac12\\
 \leq C_{\aa,\bb} \|c\|_{\dot{H}^{1-\aa}(\R^3)}^\frac{\bb}{\aa+\bb} \|c\|_{\dot{H}^{1+\bb}(\R^3)}^\frac{\aa}{\aa+\bb}. \blacksquare
\end{multline}

\subsection{Product laws}

\subsubsection{General product laws}
We refer to \cite{Dbook} for the classical product laws and to \cite{CDGG, CDGG2, CDGGbook, FCRF} for the 2D-3D version (we also refer to \cite{FCStratif1, FCStratif2}) for a 1D-3D version.
\begin{prop}
 \sl{There exists a constant $C>0$ such that for any $s,t<\frac32$ with $s+t>0$ and any $u\in \dot{H}^s(\R^3)$, $v\in \dot{H}^t(\R^3)$, then $uv\in \dot{H}^{s+t-\frac32}(\R^3)$ and we have:
 $$
 \|uv\|_{\dot{H}^{s+t-\frac32}(\R^3)} \leq C \|u\|_{\dot{H}^s(\R^3)} \|v\|_{\dot{H}^t(\R^3)}.
 $$
 }
 \label{prod3D}
\end{prop}
\begin{prop}
 \sl{There exists a constant $C>0$ such that for any $s,t<1$ with $s+t>0$ and any $u\in \dot{H}^s(\R^2)$, $v\in \dot{H}^t(\R^3)$, then $uv\in \dot{H}^{s+t-1}(\R^3)$ and we have:
 $$
 \|uv\|_{\dot{H}^{s+t-1}(\R^3)} \leq C \|u\|_{\dot{H}^s(\R^2)} \|v\|_{\dot{H}^t(\R^3)}.
 $$
 }
 \label{prod2D3D}
\end{prop}
\begin{rem}
 \sl{When applying these propositions to estimate a product of the form $u\cdot \n v$ with $\div u=0$, the second condition on the regularity exponents can be relaxed into $s+t>-\frac32$ for the classical law and into $s+t>-1$ for the 2D-3D law.}
 \label{ProdRq}
\end{rem}

\subsubsection{Estimates for the external force terms}

In this section, we prove the following estimates for the external force terms in System \eqref{SystWe}, namely
$$
\bP (\tb\cdot \n \ce +\ce\cdot \n \tb +\ce \cdot \n \ce)
$$
\begin{prop}
 \sl{Assume $\delta\in]0,\frac14]$, $\tu_0,\tb_0\in H^\delta$ and $\coe\in H^{\frac12+\delta}$ (inhomogeneous spaces). Let $(\tu,\tb)$ and $\ce$ respectively denote the corresponding unique global strong solutions of Systems \eqref{MHD2Db} and \eqref{Mag2}. There exists a constant $C>0$ (depending on $\delta,\nu',\|\tu_0\|_{L^2}$) such that for any $(\sigma,\sigma')\in[0,\frac12+\delta]\times [0,\delta]$:
 \begin{equation}
  \begin{cases}
  \vspace{0.1cm}
   \|\ce\cdot \n \ce\|_{L^1 \dot{H}^{\sigma}} \leq C \|\coe\|_{H^{\frac12+\delta}}^2,\\
   \|\tb\cdot \n \ce\|_{L^\frac43 \dot{H}^{\sigma'}} +\|\ce\cdot \n \tb\|_{L^\frac43 \dot{H}^{\sigma'}} \leq C \|\tb_0\|_{H^\delta} \|\coe\|_{H^{\frac12+\delta}}.
  \end{cases}
\label{EstimFext}
 \end{equation}
 }
\label{PropFext}
\end{prop}
\textbf{Proof:} the first term in $\eqref{EstimFext}_2$ is easily bounded thanks to Proposition \ref{prod2D3D} with $(s,t)=(\frac12,\frac12+\sigma')$ (provided that $\delta<\frac12$) and interpolation. For any $t\geq 0$, we have:
$$
\|\tb\cdot \n \ce\|_{L_t^\frac43 \dot{H}^{\sigma'} (\R^3)} \leq C\|\tb\|_{L_t^4 \dot{H}^\frac12 (\R^2)} \|\n \ce\|_{L_t^2 \dot{H}^{\frac12+\sigma'} (\R^3)} \leq \|\tb\|_{L_t^\infty L^2}^\frac12 \|\tb\|_{L_t^2 \dot{H}^1}^\frac12 \|\n \ce\|_{L_t^2 \dot{H}^{\frac12+\sigma'}}.
$$
Thanks to Theorems \ref{ThExistlim} and \ref{Thc}, we end-up with
$$
\|\tb\cdot \n \ce\|_{L_t^\frac43 \dot{H}^{\sigma'}} \leq C \|\tb_0\|_{L^2} \|\coe\|_{\dot{H}^{\frac12+\sigma'}} \leq C \|\tb_0\|_{H^\delta} \|\coe\|_{H^{\frac12+\delta}}.
$$
We could deal with $\eqref{EstimFext}_1$ using classical product laws from Proposition \ref{prod3D} but due to the condition $s,t<\frac32$ we would end-up at best with (for some $m>0$)
$$
\|\ce\cdot \n \ce\|_{\dot{H}^{\frac12+\delta}} \leq \|\ce\|_{\dot{H}^{\frac32-m\delta}} \|\n \ce\|_{\dot{H}^{\frac12+(1+m)\delta}},
$$
which requires extra-regularity from $\coe$. To avoid this problem we simply do as in \cite{FCPAA}: for any $\sigma\in[0,\frac12+\delta]$
\begin{multline}
 \|\ce\cdot \n \ce\|_{\dot{H}^{\sigma}} =\|\div (\ce\otimes \ce)\|_{\dot{H}^{\sigma}} \leq 2\|T_{\ce} \ce\|_{\dot{H}^{1+\sigma}} +\|R(\ce,\ce)\|_{\dot{H}^{1+\sigma}}\\
 \leq C\big(2\|\ce\|_{L^\infty} +\|\ce\|_{\dot{B}_{\infty,\infty}^0}\big) \|\ce\|_{\dot{H}^{1+\sigma}}.
\end{multline}
Thanks to the classical Bernstein lemma and to Lemma \ref{estimBsHs}, we have
$$
 2\|\ce\|_{L^\infty} +\|\ce\|_{\dot{B}_{\infty,\infty}^0} \leq 3 \|\ce\|_{\dot{B}_{\infty,1}^0} \leq C \|\ce\|_{\dot{B}_{2,1}^\frac32} \leq C \|\ce\|_{\dot{H}^{\frac32-\delta}}^\frac12 \|\ce\|_{\dot{H}^{\frac32+\delta}}^\frac12,
$$
so that we obtain for any $t\geq 0$:
\begin{multline}
  \|\ce\cdot \n \ce\|_{L_t^1 \dot{H}^{\sigma}} \leq C \|\ce\|_{L_t^2 \dot{H}^{\frac32-\delta}}^\frac12 \|\ce\|_{L_t^2 \dot{H}^{\frac32+\delta}}^\frac12 \|\ce\|_{L_t^2 \dot{H}^{1+\sigma}}\\
  \leq C \|\coe\|_{\dot{H}^{\frac12-\delta}}^\frac12 \|\coe\|_{\dot{H}^{\frac12+\delta}}^\frac12 \|\coe\|_{\dot{H}^{\sigma}} \leq C \|\coe\|_{H^{\frac12+\delta}}^2,
\end{multline}
which gives $\eqref{EstimFext}_1$. The last term will require most of the work. The following decomposition holds:
$$
\ce\cdot \n \tb =T_{\ce} \n \tb +T_{\n \tb} \ce +R(\ce, \n\tb),
$$
where we have defined
\begin{equation}
 \begin{cases}
 T_{\ce} \n \tb \overset{def}{=} \Sum_{q\in \Z} \dot{S}_{q-j_0} \ce \cdot \ddq^h \n \tb,\\
 T_{\n \tb} \ce \overset{def}{=} \Sum_{q\in \Z} \dot{S}_{q-j_0}^h \n \tb \cdot \ddq \ce,\\
 R(\ce, \n\tb) \overset{def}{=} \Sum_{q\in \Z} \Sum_{\alpha=-j_0}^{j_0} \dot{\D}_q \ce \cdot \n \ddqa^h \tb,
\end{cases}
\label{paraprodecp}
\end{equation}
and where $\dot{\D}_q,\dot{S}_q$ refer to the classical dyadic truncation operators in $\R^3$ (introduced in the previous subsection) and $\dot{\D}_q^h,\dot{S}_q^h$ refer to their counterpart in $\R^2$ acting on the horizontal variables $x_h=(x_1,x_2)$.\\
\begin{lem}
 \sl{If $j_0\geq 4$, for any $q\in \Z$, we have
 \begin{equation}
 \begin{cases}
  \mbox{supp}\; \left(\cF(\dot{S}_{q-j_0} \ce \cdot \ddq^h \n \tb)\right) \bigcup \mbox{supp}\; \left(\cF(\dot{S}_{q-j_0}^h \n \tb \cdot \ddq \ce)\right) \subset 2^q \cC(0,\frac1{12},\frac{10}3),\\
  \mbox{supp}\; \left(\cF(\Sum_{\alpha=-j_0}^{j_0} \dot{\D}_q \ce \cdot \n \ddqa^h \tb)\right) \subset 2^q B(0,46).
 \end{cases}
\label{LemmeSupports}
 \end{equation}
 }
\end{lem}
\textbf{Proof of the lemma :} we emphasize that contrary to the classical decomposition, there is now a slight asymmetry between the previous paraproduct terms due the 2D-3D product which will require us to fix $j_0$ large enough in order to use the classical methods. For $f=f(x)$ and $\tg=\tg(x_h)$, we formally have (denoting the same way the Fourier transform in $\R^3$ and $\R_h^2$):
\begin{equation}
 \cF(f\tg)(\xi)=\int_{\R^2} \hat{f}(\xi_h-\eta_h,\xi_3)\hat{\tg}(\eta_h) d\eta_h.
 \label{Convol2D}
\end{equation}
Thanks to this, assuming that $\mbox{supp}\; \hat{\tg}\subset \R^2$ and $\mbox{supp}\; \hat{f}\subset \R^3$ are compact sets, and if we define compact sets $K_h\subset \R^2$ and $K_v\subset \R$ such that $\mbox{supp}\; \hat{f}\subset K_h\times K_v$, it easy to see that
\begin{equation}
\mbox{supp}\; \cF(f\tg)\subset \left(K_h+\mbox{supp}\; \hat{\tg}\right)\times K_v.
\label{Supports}
\end{equation}
Thanks to this, we prove $\eqref{LemmeSupports}_3$ easily. In the first paraproduct from \eqref{paraprodecp}, as for any $q\in\Z$, we have
$$
\begin{cases}
\vspace{1mm}
 \mbox{supp}\; \cF(\dot{S}_{q-j_0} \ce) \subset 2^{q-j_0}B(0,\frac43) \subset 2^q B_h(0,\frac43 2^{-j_0}) \times 2^q [-\frac43 2^{-j_0},\frac43 2^{-j_0}],\\
  \mbox{supp}\; \cF(\ddq^h \n \tb) \subset 2^q\cC_h(0,\frac34, \frac83),
\end{cases}
$$
so that, thanks to \eqref{Supports}, if $j_0\geq 1$, we get $\mbox{supp}\; \cF(\dot{S}_{q-j_0} \ce \cdot \ddq^h \n \tb) \subset 2^q \cC(0,\frac1{12}, \frac{10}3)$.\\
Let us now turn to the last term: writing \eqref{Convol2D} for $\cF(\dot{S}_{q-j_0}^h \n \tb \cdot \ddq \ce)$, for the integrated product to be non zero we need that:
$$
\begin{cases}
\vspace{1mm}
 |\eta_h| \leq 2^{q-j_0} \frac43,\\
 2^q \frac34 \leq \sqrt{|\xi_h-\eta_h|^2+\xi_3^2} \leq 2^q \frac83.
\end{cases}
$$
This implies that if $j_0\geq 4$:
\begin{multline}
 |\xi|^2 \geq |\xi_h|^2+\xi_3^2 \geq |\xi_h-\eta_h|^2-2|\xi_h-\eta_h|\cdot |\eta_h|+|\eta_h^2|+\xi_3^2\\
 \geq 2^{2q} \left(\frac9{16}-\frac{2^{6-j_0}}9 \right) \geq 2^{2q}\frac{17}{144}>2^{2q}\frac19.
\end{multline}
So, when $j_0\geq 4$ we finally obtain that $\mbox{supp}\; \cF(\dot{S}_{q-j_0}^h \n \tb \cdot \ddq \ce) \subset 2^q \cC(0,\frac13,\frac{11}4)$, which ends the proof of the lemma $\blacksquare$.\\

We can now focus on estimating the terms from \eqref{paraprodecp}. Fixing $j_0=4$ and thanks to the previous lemma we can use Lemma 2.47 from \cite{Dbook} and write that for any $\sigma'\in[0,\delta]$:
\begin{equation}
 \|T_{\n \tb} \ce\|_{\dot{H}^{\sigma'}} \leq \left\|(2^{q\sigma'}\|\dot{S}_{q-4}^h  \n\tb \cdot \ddq \ce\|_{L^2})_{q\in\Z}\right\|_{l^2(\Z)}.
 \label{Paraprod-reduction}
\end{equation}
For any $q\in\Z$, thanks to the Bernstein lemma and convolutions estimates (see \eqref{DefBesov}), we have
$$
\|\dot{S}_{q-4}^h \n\tb \cdot \ddq \ce\|_{L^2(\R^3)} \leq \|\dot{S}_{q-4}^h \n\tb\|_{L^\infty(\R^2)} \|\ddq \ce\|_{L^2(\R^3)} \leq C 2^q\|\tb\|_{L^\infty(\R^2)} \|\ddq \ce\|_{L^2(\R^3)}.
$$
Multiplying by $2^{q\sigma'}$ and summing in $l^2(\Z)$, we obtain that thanks to \eqref{Paraprod-reduction} and Lemma \ref{estimBsHs},
\begin{multline}
  \|T_{\n \tb} \ce\|_{\dot{H}^{\sigma'}} \leq C \|\tb\|_{L^\infty(\R^2)} \|\ce\|_{\dot{H}^{1+\sigma'}} \leq C \|\tb\|_{\dot{B}_{2,1}^1(\R^2)} \|\ce\|_{\dot{H}^{1+\sigma'}}\\
  \leq C \|\tb\|_{\dot{H}^{1-\sigma'}(\R^2)}^\frac12 \|\tb\|_{\dot{H}^{1+\sigma'}(\R^2)}^\frac12 \|\ce\|_{\dot{H}^{1+\sigma'}},
\end{multline}
As $\tb_0\in H^{\delta}(\R^2)$, thanks to \eqref{estimtutbHs} we know that $\|\tb\|_{\dot{H}^{1+\sigma'}}\in L^2$. Thanks to \eqref{estimtutb} and by interpolation we have $\|\tb\|_{\dot{H}^{1-\sigma'}} \in L^\frac2{1-\sigma'}$. Thanks to Proposition \ref{Propestimce}, $\|\ce\|_{\dot{H}^{1+\sigma'}}\in L^p$ for any $p\in[2,\frac4{1-2(\delta-\sigma')}]$, so taking in particular $p=\frac4{1+\sigma'}$, we get:
\begin{equation}
 \|T_{\n \tb} \ce\|_{L^\frac43 \dot{H}^{\sigma'}} \leq C \|\tb\|_{L^\frac2{1-\sigma'}\dot{H}^{1-\sigma'}(\R^2)}^\frac12 \|\tb\|_{L^2 \dot{H}^{1+\sigma'}(\R^2)}^\frac12 \|\ce\|_{L^\frac4{1+\sigma'} \dot{H}^{1+\sigma'}}.
 \end{equation}
Thanks to Theorems \ref{ThExistlim}, \ref{Thc}, and Proposition \ref{Propestimce}, we finally obtain:
\begin{equation}
 \|T_{\n \tb} \ce\|_{L_t^\frac43 \dot{H}^{\sigma'}} \leq C \|\tb_0\|_{H^\delta} \|\coe\|_{H^{\frac12+\delta}}.
\label{EstimT1}
 \end{equation}
As $\div \ce=0$, the remainder term is classically rewritten as follows:
$$
R(\ce, \n\tb) \overset{def}{=} \Sum_{i=1}^3 \d_i\left(\Sum_{q\in \Z} \Sum_{\alpha=-4}^{4} \dot{\D}_q \ce^i \cdot \ddqa^h \tb\right),
$$
so that
$$
\|R(\ce, \n\tb)\|_{\dot{H}^{\sigma'}} \leq C \left\|\Sum_{i=1}^3 \Sum_{q\in \Z} \Sum_{\alpha=-4}^{4} \dot{\D}_q \ce^i \cdot \ddqa^h \tb\right\|_{\dot{H}^{1+\sigma'}}.
$$
Denoting as $\big(c_q(\ce)\big)_{q\in\Z}$ a positive sequence of $l^2$-norm equal to $1$, and for some fixed $N_0$ only depending on the radius from $\eqref{LemmeSupports}_3$, we have for any $q'\in \Z$,
\begin{multline}
 \left\|\dot{\D}_{q'}\Big(\Sum_{i=1}^3 \Sum_{q\in \Z} \Sum_{\alpha=-4}^{4} \dot{\D}_q \ce^i \cdot \ddqa^h \tb\Big) \right\|_{L^2} \leq C \Sum_{q\geq q'-N_0} \Sum_{\alpha=-4}^{4} \|\ddq \ce\|_{L^2(\R^3)} \|\ddqa^h \tb\|_{L^\infty(\R^2)}\\
 \leq C \Sum_{q\geq q'-N_0} 2^{-q(1+\sigma')}c_q(\ce)\|\ce\|_{\dot{H}^{1+\sigma'}} \|\tb\|_{L^\infty(\R^2)}.
\end{multline}
This entails
$$
2^{q'(1+\sigma')} \left\|\dot{\D}_{q'}\Big(\Sum_{i=1}^3 \Sum_{q\in \Z} \Sum_{\alpha=-4}^{4} \dot{\D}_q \ce^i \cdot \ddqa^h \tb\Big) \right\|_{L^2} \leq C \big(a*c(\ce)\big)_{q'} \|\ce\|_{\dot{H}^{1+\sigma'}} \|\tb\|_{L^\infty},
$$
where we define the $l^1$ sequence $(a_q)_{q\in\Z}$ as follows: for any $q\in\Z$:
$$
a_q=\begin{cases}
     2^{q(1+\sigma')} & \mbox{if }q\leq N_0,\\
     0 & \mbox{else}.
    \end{cases}
$$
Taking the $l^2$-norm, we obtain as previously that:
\begin{multline}
 \|R(\ce, \n\tb)\|_{L_t^\frac43 \dot{H}^{\sigma'}} \leq C \|\tb\|_{L^\frac2{1-\sigma'}\dot{H}^{1-\sigma'}(\R^2)}^\frac12 \|\tb\|_{L^2 \dot{H}^{1+\sigma'}(\R^2)}^\frac12 \|\ce\|_{L^\frac4{1+\sigma'} \dot{H}^{1+\sigma'}} \leq C \|\tb_0\|_{H^\delta} \|\coe\|_{H^{\frac12+\delta}}.
 \label{EstimR}
\end{multline}
Turning to the last term, using first Lemma \ref{estimBsHsAniso} and then convolutions from \eqref{DefBesov},  we have that for any $q\in \Z$,
\begin{multline}
  \|\dot{S}_{q-4} \ce \cdot \ddq^h \n \tb\|_{L^2(\R^3)} \leq \|\dot{S}_{q-4} \ce\|_{L_{h,v}^{\infty,2}(\R^3)} \|\ddq^h \n \tb\|_{L^2(\R^2)}\\
  \leq C_{\sigma'} \|\dot{S}_{q-4} \ce\|_{\dot{H}^{1-\sigma'}(\R^3)}^\frac12 \|\dot{S}_{q-4} \ce\|_{\dot{H}^{1+\sigma'}(\R^3)}^\frac12 2^q \|\ddq^h \tb\|_{L^2(\R^2)}\\
  \leq C_{\sigma'} \|\ce\|_{\dot{H}^{1-\sigma'}(\R^3)}^\frac12 \|\ce\|_{\dot{H}^{1+\sigma'}(\R^3)}^\frac12 2^q \|\ddq^h \tb\|_{L^2(\R^2)}.
\end{multline}
Thanks to $\eqref{LemmeSupports}_1$ we conclude that:
$$
\|T_{\ce} \n \tb\|_{\dot{H}^{\sigma'}} \leq C_{\sigma'} \|\ce\|_{\dot{H}^{1-\sigma'}(\R^3)}^\frac12 \|\ce\|_{\dot{H}^{1+\sigma'}(\R^3)}^\frac12 \|\tb\|_{\dot{H}^{1+\sigma'}(\R^2)}.
$$
In addition to the fact that $\|\ce\|_{\dot{H}^{1+\sigma'}}\in L_t^p$ for any $p\in[2,\frac4{1-2(\delta-\sigma')}]$, thanks once more to Proposition \ref{Propestimce} we know that when $\delta\leq \frac16$, $\|\ce\|_{\dot{H}^{1-\sigma'}}\in L_t^p$ for any $p\in[\frac2{1-\sigma'},\frac4{1-2(\delta+\sigma')}]$ which allows us to obtain that for some $C=C_{\nu',\delta,\|\tu_0\|_{L^2}}>0$:
$$
 \|T_{\ce} \n \tb\|_{L_t^\frac43\dot{H}^{\sigma'}} \leq C_{\sigma'} \|\ce\|_{L_t^4 \dot{H}^{1-\sigma'}(\R^3)}^\frac12 \|\ce\|_{L_t^4 \dot{H}^{1+\sigma'}(\R^3)}^\frac12 \|\tb\|_{L_t^2 \dot{H}^{1+\sigma'}(\R^2)} \leq C\|\coe\|_{H^{\frac12+\delta}} \|\tb_0\|_{H^\delta}.
$$
Gathering the previous estimates with \eqref{EstimT1} and \eqref{EstimR} ends the proof of $\eqref{EstimFext}_2$. $\blacksquare$

\subsection{Strichartz estimates}

\subsubsection{Isotropic estimates}

Consider the following system (external force terms $\Fe,\Ge$ depend on $(t,x)$)
\begin{equation}
 \begin{cases}
  \d_t f -\nu \D f +\frac{1}{\ee}\mathbb{P} (f\wedge e_3)= \Fe +\Ge,\\
  {f}_{|t=0}=\foe.
 \end{cases}
\label{LRF}
\end{equation}
The Strichartz estimates we use in this article are slight extensions of the ones from \cite{FCRF}:
\begin{prop}
 \sl{\begin{enumerate}
      \item (Case $\Ge=0$) For any $d\in \R$, $r\geq2$, $q\geq 1$, $\theta\in[0,1]$, and $p\in[1, \frac2{\theta(1-\frac2{r})}]$, there exists a constant $C=C_{p,\theta,r}>0$ such that for any divergence-free vectorfield $\foe$ and $\Fe$, the solution $f$ of \eqref{LRF} satisfies:
    \begin{equation}
  \||D|^d f\|_{\tilde{L}^p\dot{B}_{r, q}^0} \leq \frac{C_{p,\theta,r}}{\nu^{\frac{1}{p}-\frac{\theta}{2}(1-\frac2{r})}} \ee^{\frac{\theta}{2}(1-\frac2{r})} \left(\|\foe\|_{\dot{B}_{2, q}^{\sigma_1}} +\|\Fe\|_{L^1 \dot{B}_{2, q}^{\sigma_1}}\right),
\end{equation}
with $\sigma_1= d+\frac32-\frac3{r}-\frac2{p}+\theta(1-\frac2{r})$.
      \item (Case $(\foe, \Fe)=(0,0)$) For any $d\in \R$, $k\in]1,2]$, $r\geq2$, $q\geq 1$, $\theta\in[0,1]$, and $p\in[2, \frac2{\theta(1-\frac2{r})}[$, there exists a constant $C=C_{p,\theta,r,k}>0$ such that for any divergence-free vectorfield $\Ge$, we have:
\begin{equation}
  \||D|^d f\|_{\tilde{L}^p\dot{B}_{r, q}^0} \leq \frac{C_{p,\theta,r,k}}{\nu^{1-\frac1{k}+\frac{1}{p}-\frac{\theta}{2}(1-\frac{2}{r})}} \ee^{\frac{\theta}{2}(1-\frac{2}{r})} \|\Ge\|_{\tilde{L}^k \dot{B}_{2, q}^{\sigma_2}},
 \end{equation}
 where $\sigma_2= d-\frac12+\frac2{k}-\frac{3}{r}-\frac{2}{p}+\theta (1-\frac{2}{r})$.
\end{enumerate}
\label{EstimStri}
}
\end{prop}
The first point is easily proved from Proposition 10 from \cite{FCRF} (where we only needed to consider the homogeneous case) and the Duhamel formula. The second point is an adaptation (in the case of the rotating fluids) of Proposition 6 from \cite{FCRF}. In particular, as the dispersion estimates are a little better for the rotating fluids compared to the primitive system (see Section 5.2 from \cite{FCcompl} and Section 4.3.2 from \cite{FCRF}) this explains why, compared to Proposition 6 from \cite{FCRF}, $\theta$ is upgraded into $2\theta$.
\begin{rem}
 \sl{\begin{enumerate}
      \item Let us emphasize that as in the cases $\nu=\nu'$ of Primitive or Strongly stratified Boussinesq systems, these estimates hold for any $\ee>0$ (the threshold in $\ee$ will be set in the bootstrap argument and in Proposition \ref{PropestimStriStrong}).
      \item Of course, when the coefficients satisfy the (more restrictive) conditions from the second point, we can write that (using Proposition \ref{Propermut} and the fact that $k\leq 2$ in the last term):
 \begin{equation}
  \||D|^d f\|_{\tilde{L}^p\dot{B}_{r, q}^0} \leq \frac{\ee^{\frac{\theta}{2}(1-\frac2{r})}}{\nu^{\frac{1}{p}-\frac{\theta}{2}(1-\frac2{r})}}  \left(C_{p,\theta,r}\Big(\|\foe\|_{\dot{B}_{2, q}^{\sigma_1}} +\|\Fe\|_{L^1 \dot{B}_{2, q}^{\sigma_1}}\Big) +\frac{C_{p,\theta,r,k}}{\nu^{1-\frac1{k}}} \|\Ge\|_{L^k \dot{B}_{2, q}^{\sigma_1+\frac2{k}-2}} \right).
\label{EstimStriIsoUnif}
\end{equation}
For the sake of simplicity, we will sometimes use the previous proposition on the velocity $\ve$ splitting it into several parts, thanks to the superposition principle.
\end{enumerate}
}
\end{rem}

\subsubsection{Anisotropic estimates}
The anisotropic Lebesgue spaces recalled previously appear in the following anisotropic Strichartz estimates that we present in two versions (depending on the regularity of the external force) as in the previous section:

\begin{prop}
 \sl{\begin{enumerate}
      \item (Case $\Ge=0$) For any $d\in \R$, $m>2$, $\theta\in]0,1]$, $p\in[1, \frac4{\theta(1-\frac2{m})}]$ there exists a constant $C=C_{p,\theta,m}$ such that for any divergence-free vectorfields $\foe$ and $\Fe$, the solution $f$ of \eqref{LRF} satisfies:
    \begin{equation}
  \||D|^d f\|_{L^p L_{h,v}^{m,2}} \leq \frac{C_{p,\theta,m}}{\nu^{\frac{1}{p}-\frac{\theta}{4}(1-\frac2{m})}} \ee^{\frac{\theta}{4}(1-\frac2{m})} \left(\|\foe\|_{\dot{B}_{2, 1}^{\sigma_3}} +\|\Fe\|_{L^1 \dot{B}_{2, 1}^{\sigma_3}}\right),
\end{equation}
 with $\sigma_3= d+1-\frac2{m}-\frac2{p}+\frac{\theta}2 (1-\frac2{m})$.
 \item (Case $(\foe, \Fe)=(0,0)$) For any $d\in \R$, $k\in]1,2]$, $m>2$, $\theta\in]0,1]$, $p\in[2, \frac4{\theta(1-\frac2{m})}[$, there exists a constant $C=C_{p,\theta,m,k}>0$ such that for any divergence-free vectorfield $\Ge$, we have:
\begin{equation}
  \||D|^d f\|_{L^p L_{h,v}^{m,2}} \leq \frac{C_{p,\theta,m,k}}{\nu^{1-\frac1{k}+\frac{1}{p}-\frac{\theta}4(1-\frac2{m})}} \ee^{\frac{\theta}4(1-\frac2{m})} \|\Ge\|_{\tilde{L}^k \dot{B}_{2, 1}^{\sigma_4}},
\end{equation}
 with $\sigma_4= d+\frac2{k}-1-\frac2{m}-\frac2{p}+\frac{\theta}2 (1-\frac2{m})$.
     \end{enumerate}
}
 \label{EstimStrianiso}
\end{prop}
\textbf{Proof:} In \cite{FCRF} we only considered the homogeneous case. The cases for $\Fe$ ($k=1$) and $\Ge$ ($k\in]1,2]$) are treated adapting our arguments from \cite{FCRF} (see the proof of Proposition 17 therein).
\begin{rem}
 \sl{As previously, when the coefficients satisfy the assumptions from the second point, we can write that:
 \begin{equation}
  \||D|^d f\|_{L^p L_{h,v}^{m,2}} \leq \frac{\ee^{\frac{\theta}4(1-\frac2{m})}}{\nu^{\frac{1}{p}-\frac{\theta}4(1-\frac2{m})}} \left(C_{p,\theta,m}\Big(\|\foe\|_{\dot{B}_{2, 1}^{\sigma_3}} +\|\Fe\|_{L^1 \dot{B}_{2, 1}^{\sigma_3}}\Big) +\frac{C_{p,\theta,m,k}}{\nu^{1-\frac1{k}}} \|\Ge\|_{L_t^k \dot{B}_{2, 1}^{\sigma_3+\frac2{k}-2}} \right).
  \label{EstimStriAnisoUnif}
\end{equation}
}
\end{rem}

\subsubsection{Proof of Proposition \ref{PropestimStriStrong}}
\label{Sect:IsoStri}
We prove in this section the adapted isotropic Strichartz estimates  for System \eqref{SystWe}.

Taking advantage of Propositions \ref{Propermut}, \ref{injectionLr} and \ref{PropFext}, we can rewrite Proposition \ref{EstimStri} for System \eqref{SystWe} with $(q,k)=(2,\frac43)$ and for fixed regularity: for any $d\in \R$, $r\geq2$, $\theta\in[0,1]$, and $p\in[2, \frac2{\theta(1-\frac2{r})}[$ such that
\begin{equation}
 d+\frac32-\frac3{r}-\frac2{p}+\theta(1-\frac2{r})=\frac12+\delta,
\label{Cond:paramiso}
\end{equation}
there exists a constant $C=C_{p,\theta,r}>0$ such that the solution of \eqref{SystWe} satisfies:
\begin{multline}
  \||D|^d \We\|_{L^p L^r} \leq \||D|^d \We\|_{L^p\dot{B}_{r, 2}^0} \leq \||D|^d \We\|_{\tilde{L}^p\dot{B}_{r, 2}^0}\\
  \leq \frac{C_{p,\theta,r}}{\nu^{\frac{1}{p}-\frac{\theta}{2}(1-\frac2{r})}}  \ee^{\frac{\theta}{2}(1-\frac2{r})} \left(\|\voe\|_{\dot{H}^{\frac12+\delta}} +\|\ce\cdot \n \ce\|_{L_t^1 \dot{H}^{\frac12+\delta}} +\frac1{\nu^\frac14} \|\tb\cdot \n \ce +\ce \cdot \n \tb\|_{L^\frac43 \dot{H}^\delta} \right)\\
  \leq C_{\nu, p,\theta,r} \ee^{\frac{\theta}{2}(1-\frac2{r})} \left(\|\voe\|_{\dot{H}^{\frac12+\delta}} +\|\coe\|_{H^{\frac12+\delta}}^2 + \|\tbo\|_{H^\delta}^2 \right).
  \label{EstimStriIsoW}
\end{multline}
\begin{itemize}
 \item Taking $(d,p,r)=(1,2,3)$, we obtain that \eqref{Cond:paramiso} is satisfied if and only if $\theta=3\delta$, \eqref{EstimStriIsoW} then turns into
$$
\|\n \We\|_{L^2 L^3} \leq C_{\nu, \nu', \Co, \delta} \ee^{\frac{\delta}2} \left(\|\voe\|_{\dot{H}^{\frac12+\delta}} +\|\coe\|_{H^{\frac12+\delta}}^2 + \|\tbo\|_{H^\delta}^2 \right),
$$
together with the following conditions
\begin{equation}
 \begin{cases}
  \theta\in]0,1]& \Longleftrightarrow\quad \delta\leq \frac13,\\
  p\in[2, \frac2{\theta(1-\frac2{r})}[& \Longleftrightarrow\quad \delta<1.
 \end{cases}
\label{Condelta1}
\end{equation}
Thanks to the definitions of $\gamma,\Ko$ (see \eqref{ChoixKo}) and the assumptions from Theorem \ref{ThCVStrong},
$$
\|\n \We\|_{L^2 L^3} \leq C_{\nu, \nu', \Co, \delta} \ee^{\frac{\delta}2} \left(\Co \ee^{-\gamma} +\Big(\frac{\eta_0 \delta}{144\Do}|\ln \ee|\Big)^{\frac12} + \|\tbo\|_{H^\delta}^2 \right),
$$
and if $\ee\in ]0,\ee_0]$ with $\ee_0$ depending on $\Co,\Do,\delta,\eta_0$ so small that (we recall that $\Do$ is a constant featured in \eqref{Estimapriori3}) for any $\ee\in]0,\ee_0]$,
\begin{equation}
\Big(\frac{\eta_0 \delta}{144\Do}|\ln \ee|\Big)^{\frac12} \leq \min \left(\ee^{-\frac{\eta_0 \delta}{96}}, \ee^{-\gamma}\right),
\label{Condeps2}
\end{equation}
we finally get that:
\begin{equation}
\|\n \We\|_{L^2 L^3} \leq C_{\nu, \nu',\Co, \delta} \ee^{\frac{\delta}2-\gamma} =C_{\nu, \nu',\Co, \delta} \ee^{\eta_0 \delta},
\end{equation}
\item Choosing $(d,p,r)=(0,\frac2{1-s},6)$, \eqref{Cond:paramiso} is satisfied if and only if $\theta=\frac32(\frac12+\delta-s)$, and \eqref{EstimStriIsoW} becomes:
\begin{multline}
\|\We\|_{L^\frac2{1-s} L^6} \leq C_{\nu, \nu', \delta,s} \ee^{\frac12(\frac12+\delta-s)} \left(\|\voe\|_{\dot{H}^{\frac12+\delta}} +\|\coe\|_{H^{\frac12+\delta}}^2 + \|\tbo\|_{H^\delta}^2 \right)\\
\leq C_{\nu, \nu', \delta,s} \ee^{(1-\frac{\eta_0}{10})\frac{\delta}2} \left(\|\voe\|_{\dot{H}^{\frac12+\delta}} +\|\coe\|_{H^{\frac12+\delta}}^2 + \|\tbo\|_{H^\delta}^2 \right).
\end{multline}
with the corresponding conditions (as we have $\eta_0\leq \frac12$)
\begin{equation}
 \begin{cases}
  \theta\in]0,1]\mbox{ (for any }s\in[\frac12-\frac{\eta_0}{10} \delta, \frac12+\frac{\eta_0}{10} \delta])& \Longleftrightarrow\quad (1+\frac{\eta_0}{10})\delta\leq \frac23,\\
  p\in[2, \frac2{\theta(1-\frac2{r})}[& \Longleftrightarrow\quad \delta< \frac12.
 \end{cases}
\label{Condelta2}
\end{equation}
If $\ee\leq \ee_0$, we obtain for any $s\in[\frac12-\frac{\eta_0}{10}\delta, \frac12+\frac{\eta_0}{10}\delta]$,
$$
\|\We\|_{L^\frac2{1-s} L^6} \leq C_{\nu, \nu',\Co, \delta} \ee^{(1-\frac{\eta_0}{10})\frac{\delta}2-\gamma} =C_{\nu, \nu', \delta,\Co} \ee^{\frac{19}{20} \eta_0\delta} \leq C_{\nu, \nu', \delta,\Co} \ee^{\frac{9 \eta_0}{10}\delta}.
$$
\item The previous case when replacing $s$ by $\frac12$ gives (when $\delta<\frac12$):
\begin{equation}
 \|\We\|_{L^4 L^6} \leq C_{\nu, \nu',\Co, \delta} \ee^{\frac{\delta}2} \left(\|\voe\|_{\dot{H}^{\frac12+\delta}} +\|\coe\|_{H^{\frac12+\delta}}^2 + \|\tbo\|_{H^\delta}^2 \right) \leq C_{\nu, \nu', \delta,\Co} \ee^{\eta_0 \delta}
\end{equation}
\item Choosing $(d,p,r)=(1,\frac2{1-s(1-(\aa_0+1)\delta)}, \frac3{\frac32-\aa_0\delta})$, \eqref{Cond:paramiso} is satisfied if and only if
$$
\frac23 \aa_0\delta \theta=\frac12+(1-\aa_0)\delta -s(1-(1+\aa_0)\delta).
$$
We recall that $\aa_0$ was chosen in \eqref{Choixalpha} to ensure that $\theta \in[0,1]$. And the condition on $p$ is satisfied if, and only if $(1-\aa_0)\delta<\frac12$, which is true as $\aa_0$ is larger than $1$ and $\delta \leq \frac16$. This allows us to write that (for some constant $C_{\nu, \nu',\delta, \eta_0, s}$ bounded when $s\in[\frac12-\frac{\eta_0}{10}\delta, \frac12+\frac{\eta_0}{10}\delta]$ by some $C_{\nu, \nu',\delta, \eta_0}$):
\begin{multline}
 \|\n \We\|_{L_t ^\frac2{1-s(1-(\aa_0+1)\delta)}L^\frac3{\frac32-\aa_0\delta}} \leq C_{\nu, \nu',\delta, s, \eta_0} \ee^{\frac{\aa_0\delta \theta}3} \left(\|\voe\|_{\dot{H}^{\frac12+\delta}} +\|\coe\|_{H^{\frac12+\delta}}^2 + \|\tbo\|_{H^\delta}^2 \right)\\
 \leq C_{\nu, \nu',\delta, \eta_0,\Co} \ee^{\frac{\aa_0\delta \theta}3-\gamma}.
\end{multline}
Thanks to the notations in Theorem \ref{ThCVStrong}, we can bound from below the previous exponent as follows (as we wish a bound valid for any $s\in[\frac12-\frac{\eta_0}{10}\delta, \frac12+\frac{\eta_0}{10}\delta]$):
\begin{multline}
 \frac{\aa_0\delta \theta}3-\gamma =\frac14+(1-\aa_0)\frac{\delta}2 -s(\frac12-(1+\aa_0)\frac{\delta}2)-\frac{\delta}2(1-2\eta_0)\\
 \geq \frac14+(1-\aa_0)\frac{\delta}2 -\frac14(1+\frac{\eta_0}{5}\delta)\Big(1-(1+\aa_0)\delta\Big)-\frac{\delta}2(1-2\eta_0)\\
 \geq \frac{\delta}2 \left(\frac12+\frac{\eta_0}{10}(19+\delta)-\frac12(1-\frac{\eta_0}5 \delta) \aa_0\right) \geq \frac{\delta}{2(7+\frac35 \eta_0 \delta)}(-1+13\eta_0)\\
 \geq \frac{\delta}{16}(-1+13\eta_0) \geq \frac{\eta_0}{16} \delta,
\end{multline}
where the last estimate is true as $\eta_0\geq \frac1{12}$. This finally gives that for any $s\in[\frac12-\frac{\eta_0}{10}\delta, \frac12+\frac{\eta_0}{10}\delta]$:
$$
\|\n \We\|_{L_t ^\frac2{1-s(1-(\aa_0+1)\delta)}L^\frac3{\frac32-\aa_0\delta}} \leq C_{\nu, \nu',\delta, \eta_0} \ee^{\frac{\eta_0}{16} \delta}.
$$
\item Taking $(d,p,r,q)=(0,2,\infty,1)$, and similarly to what we did to obtain \eqref{EstimStriIsoW}, we have (thanks to Propositions \ref{Propermut} and \ref{injectionLr}),
\begin{multline}
  \|\We\|_{L^2 L^\infty} \leq \|\We\|_{L^2\dot{B}_{\infty, 1}^0} \leq \|\We\|_{\tilde{L}^2 \dot{B}_{\infty, 1}^0}\\
  \leq \frac{C_{\theta}}{\nu^{\frac{1-\theta}2}}  \ee^{\frac{\theta}{2}} \left(\|\voe\|_{\dot{B}_{2,1}^{\frac12+\theta}} +\|\ce\cdot \n \ce\|_{L_t^1 \dot{B}_{2,1}^{\frac12+\theta}} +\frac1{\nu^\frac14} \|\tb\cdot \n \ce +\ce \cdot \n \tb\|_{L^\frac43 \dot{B}_{2,1}^{\theta}} \right).
\end{multline}
\end{itemize}
Now we can use the arguments from \cite{FCRF}: applying Proposition \ref{estimBsHs} with $(\aa, \bb)=(a\theta, b\theta)$, we have:
\begin{equation}
 \|\voe\|_{\dot{B}_{2,1}^{\frac12+\theta}}\leq C_{a,b} \|\voe\|_{\dot{H}^{\frac12+(1-a)\theta}}^\frac{b}{a+b} \|\voe\|_{\dot{H}^{\frac12+(1+b)\theta}}^\frac{a}{a+b}
 \leq C_{a,b} \|\voe\|_{\dot{H}^{\frac12+\mu\delta}}^\frac{b}{a+b} \|\voe\|_{\dot{H}^{\frac12+\delta}}^\frac{a}{a+b},
 \label{Stri2inf}
\end{equation}
if we have chosen $a,b>0$ and $\theta\in]0,1]$ so that:
$$
\begin{cases}
 (1-a)\theta =\mu\delta,\\
 (1+b)\theta = \delta.
\end{cases}
$$
To do this, for some $b>0$ (that will be chosen small in a few lines), we simply choose $\theta=\frac{\delta}{1+b}$ and $a=1-\mu(1+b)$ (we recall that $\mu\in]0,1[$ so $a$ is positive if $b>0$ is small enough). Thanks to the assumption \eqref{Hypinit} we obtain:
$$
  \|\voe\|_{\dot{B}_{2,1}^{\frac12+\theta}} \leq C_{b,\mu} \|\voe\|_{\dot{H}^{\frac12+\mu\delta}}^\frac{b}{(1-\mu)(1+b)} \|\voe\|_{\dot{H}^{\frac12+\delta}}^\frac{1-\mu(1+b)}{(1-\mu)(1+b)} \leq C_{b,\mu} \Co \ee^{-\gamma}.
$$
Using the previous estimates and Proposition \ref{PropFext}, we obtain that for some positive constant $C=C(\nu',\|\tu_0\|_{L^2},\delta,b,\mu)$:
$$
 \|\ce\cdot \n \ce\|_{L^1 \dot{B}_{2,1}^{\frac12+\theta}} \leq C_{b,\mu} \|\ce\cdot \n \ce\|_{L^1 \dot{H}^{\frac12+\mu\delta}}^\frac{b}{(1-\mu)(1+b)} \|\ce\cdot \n \ce\|_{L^1 \dot{H}^{\frac12+\delta}}^\frac{1-\mu(1+b)}{(1-\mu)(1+b)} \leq C \|\coe\|_{H^{\frac12+\delta}}^2.
$$
Similarly,
$$
 \|\tb\cdot \n \ce\|_{L^\frac43 \dot{B}_{2,1}^{\theta}} \leq C_{b,\mu} \|\tb\cdot \n \ce\|_{L^\frac43 \dot{H}^{\mu\delta}}^\frac{b}{(1-\mu)(1+b)} \|\tb\cdot \n \ce\|_{L^\frac43 \dot{H}^{\delta}}^\frac{1-\mu(1+b)}{(1-\mu)(1+b)} \leq C \|\tb_0\|_{H^\delta} \|\coe\|_{H^{\frac12 +\delta}},
$$
and
$$
 \|\ce\cdot \n \tb\|_{L^\frac43 \dot{B}_{2,1}^{\theta}} \leq C \|\tb_0\|_{H^\delta} \|\coe\|_{H^{\frac12 +\delta}}.
$$
Plugging these estimates into \eqref{Stri2inf}, there exists a constant $\Eo=\Eo(\nu, \nu',\Co, \delta, \eta_0, \mu)>0$ such that (the last estimate being true because $\ee\leq \ee_0$, where $\ee_0$ has been defined in the proof of the first estimates):
$$
 \|\We\|_{L^2 L^\infty} \leq C \ee^{\frac{\delta}{2(1+b)}} \left(\Co \ee^{-\gamma} +\|\coe\|_{H^{\frac12+\delta}}^2 +\|\tb_0\|_{H^\delta}^2 \right) \leq C \ee^{\frac{\delta}{2(1+b)}-\gamma} =\Eo \ee^{\frac{\delta}2 \left(2\eta_0 +\frac1{1+b}-1\right)}.
$$
As we have:
$$
\frac1{1+b}-1\geq -\eta_0\quad \Longleftrightarrow \quad b\leq \frac{\eta_0}{1-\eta_0},
$$
Finally choosing $b=\frac{\eta_0}{1-\eta_0}$, we obtain
\begin{equation}
 \|\We\|_{L^2 L^\infty} \leq \Eo \ee^{\frac{\eta_0}2 \delta}.
\end{equation}

\subsubsection{Proof of the anisotropic Strichartz estimates for System \eqref{SystWe}}
\label{Sect:AnisoStri}
We prove in this section the last two estimates from Proposition \ref{PropestimStriStrong}.

As in the previous section, writing \eqref{EstimStriAnisoUnif} for $k=\frac43$ we get that for any $d\in \R$, $m>2$, $\theta\in]0,1]$, $p\in[2, \frac4{\theta(1-\frac2{m})}[$, there exists a constant $C=C_{p,\theta,m}>0$ such that:
\begin{multline}
  \||D|^d \We\|_{L^p L_{h,v}^{m,2}}\\
  \leq \frac{C_{p,\theta,m}}{\nu^{\frac{1}{p}-\frac{\theta}4(1-\frac2{m})}} \ee^{\frac{\theta}4(1-\frac2{m})} \left(\|\voe\|_{\dot{B}_{2,1}^{\sigma_3}} +\|\ce\cdot \n \ce\|_{L^1 \dot{B}_{2,1}^{\sigma_3}} +\frac1{\nu^\frac14} \|\tb\cdot \n \ce +\ce \cdot \n \tb\|_{L_t^\frac43 \dot{B}_{2,1}^{\sigma_3-\frac12}} \right).
  \label{EstimStriAnisoW}
\end{multline}
\begin{rem}
 \sl{Notice that here $\sigma_3$ is kept as in Proposition \ref{EstimStrianiso} as we have to deal, as in the end of the previous section, with the third Besov index equal to $1$ through the methods from \cite{FCRF}.}
\end{rem}
\begin{itemize}
 \item Choosing $(d,p,m)=(0,\frac2{1-s}, \infty)$ we obtain that $\sigma_3=s+\frac{\theta}2$. Reproducing the arguments leading to \eqref{Stri2inf}, we define $a,b>0$ and $\theta\in]0,1]$ so that:
$$
\begin{cases}
s+ (1-a)\frac{\theta}2 =\frac12 +\mu\delta,\\
s+ (1+b)\frac{\theta}2 =\frac12 +\delta.
\end{cases}
$$
First consider some small $b>0$ (whose value will be precised in the following lines), and simply define (as above, the second term is positive when $b>0$ is small enough)
\begin{equation}
 \theta=\frac2{1+b}(\frac12+\delta-s)\quad \mbox{and}\quad a=1-(1+b)\frac{\frac12+\mu\delta-s}{\frac12+\delta-s}.
\label{Valeursatheta}
\end{equation}
Then, thanks to Proposition \ref{estimBsHs}, we get:
\begin{multline}
 \|\voe\|_{\dot{B}_{2,1}^{s+\frac{\theta}2}}\leq C_{a,b} \|\voe\|_{\dot{H}^{s+ (1-a)\frac{\theta}2}}^\frac{b}{a+b} \|\voe\|_{s+ (1+b)\frac{\theta}2}^\frac{a}{a+b}\\
 \leq C_{a,b} \|\voe\|_{\dot{H}^{\frac12+\mu\delta}}^\frac{b}{a+b} \|\voe\|_{\dot{H}^{\frac12+\delta}}^\frac{a}{a+b} \leq C_{b,\mu,\delta,s} \Co \ee^{-\gamma}.
 \label{Stri2inf2}
\end{multline}
Similarly, we obtain that:
$$
\begin{cases}
 \|\ce\cdot \n \ce\|_{L^1 \dot{B}_{2,1}^{s+\frac{\theta}2}} \leq C_{\nu',\delta,s,\mu,b} \|\coe\|_{H^{\frac12+\delta}}^2,\\
 \|\tb\cdot \n \ce +\ce\cdot \n \tb\|_{L^\frac43 \dot{B}_{2,1}^{s+\frac{\theta}2-\frac12}} \leq C_{\nu',\delta,s,\mu,b} \|\tb_0\|_{H^\delta} \|\coe\|_{H^{\frac12 +\delta}}.
\end{cases}
$$
Plugging these estimates and \eqref{Valeursatheta} into \eqref{EstimStriAnisoW}, thanks to the fact that $\ee\leq \ee_0$, we obtain that for all $s\in[\frac12-\frac{\eta_0}{10}\delta, \frac12+\frac{\eta_0}{10}\delta]$,
\begin{multline}
\|\We\|_{L^\frac2{1-s} L_{h,v}^{\infty, 2}} \leq C_{\nu,\nu',\Co,\delta,s,\mu,b} \ee^{\frac1{2(1+b)}(\frac12+\delta-s)} \left(\Co \ee^{-\gamma} +\|\coe\|_{H^{\frac12+\delta}}^2 +\|\tb_0\|_{H^\delta}^2 \right)\\
\leq C_{\nu,\nu',\Co,\delta,s,\mu,b} \ee^{\frac1{2(1+b)}(\frac12+\delta-s)-\gamma} \leq \Eo \ee^{\frac{\delta}2 \left(2\eta_0 +\frac1{1+b}(1-\frac{\eta_0}{10})-1\right)}.
\end{multline}
As we have
$$
2\eta_0 +\frac1{1+b}(1-\frac{\eta_0}{10})-1 \geq \frac95\eta_0 \quad \Longleftrightarrow \quad b\leq \frac{\eta_0}{10-2\eta_0},
$$
we can choose $b=\frac{\eta_0}{10-2\eta_0}$ and finally obtain that for all $s\in[\frac12-\frac{\eta_0}{10}\delta, \frac12+\frac{\eta_0}{10}\delta]$,
$$
\|\We\|_{L^\frac2{1-s} L_{h,v}^{\infty, 2}} \leq \Eo \ee^{\frac{9\eta_0}{10}\delta}.
$$
Note that the conditions on $p,\theta$ are satisfied as we already asked that $\delta\leq \frac16$.
\item For the last estimate, choosing $(d,p,m)=(1,2,\frac2{1-s})$, we obtain $\sigma_3=s+s\frac{\theta}2$. As previously, for some $b>0$ that will be precised later, choosing $a,\theta$ as follows,
$$
 \theta=\frac2{s+b}(\frac12+\delta-s)\quad \mbox{and}\quad a=s-(s+b)\frac{\frac12+\mu\delta-s}{\frac12+\delta-s},
$$
ensures that:
$$
\begin{cases}
s+ (s-a)\frac{\theta}2 =\frac12 +\mu\delta,\\
s+ (s+b)\frac{\theta}2 =\frac12 +\delta,
\end{cases}
$$
which allows us to reproduce the previous arguments and leads for all $s\in[\frac12-\frac{\eta_0}{10}\delta, \frac12+\frac{\eta_0}{10}\delta]$ to
\begin{multline}
 \|\n \We\|_{L^2 L_{h,v}^{\frac2{1-s}, 2}} \leq C_{\nu,\nu',\Co,\delta,s,\mu,b} \ee^{\frac{s}{2(s+b)}(\frac12+\delta-s)-\gamma}\\
 \leq C_{\nu,\nu',\Co,\delta,s,\mu,b} \ee^{\frac{\delta}2 \left(2\eta_0 +\frac{s}{s+b}(1-\frac{\eta_0}{10})-1\right)} \leq C_{b,c,\nu,\nu',\delta,\Co} \ee^{\frac{\delta}2 \left(2\eta_0 +\frac{\frac12-\frac{\eta_0}{10}\delta}{\frac12-\frac{\eta_0}{10}\delta+b}(1-\frac{\eta_0}{10})-1\right)}.
\end{multline}
Choosing $b=(\frac12-\frac{\eta_0}{10}) \frac{\eta_0}{10-2\eta_0}$ we finally obtain that (the conditions on $p,\theta$ being satisfied as $\delta\leq \frac16$.)
$$
\|\n \We\|_{L^2 L_{h,v}^{\frac2{1-s}, 2}} \leq \Eo \ee^{\frac{9\eta_0}{10}\delta},
$$
which ends the proof of Proposition \ref{PropestimStriStrong}. $ \blacksquare$
\end{itemize}




\begin{thebibliography}{}
	

\bibitem{AP2008} H. Abidi, M. Paicu, Global existence for the magnetohydrodynamic system in critical spaces, \textit{Proc. Roy. Soc. Edinburgh Sect. A}, \textbf{138} (2008), no.3, p.447-476.

\bibitem{AKL2021} J. Ahn, J. Kim, J. Lee, Global solutions to 3D incompressible rotational MHD system, \textit{J. Evol. Equ.}, \textbf{21} (2021), no. 1, p.235-246.

\bibitem{A2019} D. Ars\'enio, Recent progress on the mathematical theory of plasmas, \textit{Boletim da Sociedade Portuguesa de Matem\'atica}, \textbf{77} (2019), p.27–38.

\bibitem{ASR2019} D. Ars\'enio and L. Saint-Raymond, From the Vlasov–Maxwell–Boltzmann system to incompressible viscous electro-magnetohydrodynamics, Vol. 1, \textit{EMS Monographs in Mathematics, European Mathematical Society (EMS), Zürich}, 2019.

\bibitem{BMN1} A. Babin, A. Mahalov and B. Nicolaenko, Global Splitting, Integrability and Regularity of 3D Euler and Navier-Stokes Equations for Uniformly Rotating Fluids, \textit{European Journal of Mechanics}, \textbf{15}, 1996, p.291-300.



\bibitem{Dbook} H. Bahouri, J.-Y. Chemin, R. Danchin. Fourier analysis and nonlinear partial differential equations, \textit{Grundlehren der mathematischen Wissenschaften}, \textit{343}, \textit{Springer Verlag}, 2011.

\bibitem{BIM2005} J. Benameur, S. Ibrahim, M. Majdoub, Asymptotic study of a magneto-hydrodynamic system, \textit{Differential Integral Equations}, \textbf{18} (2005), no. 3, p.299-324.

\bibitem{CKS1997} R.E. Caflisch, I. Klapper, G. Steele, Remarks on singularities, dimension and energy dissipation for indeal hydrodynamics and MHD, \textit{Comm. Math. Phys.}, \textbf{184} (1997) p.443-455.

\bibitem{CW2010} C. Cao, J. Wu, Two regularity criteria for the 3D MHD equations, \textit{J. Differential Equations}, \textbf{248} (2010), no.9, p.2263-2274.

\bibitem{CW2011} Cao C., Wu J., Global regularity for the 2D MHD equations with mixed partial dissipation and magnetic diffusion, \textit{Adv. Math.}, \textbf{226} (2011), no.2, p.1803-1822.

\bibitem{CWY2014} Cao C., Wu J., B. Yuan, The 2D incompressible magnetohydrodynamics equations with only magnetic diffusion, \textit{SIAM J. Math. Anal.}, \textbf{46} (2014), no.1, p.588-602.








\bibitem{FC1} F. Charve, Convergence of weak solutions for the primitive system of the quasi-geostrophic equations, \textit{Asymptotic Analysis}, \textbf{42} (2005), p.173-209.

\bibitem{FC2} F. Charve, Global well-posedness and asymptotics for a geophysical fluid system,  \textit{Communications in Partial Differential Equations}, \textbf{29 (11 \& 12)} (2004), p.1919-1940.



\bibitem{FCVSN} F. Charve, V-S. Ngo, Asymptotics for the primitive equations with small anisotropic viscosity, \textit{Revista Matem\'atica Iberoamericana}, \textbf{27} (1) (2011), p.1-38.




\bibitem{FCestimLp} F. Charve, A priori estimates for the 3D quasi-geostrophic system, \textit{J. Math. Anal. Appl}, \textbf{444} (2016), no. 2, p.911-946.



\bibitem{FCPAA} F. Charve, Enhanced convergence rates and asymptotics for a dispersive Boussinesq-type system with large ill-prepared data, \textit{Pure and Appl. Anal.}, \textbf{2} (2020), no. 2, p.477-517.

\bibitem{FCcompl} F. Charve, Sharper dispersive estimates and asymptotics for a Boussinesq-type system with larger ill-prepared initial data, \textit{Asymptotic Analysis}, \textbf{131} (2023), no. 3-4, p. 443-470.

\bibitem{FCRF} F. Charve, Asymptotics for the rotating fluids and primitive systems with large ill-prepared initial data in critical spaces, \textit{Tunis. Journ. Math}, \textbf{Vol. 5} (2023), no. 1, p 171–213.

\bibitem{FCStratif1} F. Charve, Hidden asymptotics for the weak solutions of the strongly stratified Boussinesq system without rotation, to appear in \textit{Journal de math\'ematiques pures et appliqu\'ees}, \textbf{202} (2025), https://doi.org/10.1016/j.matpur.2025.103750.

\bibitem{FCStratif2} F. Charve, New asymptotics for strong solutions of the strongly stratified Boussinesq system without rotation and for ill-prepared initial data, \textit{Communications in Partial Differential Equations}, \textbf{49 (7 \& 8)} (2024), p.735-779.


\bibitem{Chemin2} J.-Y. Chemin, A propos d'un probl{\`e}me de p{\'e}nalisation de type antisym{\'e}trique, \textit{Journal de Math{\'e}matiques pures et appliqu{\'e}es}, \textbf{76} (1997), p.739-755.

\bibitem{CDGG} J.-Y. Chemin, B.Desjardins, I. Gallagher and E. Grenier, Anisotropy and dispersion in rotating fluids, \textit{Nonlinear Partial Differential Equations and their application, Coll\`ege de France Seminar}, Studies in Mathematics and its Applications, \textbf{31} (2002), p.171-191.

\bibitem{CDGG2} J.-Y. Chemin, B.Desjardins, I. Gallagher and E. Grenier, Fluids with anisotropic viscosity, Special issue for R. Temam's 60th birthday, \textit{M2AN. Mathematical Modelling and Numerical Analysis}, \textbf{34} (2000), no. 2, p.315-335.


\bibitem{CDGGbook} J.-Y. Chemin, B. Desjardins, I. Gallagher and E. Grenier, Mathematical Geophysics: An introduction to rotating fluids and to the Navier-Stokes equations, \textit{Oxford University Press}, 2006.

\bibitem{CW2002} G-Q. Chen, D. Wang, Global solutions of nonlinear magnetohydrodynamics with large initial data, \textit{J. Differential Equations}, \textbf{182} (2002), no. 2, p.344-376.




\bibitem{CF2021} D. Cobb, F. Fanelli, On the fast rotation asymptotics of a non-homogeneous incompressible MHD system, \textit{Nonlinearity}, \textbf{34} (2021), no. 4, p.2483-2526.






\bibitem{D2001} P. A. Davidson, An introduction to magnetohydrodynamics, \textit{Cambridge Texts in Applied Mathematics, Cambridge University Press, Cambridge}, 2001.

\bibitem{DDG1999} B. Desjardins, E. Dormy, E. Grenier, Stability of mixed Ekman-Hartmann boundary layers, \textit{Nonlinearity}, \textbf{12} (1999), no. 2, p.181-199.


\bibitem{Dutrifoy2} A. Dutrifoy, Examples of dispersive effects in non-viscous rotating fluids, \textit{Journal de Math\'ematiques Pures et Appliqu\'ees}, \textbf{84} (9) (2005), no. 3, p.331-356.

\bibitem{DL1972} In\'equations en thermo\'elasticit\'e et magn\'etohydrodynamique, \textit{Arch. Ration. Mech. Anal.}, \textbf{46} (1972), p.241-279.



\bibitem{FMRR2014} C. L. Fefferman, D. S. McCormick, J. C. Robinson, J. L. Rodrigo, Higher order commutator estimates and local existence for the non-resistive MHD equations and related models, \textit{J. Funct. Anal.}, \textbf{267} (2014), no.4, p.1035-1056.


\bibitem{IG1} I. Gallagher, Applications of Schochet's Methods to Parabolic Equation, \textit{ Journal de Math\'ematiques Pures et Appliqu\'ees}, \textbf{77} (1998), p.989-1054.




\bibitem{G1997} E. Grenier, Oscillatory perturbations of the Navier-Stokes equations, \textit{J. Math. Pures Appli.}, \textbf{76} (1997), p.477-498.


\bibitem{HX2005} C. He, Z. Xin, On the regularity of weak solutions to the magnetohydrodynamic equations, \textit{J. Differential Equations}, \textbf{213} (2005), no.2, p.235-254.

\bibitem{HW2010} X. Hu, D. Wang, Global existence and large-time behavior of solutions to the three-dimensional equations of compressible magnetohydrodynamic flows, \textit{Arch. Ration. Mech. Anal.}, \textbf{197} (2010), no. 1, p.203-238.




















\bibitem{KimJ} J. Kim, Rotational effect on the asymptotic stability of the MHD system, \textit{Journ. Diff. Eq.}, \textbf{319} (15) (2022), p.288-311.

\bibitem{KLT} Y. Koh, S. Lee and R. Takada, Dispersive estimates for the Navier-Stokes equations in the rotational framework, \textit{Advances in differential equations}, \textbf{19} (9-10) (2014), p.857-878.




\bibitem{LR2002} P. G. Lemari\'e-Rieusset, Recent developments in the Navier-Stokes problem, \textit{Boca Raton, FL: CRC Press}, 2002.

\bibitem{LZ2014} F. Lin, P. Zhang, Global small solutions to MHD type system (I): 3-D case, \textit{Comm. Pure Appl. Math.}, \textbf{67} (2014), p.531-580.





%



\bibitem{MuWei} P. Mu, Z. Wei, Rotation-dominant three-scale limit of the Cauchy problem to the inviscid rotating stratified Boussinesq equations, \textit{J. Differential Equations} \textbf{353} (2023), 385–419.


\bibitem{N2017} V.-S. Ngo, A global existence result for the anisotropic rotating magnetohydrodynamical systems, \textit{Acta Applicandae Mathematicae}, \textbf{150} (2017), p.1–42.


\bibitem{R2005} F. Rousset, Stability of large amplitude Ekman-Hartmann boundary layers in MHD: the case of ill-prepared data, \textit{Comm. Math. Phys.}, \textbf{259} (2005), n. 1, p.223–256.

\bibitem{ST1983} M. Sermange, R. Temam, Some mathematical questions related to the MHD equations, \textit{Comm. Pure Appl. Math.}, \textbf{36} (1983), p.635-664.







\bibitem{Scro3} S. Scrobogna, Global existence and convergence of nondimensionalized incompressible Navier-Stokes equations in low-Froude number regime, \textit{Discrete Contin. Dyn. Syst.} \textbf{40} (2020), no. 9, p.5471-5511.



\bibitem{Stein2} E. Stein, Harmonic analysis, real-variable methods, orthogonality and oscillatory integrals, \textit{Princeton Mathematical Series}, \textbf{43}, \textit{Princeton University Press}, 1993.


\bibitem{T2} R. Takada, Strongly statified limit for the 3D Boussinesq equations, \textit{Arch. Ration. Mech. Math.}, \textbf{232} (2019), no. 3, p.1475-1503.

\bibitem{TY} R. Takada, K. Yoneda, Global solutions for the rotating magnetohydrodynamics system in the scaling critical Sobolev space, \textit{Funkcial. Ekvac.} \textbf{67} (2024), no. 1, p.29–59.


\bibitem{Wid} K. Widmayer, Convergence to stratified flow for an inviscid 3D Boussinesq system, \textit{Commun. math. sci.} \textbf{16} (2018), no. 6, p.1713-1728.

\bibitem{W1997} J. Wu, Viscous and inviscid magneto-hydrodynamics equations, \textit{J. Analyse Math.}, \textbf{73} (1997), p.251-265.




\end{thebibliography}
\end{document}